\documentclass[reqno]{amsart}

\makeatletter
\renewcommand{\paragraph}{%
  \@startsection{paragraph}{4}{\z@}%
  {.7\baselineskip \@plus .2\baselineskip \@minus .1\baselineskip}%
  {.25\baselineskip}%
  {\normalfont\normalsize\bfseries}}
\makeatother
\usepackage{graphicx}

\makeatletter 
\newcommand*\bigcdot{\mathpalette\bigcdot@{.5}}
\newcommand*\bigcdot@[2]{\mathbin{\vcenter{\hbox{\scalebox{#2}{$\m@th#1\bullet$}}}}}
\makeatother

\usepackage{amsmath}
\usepackage{amsthm}
\usepackage{amsfonts}
\usepackage{mathrsfs}
\usepackage{xcolor}

\usepackage{enumerate}
\usepackage{amssymb}
\usepackage{bm}

\usepackage{latexsym, tikz-cd}
\usetikzlibrary{patterns}
\usepackage{placeins} 

\usepackage[colorlinks=true,linkcolor=blue]{hyperref}

\newtheorem{theorem}{Theorem}[section]
\newtheorem{proposition}[theorem]{Proposition}

\newtheorem{lemma}[theorem]{Lemma}

\theoremstyle{definition}
\newtheorem{definition}[theorem]{Definition}
\theoremstyle{remark}
\newtheorem{remark}[theorem]{Remark}

\newcommand{\nwc}{\newcommand}

\nwc{\Oph}{\operatorname{Op}_\hbar}

\renewcommand{\phi}{\varphi}

\nwc{\gl}{\langle}
\nwc{\gr}{\rangle}

\newcommand{\R}{{\mathbb R}}

\newcommand{\C}{{\mathbb C}}

\newcommand{\bd}{\partial}
\newcommand{\wind}{\operatorname{wind}}
\newcommand{\degmap}{\operatorname{deg}}

\nwc{\rest}{\restriction}

\nwc{\defeq}{\stackrel{\rm{def}}{=}}

\newcommand{\bdr}{\partial} 
\newcommand{\region}{\Omega}
\newcommand{\disk}{\mathbb{D}}

\newcommand{\conj}{\bar}

\newcommand{\restrict}[2]{\left.#1\right|_{#2}}

\newcommand{\abs}[1]{\left|#1\right|}

\newcommand{\peri}[1]{\abs{\bdr{#1}}}

        \definecolor{pink}{rgb}{1,0,1}

\makeatletter

\@addtoreset{equation}{section}
\makeatother

\title[Third Robin Eigenvalue on Curved Surfaces]{On the Isoperimetric Inequality for the Third Robin Eigenvalue on Surfaces 
}
\date{}

\author{Xiaoliang Wang, Zuoqin Wang}
\thanks{Partially supported by  NNSFC no. 12571064.}
\address{School of Mathematical Sciences\\
	University of Science and Technology of China\\
	Hefei, 230026\\ P.R. China}
\email{wangxl514@mail.ustc.edu.cn}

\address{School of Mathematical Sciences\\
	University of Science and Technology of China\\
	Hefei, 230026\\ P.R. China}
\email{wangzuoq@ustc.edu.cn}

\begin{document}

\begin{abstract}
We investigate the third Robin eigenvalue on simply connected domains
of area $A$ in surfaces with Gaussian curvature bounded above by $K$.
For $K>0$, we discover a transition at the critical parameter
$\alpha=-4\pi+KA$: two equal geodesic disks give a sharp upper bound
when $\alpha=-4\pi$ or $-4\pi+KA\leq\alpha\leq4\pi$, but the two disk bound fails for $-4\pi<\alpha<-4\pi+KA$.
In the latter range, we establish a new upper bound involving a single
geodesic disk with a modified Robin parameter and prove its sharpness
in an explicit subrange. The sharp two-disk inequality persists for
nonpositive curvature bounds, and Steklov--Helmholtz analogues are also obtained. 
Both the two-disk and one-disk bounds arise from a unified two-pole Green-function method, with the one-disk regime emerging from a degenerate area split that produces the $2\pi$ shift in the effective Robin parameter, while sharpness in the stated subrange is established by boundary oscillations that redistribute Robin boundary mass.
\end{abstract}

\maketitle

\section{Introduction}

Let $(M,g)$ be a complete Riemannian surface and let $\region\subset M$ be a
bounded simply connected domain with Lipschitz boundary. Let $\nu$
denote the outward unit normal to $\partial\region$. For a real
parameter $\alpha$, consider the Robin eigenvalue problem, interpreted
in the weak sense,
\begin{equation}\label{eq:intro-Robin-problem}
	\begin{cases}
		-\Delta_g u=\lambda u & \text{in }\region,\\[2mm]
		\displaystyle\partial_\nu u+\frac{\alpha}{L}u=0
		& \text{on }\partial\region,
	\end{cases}
	\qquad L:=|\partial\region|_g.
\end{equation}
We enumerate the eigenvalues, with multiplicities, as
\[
\lambda_1\left(\region;\frac{\alpha}{L}\right)
\leq\lambda_2\left(\region;\frac{\alpha}{L}\right)
\leq\lambda_3\left(\region;\frac{\alpha}{L}\right)
\leq\cdots.
\]
The normalization by the boundary length makes $\alpha$ dimensionless
under rescaling of the metric. In particular, $\alpha=0$ is the Neumann case.
We refer to Bucur, Freitas and Kennedy~\cite{bucurfreitaskennedy2017}
for background on Robin spectral optimization.

The classical starting point for maximization of a positive Neumann
eigenvalue is the Szeg\"o--Weinberger inequality. Szeg\"o proved
that the first positive Neumann eigenvalue of a simply connected planar
domain of prescribed area is maximized by a disk~\cite{szego1954};
Weinberger removed the topological restriction and extended the result
to higher dimensions~\cite{weinberger1956}. Bandle established
curvature-dependent disk comparisons on surfaces by conformal
methods~\cite{bandle1972}. For the second Robin eigenvalue, Freitas and
Laugesen proved the corresponding sharp planar disk inequality for
$-2\pi\leq\alpha\leq2\pi$~\cite{freitaslaugesen2020}; see
also~\cite{freitaslaugesen2021}. Langford and Laugesen subsequently
obtained comparisons on spherical and hyperbolic domains, and on
surfaces with an upper Gaussian-curvature bound, in appropriate
area--parameter ranges~\cite{langfordlaugesen2025robin}.

For the \emph{third}  eigenvalue, which is the main scope of this paper, 
the natural candidate is a disjoint union of two equal disks, to be
approached by connected domains with a narrow neck. In the planar
Neumann problem this phenomenon was established by Girouard,
Nadirashvili and Polterovich~\cite{girouardnadirashvilipolterovich2009}.
Girouard and Polterovich extended the result to log-subharmonic
densities~\cite{girouardpolterovich2010shape}, 
while Bucur and Henrot
proved the Euclidean inequality in all dimensions without a topological
assumption~\cite{bucurhenrot2019}. The corresponding constant-curvature
comparisons were established in hyperbolic space by Freitas and
Laugesen~\cite{freitaslaugesen2022} and on the sphere by Bucur,
Martinet and Nahon~\cite{bucurmartinetnahon2025}.
For Robin boundary conditions, Girouard and Laugesen obtained the
sharp planar two-disk inequality for $-4\pi\leq\alpha\leq0$
~\cite{girouardlaugesen2021};  Kim and Laugesen extended it to the
entire interval $[-4\pi,4\pi]$, with a degree-theoretic orthogonality
argument accounting for the positive-parameter range
~\cite{kimlaugesen2025}. Our main result reveals a new phenomenon: in the positive-curvature case, the Robin parameter determines whether the sharp comparison is provided by two equal disks or by a different one-disk upper bound.

The variable-curvature counterpart is more subtle. We assume only
\begin{equation}\label{eq:intro-curvature-bound}
	K_g\leq K\qquad\text{on }\region,
\end{equation}
where $K\in\mathbb R$ is fixed. The sharp two-disk inequality under this curvature hypothesis was
conjectured by Langford and Laugesen
~\cite[Conjecture~1.4]{langfordlaugesen2023neumann}; see also the
geometric spectral theory workshop report
~\cite[p.~2099]{hassannezhad2023geometric}. The conjecture was recently verified by Liu, Long and Zou 
~\cite{liulongzou2026} through a stronger reciprocal inequality, using two-pole Green coordinates. 
Chen and Yang also obtained related harmonic-mean Neumann and
Aharonov--Bohm inequalities through reciprocal variational arguments
and Green-level comparison
~\cite{chenyang2026}.
Thus the Neumann endpoint of our first inequality overlaps with the
independent result of Liu, Long and Zou; the present paper concerns the
Robin family for both signs of $\alpha$ and the different optimal
comparisons arising in positive curvature.

The purpose of this paper is twofold. First, writing \(A=|\Omega|_g\), we establish upper bounds for the third Robin eigenvalue \(\lambda_3\bigl(\Omega;\alpha/|\partial\Omega|_g\bigr)\) throughout \(-4\pi\leq\alpha\leq4\pi\), under the curvature assumption (1.2) and the area restriction \(KA<4\pi\) when \(K>0\). A sharp comparison with two equal geodesic disks holds when \(\alpha=-4\pi\) or \(-4\pi+KA\leq\alpha\leq4\pi\). In the complementary range \(-4\pi<\alpha<-4\pi+KA\), which occurs only in positive curvature, the two-disk bound fails. We instead obtain an upper bound given by the first angular Robin eigenvalue of a single geodesic disk of area \(A\), with normalized Robin parameter \(\alpha+2\pi\). This one-disk bound is sharp when \(\alpha\leq-2\pi\) and the first angular eigenvalue does not exceed the second radial eigenvalue for the same Robin parameter. Second, we transfer these results to the Steklov–Helmholtz problem through a version of Robin–Steklov duality that incorporates the Dirichlet index shift.

\subsection{Main results}\label{sec:intro-main-results}
Write $\mathbb M_K^2$ for the simply connected space form of constant
curvature $K$. For $m>0$ (with $Km<4\pi$ if $K>0$), let $D_K(m)$
be its geodesic disk of area $m$, and set \( \ell_K(m):=|\partial D_K(m)|\). We denote by
$\lambda_{1,\mathrm{ang}}(D_K(m);\gamma)$ the first angular Robin eigenvalue with coefficient $\gamma$, whose eigenfunction has the form $f(r)e^{\pm i\theta}$ in geodesic polar coordinates. We also write $\lambda_{2,\mathrm{rad}}(D_K(m);\gamma)$ for the second radial Robin eigenvalue.

For $A>0$, with $KA<4\pi$ if $K>0$, define the following sets of
curvature, area, and Robin parameters:
\begin{align*}
	\mathcal P_1
	&:=\{(K,A,\alpha):K\leq0,\ -4\pi\leq\alpha\leq4\pi\}\\
	&\quad\cup\{(K,A,\alpha):K>0,\ \alpha=-4\pi
	\ \text{or}\ -4\pi+KA\leq\alpha\leq4\pi\},\\
	\mathcal P_2
	&:=\{(K,A,\alpha):K>0,\ -4\pi\leq\alpha\leq-4\pi+KA\}.
\end{align*}
 
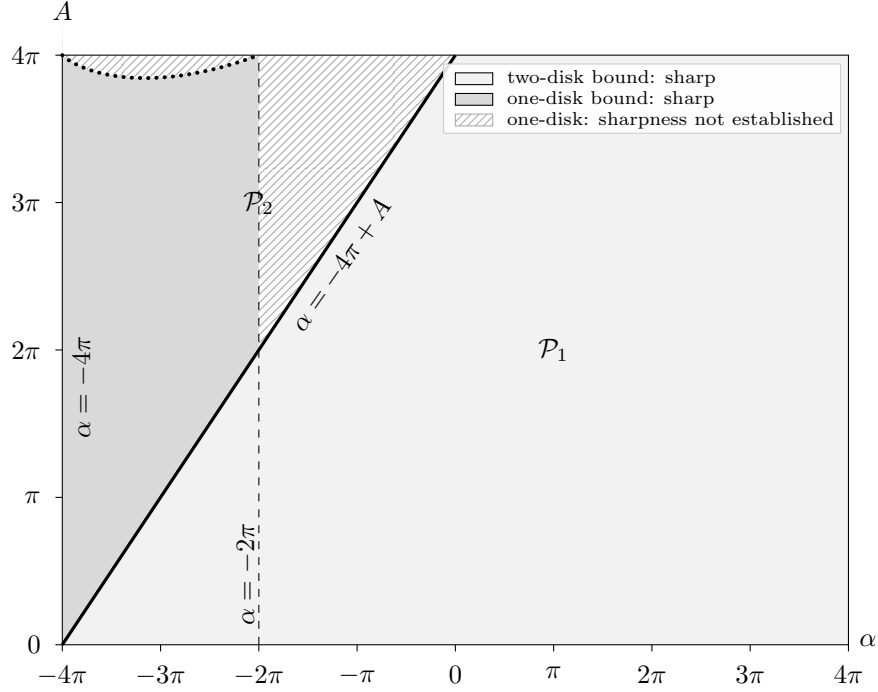
\begin{figure}[htbp]
	\centering
	\begin{tikzpicture}[x=1.30cm,y=1.95cm]
		\fill[gray!10] (-4,0) rectangle (4,4);
		
		\fill[gray!30]
		(-4,0) -- (-2,2) -- (-2,4)
		plot[smooth] coordinates {
			(-2.0000,4.0000) (-2.0978,3.9765) (-2.2175,3.9500)
			(-2.3357,3.9266) (-2.4524,3.9060) (-2.5674,3.8885)
			(-2.6807,3.8738) (-2.7922,3.8620) (-2.9021,3.8533)
			(-3.0113,3.8474) (-3.1363,3.8445) (-3.1814,3.8445)
			(-3.3009,3.8474) (-3.3860,3.8523) (-3.4746,3.8601)
			(-3.5549,3.8699) (-3.6300,3.8816) (-3.7008,3.8953)
			(-3.7677,3.9109) (-3.8274,3.9275) (-3.8810,3.9452)
			(-3.9291,3.9637) (-3.9694,3.9823) (-3.9924,3.9950)
			(-4.0000,4.0000) (-4,0)
		} -- cycle;
		
		\fill[white]
		(-4,4)
		plot[smooth] coordinates {
			(-4.0000,4.0000) (-3.9924,3.9950) (-3.9694,3.9823)
			(-3.9291,3.9637) (-3.8810,3.9452) (-3.8274,3.9275)
			(-3.7677,3.9109) (-3.7008,3.8953) (-3.6300,3.8816)
			(-3.5549,3.8699) (-3.4746,3.8601) (-3.3860,3.8523)
			(-3.3009,3.8474) (-3.1814,3.8445) (-3.1363,3.8445)
			(-3.0113,3.8474) (-2.9021,3.8533) (-2.7922,3.8620)
			(-2.6807,3.8738) (-2.5674,3.8885) (-2.4524,3.9060)
			(-2.3357,3.9266) (-2.2175,3.9500) (-2.0978,3.9765)
			(-2.0000,4.0000)
		}
		-- (-2,2) -- (0,4) -- cycle;
		\fill[pattern=north east lines,pattern color=gray!65]
		(-4,4)
		plot[smooth] coordinates {
			(-4.0000,4.0000) (-3.9924,3.9950) (-3.9694,3.9823)
			(-3.9291,3.9637) (-3.8810,3.9452) (-3.8274,3.9275)
			(-3.7677,3.9109) (-3.7008,3.8953) (-3.6300,3.8816)
			(-3.5549,3.8699) (-3.4746,3.8601) (-3.3860,3.8523)
			(-3.3009,3.8474) (-3.1814,3.8445) (-3.1363,3.8445)
			(-3.0113,3.8474) (-2.9021,3.8533) (-2.7922,3.8620)
			(-2.6807,3.8738) (-2.5674,3.8885) (-2.4524,3.9060)
			(-2.3357,3.9266) (-2.2175,3.9500) (-2.0978,3.9765)
			(-2.0000,4.0000)
		}
		-- (-2,2) -- (0,4) -- cycle;
		
		\draw[very thick] (-4,0) -- (0,4);
		\draw[dashed] (-2,0) -- (-2,4);
	\draw[
	black,
	line width=1.5pt,
	dash pattern=on 0pt off 3pt,
	line cap=round
	]
		plot[smooth] coordinates {
			(-4.0000,4.0000) (-3.9924,3.9950) (-3.9694,3.9823)
			(-3.9291,3.9637) (-3.8810,3.9452) (-3.8274,3.9275)
			(-3.7677,3.9109) (-3.7008,3.8953) (-3.6300,3.8816)
			(-3.5549,3.8699) (-3.4746,3.8601) (-3.3860,3.8523)
			(-3.3009,3.8474) (-3.1814,3.8445) (-3.1363,3.8445)
			(-3.0113,3.8474) (-2.9021,3.8533) (-2.7922,3.8620)
			(-2.6807,3.8738) (-2.5674,3.8885) (-2.4524,3.9060)
			(-2.3357,3.9266) (-2.2175,3.9500) (-2.0978,3.9765)
			(-2.0000,4.0000)
		};
		
		\draw (-4,0) rectangle (4,4);
		
		\foreach \x/\lab in {
			-4/{-4\pi},-3/{-3\pi},-2/{-2\pi},-1/{-\pi},0/{0},
			1/{\pi},2/{2\pi},3/{3\pi},4/{4\pi}}
		\draw (\x,0) -- ++(0,-0.045) node[below=2pt] {$\lab$};
		\foreach \y/\lab in {0/{0},1/{\pi},2/{2\pi},3/{3\pi},4/{4\pi}}
		\draw (-4,\y) -- ++(-0.07,0) node[left=2pt] {$\lab$};
		\node[below=15pt] at (4.2,0.4) {$\alpha$};
		\node at (-4,4.3) {$A$};

        \node at (1,2) {$\mathcal P_1$};
        \node at (-2,3) {$\mathcal P_2$};
		
		\node[rotate=90,anchor=south] at (-3.60,1.75)
		{$\alpha=-4\pi$};
		\node[rotate=90,anchor=south] at (-1.93,0.48) {$\alpha=-2\pi$};
		\node[rotate=55] at (-1.16,2.58) {$\alpha=-4\pi+A$};
		
		\fill[white,opacity=.88] (-0.12,3.94) rectangle (3.92,3.47);
		\draw[gray!40] (-0.12,3.94) rectangle (3.92,3.47);
		\fill[gray!10] (0.00,3.80) rectangle (0.38,3.88);
		\draw (0.00,3.80) rectangle (0.38,3.88);
		\node[anchor=west,font=\scriptsize] at (0.43,3.84)
		{two-disk bound: sharp};
		\fill[gray!30] (0.00,3.66) rectangle (0.38,3.74);
		\draw (0.00,3.66) rectangle (0.38,3.74);
		\node[anchor=west,font=\scriptsize] at (0.43,3.70)
		{one-disk bound: sharp};
		\fill[white] (0.00,3.52) rectangle (0.38,3.60);
		\fill[pattern=north east lines,pattern color=gray!65]
		(0.00,3.52) rectangle (0.38,3.60);
		\draw[gray!65] (0.00,3.52) rectangle (0.38,3.60);
		\node[anchor=west,font=\scriptsize] at (0.43,3.56)
		{one-disk: sharpness not established};
	\end{tikzpicture}
	\caption{The positive-curvature parameter ranges after scaling to $K=1$.
The horizontal axis is $\alpha$; the vertical axis is the rescaled area
$A$ (equal to $KA$ before scaling). Both axes are marked in multiples
of $\pi$, and the upper edge $A=4\pi$ is excluded.
The diagonal is $\alpha=-4\pi+A$. On and to its right, and also on
the line $\alpha=-4\pi$, the equal two-disk bound
\eqref{eq:intro-Robin-main} is sharp.
In the wedge $-4\pi<\alpha<-4\pi+A$, the one-disk bound
\eqref{eq:intro-Robin-one-disk} holds. Dark gray indicates the subrange
where sharpness is proved under \eqref{eq:intro-Robin-one-disk-sharp};
hatching indicates where sharpness has not been established.
The dotted curve schematically depicts the angular--radial crossing
defined by \eqref{eq:intro-angular-radial-crossing}.}
\label{fig:intro-Robin-parameter-ranges}
\end{figure}
\FloatBarrier

\begin{theorem}\label{thm:intro-Robin-main}
	Let $\region$ be a bounded simply connected Lipschitz domain in a
	Riemannian surface $(M,g)$, and put
	$A=|\region|_g$, $L=|\partial\region|_g$.
	Assume $K_g\leq K$ on $\region$, with $KA<4\pi$ when $K>0$.
	
	If $(K,A,\alpha)\in\mathcal P_1$, then
	\begin{equation}\label{eq:intro-Robin-main}
		\lambda_3\left(\region;\frac{\alpha}{L}\right)
		<
		\lambda_2\left(D_K(A/2);\frac{\alpha}{2\ell_K(A/2)}\right).
	\end{equation}
	This inequality is sharp among bounded simply connected Lipschitz
	domains of area $A$ in $\mathbb M_K^2$.
	
	If $(K,A,\alpha)\in\mathcal P_2$, then
	\begin{equation}\label{eq:intro-Robin-one-disk}
		\lambda_3\left(\region;\frac{\alpha}{L}\right)
		<
		\lambda_{1,\mathrm{ang}}\left(
		D_K(A);\frac{\alpha+2\pi}{\ell_K(A)}\right).
	\end{equation}
	The latter bound is sharp whenever
	\begin{equation}\label{eq:intro-Robin-one-disk-sharp}
		\alpha\leq-2\pi,
		\qquad
		\lambda_{1,\mathrm{ang}}\left(
		D_K(A);\frac{\alpha+2\pi}{\ell_K(A)}\right)
		\leq
		\lambda_{2,\mathrm{rad}}\left(
		D_K(A);\frac{\alpha+2\pi}{\ell_K(A)}\right).
	\end{equation}
\end{theorem}

Figure~\ref{fig:intro-Robin-parameter-ranges} displays the
positive-curvature parameter ranges after scaling to $K=1$.

The right-hand side of \eqref{eq:intro-Robin-main} is the third
Robin eigenvalue of two equal disjoint geodesic disks, with the Robin
coefficient normalized by their total boundary length. At $\alpha=0$
this gives the sharp Neumann inequality
\[
\lambda_3^N(\region)<\lambda_2^N\bigl(D_K(A/2)\bigr).
\]
When $K>0$, the two comparison values in Theorem~\ref{thm:intro-Robin-main}
also agree at $\alpha=-4\pi$ and $\alpha=-4\pi+KA$, where they equal
$0$ and $2K$, respectively.

\begin{figure}[!ht]
\centering
\begin{tikzpicture}[line width=0.55pt]
  \begin{scope}[xshift=0cm]
    \fill[black!3] (-0.66,0) circle (0.72);
    \fill[black!3] ( 0.66,0) circle (0.72);
    \begin{scope}[even odd rule]
      \clip (-1.60,-1.05) rectangle (1.60,1.05)
        (0.66,0) circle (0.72);
      \draw (-0.66,0) circle (0.72);
    \end{scope}
    \begin{scope}[even odd rule]
      \clip (-1.60,-1.05) rectangle (1.60,1.05)
        (-0.66,0) circle (0.72);
      \draw (0.66,0) circle (0.72);
    \end{scope}
    \node[font=\scriptsize] at (-0.88,-0.98) {$A/2$};
    \node[font=\scriptsize] at ( 0.88,-0.98) {$A/2$};
    \node[font=\small,anchor=north] at (0,-1.42) {(a) two equal disks};
  \end{scope}

  \begin{scope}[xshift=3.45cm,xscale=-1]
    \fill[black!3]
      plot[smooth cycle,samples=240,domain=0:360]
      ({-0.66+(0.68+0.030*sin(14*\x))*cos(\x)},
       {(0.68+0.030*sin(14*\x))*sin(\x)});
    \fill[black!3] (0.66,0) circle (0.73);
    \begin{scope}[even odd rule]
      \clip (-1.75,-1.05) rectangle (1.75,1.05)
        (0.66,0) circle (0.73);
      \draw
        plot[smooth cycle,samples=240,domain=0:360]
        ({-0.66+(0.68+0.030*sin(14*\x))*cos(\x)},
         {(0.68+0.030*sin(14*\x))*sin(\x)});
    \end{scope}
    \begin{scope}[even odd rule]
      \clip (-1.75,-1.05) rectangle (1.75,1.05)
        plot[smooth cycle,samples=240,domain=0:360]
        ({-0.66+(0.68+0.030*sin(14*\x))*cos(\x)},
         {(0.68+0.030*sin(14*\x))*sin(\x)});
      \draw (0.66,0) circle (0.73);
    \end{scope}
    \node[font=\scriptsize,align=center] at (-0.88,-0.98)
      {$A/2-\delta$\\oscillated};
    \node[font=\scriptsize,align=center] at (0.90,-1.04)
      {$A/2+\delta$\\smooth};
    \node[font=\small,align=center,anchor=north] at (0,-1.42)
      {(b) two unequal disks\\{\scriptsize with an oscillating smaller disk}};
  \end{scope}

  \begin{scope}[xshift=7.55cm]
    \fill[black!3] (-0.86,0) circle (0.90);
    \fill[black!3]
      plot[smooth cycle,samples=220,domain=0:360]
      ({0.30+(0.34+0.035*sin(14*\x))*cos(\x)},
       {(0.34+0.035*sin(14*\x))*sin(\x)});
    \begin{scope}[even odd rule]
      \clip (-1.95,-1.20) rectangle (1.10,1.20)
        plot[smooth cycle,samples=220,domain=0:360]
        ({0.30+(0.34+0.035*sin(14*\x))*cos(\x)},
         {(0.34+0.035*sin(14*\x))*sin(\x)});
      \draw (-0.86,0) circle (0.90);
    \end{scope}
    \begin{scope}[even odd rule]
      \clip (-1.95,-1.20) rectangle (1.10,1.20)
        (-0.86,0) circle (0.90);
      \draw
        plot[smooth cycle,samples=220,domain=0:360]
        ({0.30+(0.34+0.035*sin(14*\x))*cos(\x)},
         {(0.34+0.035*sin(14*\x))*sin(\x)});
    \end{scope}
    \node[font=\scriptsize] at (-1.02,-1.17) {$A-m$};
    \node[font=\scriptsize] at ( 0.34,-0.62) {$m\downarrow0$};
    \node[font=\scriptsize,align=center] at (0.52,0.78)
      {oscillated\\small disk};
    \node[font=\small,anchor=north] at (-0.20,-1.42) {(c) shrinking smaller disk};
  \end{scope}
\end{tikzpicture}
\caption{Schematics for (a) sharpness of the two-disk bound
\eqref{eq:intro-Robin-main}; (b) sharpness at $\alpha=-4\pi+KA$
and counterexamples to \eqref{eq:intro-Robin-main} for
$-4\pi<\alpha<-4\pi+KA$ ($K>0$); and (c) sharpness of the
one-disk bound \eqref{eq:intro-Robin-one-disk} under
\eqref{eq:intro-Robin-one-disk-sharp}. Oscillations are exaggerated
and corners are left unsmoothed for clarity.}
\label{fig:sharpness-geometries}
\end{figure}
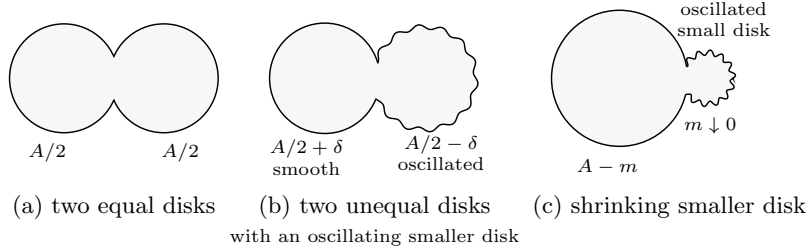
\FloatBarrier

The sharpness mechanism also changes as $\alpha$ crosses $-4\pi+KA$. For $K>0$ and $-4\pi+KA\leq\alpha\leq4\pi$, the bound
\eqref{eq:intro-Robin-main} is approached by two equal-area geodesic disks joined by a short neck; see Figure~\ref{fig:sharpness-geometries}(a). At $\alpha=-4\pi+KA$, this extremal configuration becomes unstable: for any fixed $\delta>0$ such that $0<A/2-\delta<A/2+\delta<2\pi/K$, the same bound is approached by a large disk $D_K(A/2+\delta)$ joined to a small disk $D_K(A/2-\delta)$ with an oscillating boundary; see Figure~\ref{fig:sharpness-geometries}(b).
The limiting configuration consists of the two disjoint disks $D_K(A/2+\delta)$ and $D_K(A/2-\delta)$, with the boundary measure of the smaller disk multiplied by a constant $c>1$. When $-4\pi<\alpha<-4\pi+KA$, similar unequal-disk constructions instead give counterexamples to \eqref{eq:intro-Robin-main}, and \eqref{eq:intro-Robin-one-disk} provides the new upper bound.

In the range $-4\pi < \alpha < -4\pi + KA$, under \eqref{eq:intro-Robin-one-disk-sharp},
the sharpness sequence again consists of a large disk joined by a
short neck to a small disk with an oscillating boundary, but now
their areas are $A-m$ and $m$, with $m\downarrow0$; see Figure~\ref{fig:sharpness-geometries}(c).
For $\alpha<-2\pi$, the small disk shrinks while its boundary
length tends to a positive finite limit. The domains converge to
$D_K(A)$, whereas their boundary measures converge weakly, as
measures on the ambient surface, to
\[
ds_{\partial\Omega_j}\rightharpoonup
ds_{\partial D_K(A)}
+\frac{2\pi\ell_K(A)}{-\alpha-2\pi}\,\delta_p,
\qquad p\in\partial D_K(A).
\]
Thus the disappearing small disk leaves a positive boundary mass at
the attachment point. At $\alpha=-2\pi$, the small disk still
shrinks, but its boundary length, and hence the total boundary
length, diverges. There is then no finite unnormalized boundary
measure limit; instead,
$ds_{\partial\Omega_j}/|\partial\Omega_j|\rightharpoonup\delta_p$.
These constructions are established in Section~\ref{sec:sharpness}.

For $K>0$, the switch between the two extremal geometries at
$\alpha=-4\pi+KA$ is driven by the reversal of convexity of the
model Riccati function $R_{K,\mu}$ as $\mu$ crosses $2K$. For the definition of $R_{K,\mu}$, see \eqref{eq:model-R-definition}, equivalently,
\[
R_{K,\mu}(x)
=\ell_K(x)\,\sigma_{1,\mathrm{ang}}\bigl(D_K(x),\mu\bigr),
\]
where $\sigma_{1,\mathrm{ang}}(D_K(x),\lambda)$ denotes the eigenvalue
$\sigma$ of \eqref{eq:intro-SH-problem} on $D_K(x)$ corresponding to
an eigenfunction of the form $u(r,\theta)=f(r)e^{i\theta}$,
where $f(r)>0$ for $r>0$; see \eqref{eq:model-angular-SH}.

\begin{figure}[htbp]
\centering
\begin{tikzpicture}[
  x=8.5cm,y=1.28cm,
  model curve/.style={line width=.85pt,line cap=round,line join=round},
  every node/.style={font=\small}]
  \draw[gray!55,densely dotted] (1,1.18) -- (1,-3.18);
  \draw[->,line width=.55pt] (0,-3.23) -- (0,1.43)
    node[above] {$R_{K,\mu}(x)$};
  \draw[->,line width=.55pt] (-.025,-1) -- (1.10,-1)
    node[right] {$x$};
  \draw (-.009,1) -- (.009,1);
  \node[left] at (-.012,1) {$2\pi$};
  \node[left] at (-.012,-1) {$-2\pi$};
  \node[below right,inner sep=2pt] at (0,-1) {$0$};
  \draw (1,-.96) -- (1,-1.04);
  \node[below,fill=white,inner sep=2pt] at (1,-1.10) {$4\pi/K$};

  \draw[model curve] (0,1) -- (1,1);
  \node[above] at (.46,1.04) {$\mu=0$};
  \draw[model curve] (0,1) -- (1,-1)
    node[pos=.64,below,sloped,inner sep=3pt] {$\mu=2K$};
  \draw[model curve] (0,1)
    .. controls (0.006666667,0.990000000) and (0.013333333,0.980016661) .. (0.020000000,0.970050891)
    .. controls (0.026666667,0.960085120) and (0.033333333,0.950136918) .. (0.040000000,0.940207257)
    .. controls (0.046666667,0.930277597) and (0.053333333,0.920366479) .. (0.060000000,0.910474952)
    .. controls (0.066666667,0.900583424) and (0.073333333,0.890711488) .. (0.080000000,0.880860273)
    .. controls (0.086666667,0.871009059) and (0.093333333,0.861178565) .. (0.100000000,0.851370013)
    .. controls (0.106666667,0.841561461) and (0.113333333,0.831774851) .. (0.120000000,0.822011501)
    .. controls (0.126666667,0.812248152) and (0.133333333,0.802508063) .. (0.140000000,0.792792664)
    .. controls (0.146666667,0.783077265) and (0.153333333,0.773386556) .. (0.160000000,0.763722086)
    .. controls (0.166666667,0.754057616) and (0.173333333,0.744419385) .. (0.180000000,0.734809078)
    .. controls (0.186666667,0.725198771) and (0.193333333,0.715616388) .. (0.200000000,0.706063762)
    .. controls (0.206666667,0.696511136) and (0.213333333,0.686988268) .. (0.220000000,0.677497156)
    .. controls (0.226666667,0.668006045) and (0.233333333,0.658546691) .. (0.240000000,0.649121281)
    .. controls (0.246666667,0.639695870) and (0.253333333,0.630304404) .. (0.260000000,0.620949275)
    .. controls (0.266666667,0.611594147) and (0.273333333,0.602275356) .. (0.280000000,0.592995533)
    .. controls (0.286666667,0.583715709) and (0.293333333,0.574474853) .. (0.300000000,0.565275856)
    .. controls (0.306666667,0.556076859) and (0.313333333,0.546919721) .. (0.320000000,0.537807633)
    .. controls (0.326666667,0.528695546) and (0.333333333,0.519628508) .. (0.340000000,0.510610049)
    .. controls (0.346666667,0.501591590) and (0.353333333,0.492621710) .. (0.360000000,0.483704321)
    .. controls (0.366666667,0.474786932) and (0.373333333,0.465922035) .. (0.380000000,0.457113981)
    .. controls (0.386666667,0.448305927) and (0.393333333,0.439554717) .. (0.400000000,0.430865204)
    .. controls (0.406666667,0.422175692) and (0.413333333,0.413547878) .. (0.420000000,0.404987195)
    .. controls (0.426666667,0.396426512) and (0.433333333,0.387932960) .. (0.440000000,0.379512642)
    .. controls (0.446666667,0.371092323) and (0.453333333,0.362745238) .. (0.460000000,0.354478263)
    .. controls (0.466666667,0.346211288) and (0.473333333,0.338024424) .. (0.480000000,0.329925453)
    .. controls (0.486666667,0.321826482) and (0.493333333,0.313815405) .. (0.500000000,0.305901063)
    .. controls (0.506666667,0.297986721) and (0.513333333,0.290169115) .. (0.520000000,0.282458339)
    .. controls (0.526666667,0.274747563) and (0.533333333,0.267143616) .. (0.540000000,0.259658074)
    .. controls (0.546666667,0.252172532) and (0.553333333,0.244805396) .. (0.560000000,0.237570008)
    .. controls (0.566666667,0.230334621) and (0.573333333,0.223230983) .. (0.580000000,0.216274564)
    .. controls (0.586666667,0.209318145) and (0.593333333,0.202508945) .. (0.600000000,0.195865006)
    .. controls (0.606666667,0.189221067) and (0.613333333,0.182742389) .. (0.620000000,0.176450152)
    .. controls (0.626666667,0.170157916) and (0.633333333,0.164052120) .. (0.640000000,0.158157810)
    .. controls (0.646666667,0.152263500) and (0.653333333,0.146580676) .. (0.660000000,0.141139183)
    .. controls (0.666666667,0.135697690) and (0.673333333,0.130497528) .. (0.680000000,0.125574576)
    .. controls (0.686666667,0.120651624) and (0.693333333,0.116005881) .. (0.700000000,0.111680890)
    .. controls (0.706666667,0.107355900) and (0.713333333,0.103351662) .. (0.720000000,0.099721590)
    .. controls (0.726666667,0.096091519) and (0.733333333,0.092835613) .. (0.740000000,0.090020180)
    .. controls (0.746666667,0.087204747) and (0.753333333,0.084829787) .. (0.760000000,0.082978731)
    .. controls (0.766666667,0.081127675) and (0.773333333,0.079800522) .. (0.780000000,0.079103875)
    .. controls (0.786666667,0.078407228) and (0.793333333,0.078341087) .. (0.800000000,0.079044073)
    .. controls (0.806666667,0.079747059) and (0.813333333,0.081219173) .. (0.820000000,0.083644460)
    .. controls (0.826666667,0.086069747) and (0.833333333,0.089448207) .. (0.840000000,0.094029998)
    .. controls (0.846666667,0.098611790) and (0.853333333,0.104396913) .. (0.860000000,0.111736120)
    .. controls (0.866666667,0.119075328) and (0.873333333,0.127968621) .. (0.880000000,0.138922944)
    .. controls (0.886666667,0.149877268) and (0.893333333,0.162892623) .. (0.900000000,0.178745708)
    .. controls (0.906666667,0.194598794) and (0.913333333,0.213289609) .. (0.920000000,0.236040441)
    .. controls (0.926666667,0.258791272) and (0.933333333,0.285602119) .. (0.940000000,0.318715281)
    .. controls (0.943587947,0.336536520) and (0.947175893,0.356183228) .. (0.950763840,0.378257084)
    .. controls (0.953708118,0.396370948) and (0.956652397,0.416119227) .. (0.959596675,0.437973481)
    .. controls (0.962012758,0.455907141) and (0.964428841,0.475258941) .. (0.966844924,0.496409588)
    .. controls (0.968827568,0.513765866) and (0.970810212,0.532333467) .. (0.972792856,0.552410953)
    .. controls (0.974419819,0.568886592) and (0.976046782,0.586378970) .. (0.977673745,0.605119749)
    .. controls (0.979008835,0.620498483) and (0.980343925,0.636717875) .. (0.981679016,0.653955698)
    .. controls (0.982774594,0.668101102) and (0.983870173,0.682932306) .. (0.984965751,0.698584823)
    .. controls (0.985864786,0.711429321) and (0.986763821,0.724826883) .. (0.987662855,0.738880441)
    .. controls (0.988400606,0.750412828) and (0.989138356,0.762386957) .. (0.989876106,0.774880978)
    .. controls (0.990481506,0.785133605) and (0.991086906,0.795736321) .. (0.991692306,0.806748594)
    .. controls (0.992189098,0.815785294) and (0.992685891,0.825097785) .. (0.993182684,0.834731508)
    .. controls (0.993590353,0.842636966) and (0.993998023,0.850758738) .. (0.994405692,0.859131750)
    .. controls (0.994740227,0.866002665) and (0.995074761,0.873042763) .. (0.995409296,0.880279063)
    .. controls (0.995683816,0.886217190) and (0.995958336,0.892287439) .. (0.996232856,0.898510862)
    .. controls (0.996458128,0.903617819) and (0.996683400,0.908827923) .. (0.996908672,0.914157690)
    .. controls (0.997093530,0.918531311) and (0.997278389,0.922985513) .. (0.997463248,0.927533339)
    .. controls (0.997614943,0.931265297) and (0.997766639,0.935060301) .. (0.997918334,0.938928713)
    .. controls (0.998042816,0.942103143) and (0.998167298,0.945327005) .. (0.998291780,0.948608575)
    .. controls (0.998393930,0.951301440) and (0.998496080,0.954033166) .. (0.998598230,0.956810391)
    .. controls (0.998682054,0.959089390) and (0.998765879,0.961399027) .. (0.998849704,0.963744652)
    .. controls (0.998918490,0.965669478) and (0.998987277,0.967618536) .. (0.999056064,0.969596150)
    .. controls (0.999112510,0.971218985) and (0.999168957,0.972861048) .. (0.999225403,0.974525841)
    .. controls (0.999271724,0.975891974) and (0.999318044,0.977273413) .. (0.999364364,0.978673002)
    .. controls (0.999402374,0.979821508) and (0.999440385,0.980982236) .. (0.999478395,0.982157498)
    .. controls (0.999509587,0.983121922) and (0.999540778,0.984096133) .. (0.999571970,0.985082017)
    .. controls (0.999597566,0.985891036) and (0.999623161,0.986707916) .. (0.999648757,0.987534193)
    .. controls (0.999669761,0.988212239) and (0.999690765,0.988896613) .. (0.999711769,0.989588572)
    .. controls (0.999729005,0.990156396) and (0.999746241,0.990729326) .. (0.999763477,0.991308391)
    .. controls (0.999777621,0.991783574) and (0.999791765,0.992262886) .. (0.999805909,0.992747170)
    .. controls (0.999817515,0.993144574) and (0.999829122,0.993545325) .. (0.999840728,0.993950110)
    .. controls (0.999850252,0.994282278) and (0.999859777,0.994617163) .. (0.999869301,0.994955327)
    .. controls (0.999877117,0.995232825) and (0.999884932,0.995512531) .. (0.999892748,0.995794907)
    .. controls (0.999899162,0.996026625) and (0.999905575,0.996260141) .. (0.999911989,0.996495831)
    .. controls (0.999917252,0.996689239) and (0.999922515,0.996884112) .. (0.999927778,0.997080760)
    .. controls (0.999932097,0.997242129) and (0.999936415,0.997404693) .. (0.999940734,0.997568705)
    .. controls (0.999944278,0.997703294) and (0.999947822,0.997838858) .. (0.999951366,0.997975606)
    .. controls (0.999954275,0.998087821) and (0.999957183,0.998200834) .. (0.999960091,0.998314814)
    .. controls (0.999962478,0.998408346) and (0.999964864,0.998502529) .. (0.999967251,0.998597504)
    .. controls (0.999969209,0.998675441) and (0.999971167,0.998753911) .. (0.999973126,0.998833028)
    .. controls (0.999974733,0.998897952) and (0.999976340,0.998963312) .. (0.999977947,0.999029202)
    .. controls (0.999979266,0.999083272) and (0.999980584,0.999137698) .. (0.999981903,0.999192559)
    .. controls (0.999982985,0.999237578) and (0.999984068,0.999282890) .. (0.999985150,0.999328557)
    .. controls (0.999986038,0.999366032) and (0.999986926,0.999403746) .. (0.999987814,0.999441752)
    .. controls (0.999988543,0.999472940) and (0.999989271,0.999504323) .. (0.999990000,0.999535946);
  \draw[model curve, ->] (0,1)
    .. controls (0.004361010,0.989097475) and (0.008722020,0.978183067) .. (0.013083030,0.967256234)
    .. controls (0.017444039,0.956329401) and (0.021805049,0.945390144) .. (0.026166059,0.934437889)
    .. controls (0.030527069,0.923485634) and (0.034888079,0.912520381) .. (0.039249089,0.901541521)
    .. controls (0.043610098,0.890562662) and (0.047971108,0.879570195) .. (0.052332118,0.868563474)
    .. controls (0.056693128,0.857556753) and (0.061054138,0.846535778) .. (0.065415148,0.835499859)
    .. controls (0.069776157,0.824463940) and (0.074137167,0.813413078) .. (0.078498177,0.802346539)
    .. controls (0.082859187,0.791280000) and (0.087220197,0.780197783) .. (0.091581207,0.769099107)
    .. controls (0.095942216,0.758000431) and (0.100303226,0.746885296) .. (0.104664236,0.735752866)
    .. controls (0.109025246,0.724620435) and (0.113386256,0.713470711) .. (0.117747266,0.702302800)
    .. controls (0.122108276,0.691134888) and (0.126469285,0.679948791) .. (0.130830295,0.668743552)
    .. controls (0.135191305,0.657538314) and (0.139552315,0.646313935) .. (0.143913325,0.635069393)
    .. controls (0.148274335,0.623824852) and (0.152635344,0.612560149) .. (0.156996354,0.601274188)
    .. controls (0.161357364,0.589988227) and (0.165718374,0.578681009) .. (0.170079384,0.567351359)
    .. controls (0.174440394,0.556021708) and (0.178801403,0.544669625) .. (0.183162413,0.533293846)
    .. controls (0.187523423,0.521918067) and (0.191884433,0.510518593) .. (0.196245443,0.499094063)
    .. controls (0.200606453,0.487669534) and (0.204967462,0.476219950) .. (0.209328472,0.464743848)
    .. controls (0.213689482,0.453267745) and (0.218050492,0.441765123) .. (0.222411502,0.430234402)
    .. controls (0.226772512,0.418703682) and (0.231133522,0.407144862) .. (0.235494531,0.395556237)
    .. controls (0.239855541,0.383967611) and (0.244216551,0.372349180) .. (0.248577561,0.360699095)
    .. controls (0.252938571,0.349049010) and (0.257299581,0.337367273) .. (0.261660590,0.325651879)
    .. controls (0.266021600,0.313936486) and (0.270382610,0.302187437) .. (0.274743620,0.290402560)
    .. controls (0.279104630,0.278617682) and (0.283465640,0.266796976) .. (0.287826649,0.254938077)
    .. controls (0.292187659,0.243079179) and (0.296548669,0.231182088) .. (0.300909679,0.219244230)
    .. controls (0.305270689,0.207306372) and (0.309631699,0.195327747) .. (0.313992708,0.183305546)
    .. controls (0.318353718,0.171283344) and (0.322714728,0.159217566) .. (0.327075738,0.147105139)
    .. controls (0.331436748,0.134992713) and (0.335797758,0.122833637) .. (0.340158767,0.110624547)
    .. controls (0.344519777,0.098415456) and (0.348880787,0.086156351) .. (0.353241797,0.073843536)
    .. controls (0.357602807,0.061530721) and (0.361963817,0.049164195) .. (0.366324827,0.036739893)
    .. controls (0.370685836,0.024315591) and (0.375046846,0.011833513) .. (0.379407856,-0.000710826)
    .. controls (0.383768866,-0.013255165) and (0.388129876,-0.025861765) .. (0.392490886,-0.038535583)
    .. controls (0.396851895,-0.051209401) and (0.401212905,-0.063950438) .. (0.405573915,-0.076764187)
    .. controls (0.409934925,-0.089577936) and (0.414295935,-0.102464399) .. (0.418656945,-0.115429679)
    .. controls (0.423017954,-0.128394960) and (0.427378964,-0.141439059) .. (0.431739974,-0.154568781)
    .. controls (0.436100984,-0.167698503) and (0.440461994,-0.180913849) .. (0.444823004,-0.194222422)
    .. controls (0.449184013,-0.207530995) and (0.453545023,-0.220932796) .. (0.457906033,-0.234436351)
    .. controls (0.462267043,-0.247939907) and (0.466628053,-0.261545218) .. (0.470989063,-0.275261875)
    .. controls (0.475350073,-0.288978532) and (0.479711082,-0.302806536) .. (0.484072092,-0.316756714)
    .. controls (0.488433102,-0.330706892) and (0.492794112,-0.344779245) .. (0.497155122,-0.358986041)
    .. controls (0.501516132,-0.373192836) and (0.505877141,-0.387534075) .. (0.510238151,-0.402023713)
    .. controls (0.514599161,-0.416513352) and (0.518960171,-0.431151389) .. (0.523321181,-0.445953770)
    .. controls (0.527682191,-0.460756150) and (0.532043200,-0.475722873) .. (0.536404210,-0.490872234)
    .. controls (0.540765220,-0.506021595) and (0.545126230,-0.521353593) .. (0.549487240,-0.536889326)
    .. controls (0.553848250,-0.552425058) and (0.558209259,-0.568164524) .. (0.562570269,-0.584132174)
    .. controls (0.566931279,-0.600099825) and (0.571292289,-0.616295659) .. (0.575653299,-0.632748177)
    .. controls (0.580014309,-0.649200694) and (0.584375319,-0.665909895) .. (0.588736328,-0.682909189)
    .. controls (0.593097338,-0.699908482) and (0.597458348,-0.717197868) .. (0.601819358,-0.734816801)
    .. controls (0.606180368,-0.752435735) and (0.610541378,-0.770384215) .. (0.614902387,-0.788709035)
    .. controls (0.619263397,-0.807033854) and (0.623624407,-0.825735013) .. (0.627985417,-0.844868928)
    .. controls (0.630715154,-0.856845641) and (0.633444892,-0.868991909) .. (0.636174629,-0.881324023)
    .. controls (0.636730184,-0.883833852) and (0.637285740,-0.886351379) .. (0.637841296,-0.888876753)
    .. controls (0.638396851,-0.891402127) and (0.638952407,-0.893935349) .. (0.639507962,-0.896476567)
    .. controls (0.640028124,-0.898855884) and (0.640548285,-0.901242212) .. (0.641068446,-0.903635674)
    .. controls (0.641103841,-0.903798536) and (0.641139235,-0.903961432) .. (0.641174629,-0.904124360)
    .. controls (0.641730184,-0.906681722) and (0.642285740,-0.909247234) .. (0.642841296,-0.911821052)
    .. controls (0.643396851,-0.914394870) and (0.643952407,-0.916976994) .. (0.644507962,-0.919567586)
    .. controls (0.645063518,-0.922158177) and (0.645619073,-0.924757236) .. (0.646174629,-0.927364928)
    .. controls (0.646730184,-0.929972620) and (0.647285740,-0.932588944) .. (0.647841296,-0.935214070)
    .. controls (0.648396851,-0.937839197) and (0.648952407,-0.940473126) .. (0.649507962,-0.943116031)
    .. controls (0.650063518,-0.945758936) and (0.650619073,-0.948410818) .. (0.651174629,-0.951071855)
    .. controls (0.651730184,-0.953732891) and (0.652285740,-0.956403083) .. (0.652841296,-0.959082613)
    .. controls (0.653278022,-0.961189014) and (0.653714749,-0.963301185) .. (0.654151476,-0.965419219)
    .. controls (0.654270305,-0.965995514) and (0.654389133,-0.966572242) .. (0.654507962,-0.967149407)
    .. controls (0.655063518,-0.969847802) and (0.655619073,-0.972555724) .. (0.656174629,-0.975273366)
    .. controls (0.656730184,-0.977991008) and (0.657285740,-0.980718370) .. (0.657841296,-0.983455651)
    .. controls (0.658396851,-0.986192932) and (0.658952407,-0.988940132) .. (0.659507962,-0.991697454)
    .. controls (0.660063518,-0.994454777) and (0.660619073,-0.997222222) .. (0.661174629,-1.000000000)
    .. controls (0.661730184,-1.002777778) and (0.662285740,-1.005565888) .. (0.662841296,-1.008364546)
    .. controls (0.663396851,-1.011163204) and (0.663952407,-1.013972410) .. (0.664507962,-1.016792386)
    .. controls (0.665063518,-1.019612361) and (0.665619073,-1.022443106) .. (0.666174629,-1.025284848)
    .. controls (0.666527921,-1.027091986) and (0.666881213,-1.028903570) .. (0.667234505,-1.030719662)
    .. controls (0.667436769,-1.031759393) and (0.667639032,-1.032800602) .. (0.667841296,-1.033843299)
    .. controls (0.668396851,-1.036707269) and (0.668952407,-1.039582471) .. (0.669507962,-1.042469144)
    .. controls (0.670063518,-1.045355817) and (0.670619073,-1.048253962) .. (0.671174629,-1.051163827)
    .. controls (0.671730184,-1.054073692) and (0.672285740,-1.056995276) .. (0.672841296,-1.059928835)
    .. controls (0.673396851,-1.062862394) and (0.673952407,-1.065807927) .. (0.674507962,-1.068765696)
    .. controls (0.675063518,-1.071723466) and (0.675619073,-1.074693473) .. (0.676174629,-1.077675986)
    .. controls (0.676730184,-1.080658498) and (0.677285740,-1.083653518) .. (0.677841296,-1.086661322)
    .. controls (0.678396851,-1.089669125) and (0.678952407,-1.092689713) .. (0.679507962,-1.095723372)
    .. controls (0.679777820,-1.097196952) and (0.680047677,-1.098673617) .. (0.680317535,-1.100153399)
    .. controls (0.680603233,-1.101720043) and (0.680888931,-1.103290181) .. (0.681174629,-1.104863854)
    .. controls (0.681730184,-1.107923948) and (0.682285740,-1.110997407) .. (0.682841296,-1.114084535)
    .. controls (0.683396851,-1.117171663) and (0.683952407,-1.120272459) .. (0.684507962,-1.123387237)
    .. controls (0.685063518,-1.126502015) and (0.685619073,-1.129630774) .. (0.686174629,-1.132773837)
    .. controls (0.688583274,-1.146400780) and (0.690991919,-1.160296593) .. (0.693400565,-1.174488895)
    .. controls (0.697761574,-1.200184987) and (0.702122584,-1.226853015) .. (0.706483594,-1.254696516)
    .. controls (0.710844604,-1.282540017) and (0.715205614,-1.311558991) .. (0.719566624,-1.342019290)
    .. controls (0.723927633,-1.372479589) and (0.728288643,-1.404381214) .. (0.732649653,-1.438075617)
    .. controls (0.737010663,-1.471770020) and (0.741371673,-1.507257203) .. (0.745732683,-1.545012127)
    .. controls (0.750093692,-1.582767051) and (0.754454702,-1.622789716) .. (0.758815712,-1.665737738)
    .. controls (0.763176722,-1.708685760) and (0.767537732,-1.754559139) .. (0.771898742,-1.804294104)
    .. controls (0.776259751,-1.854029069) and (0.780620761,-1.907625622) .. (0.784981771,-1.966465067)
    .. controls (0.789342781,-2.025304512) and (0.793703791,-2.089386851) .. (0.798064801,-2.160826596)
    .. controls (0.802425811,-2.232266341) and (0.806786820,-2.311063493) .. (0.811147830,-2.400659288)
    .. controls (0.815508840,-2.490255083) and (0.819869850,-2.590649521) .. (0.824230860,-2.707680931)
    .. controls (0.828591870,-2.824712340) and (0.832952879,-2.958380721) .. (0.837313889,-3.120000000);
  \draw[model curve,dash pattern=on 4pt off 3pt]
    (0.661174629,-1) -- (1,-2.694126853);
  \fill (0.661174629,-1) circle[radius=1.4pt];
  \node[above,inner sep=2pt] at (.67,.28) {$0<\mu<2K$};
  \node[below left,inner sep=2pt] at (.46,-.30) {$\mu>2K$};
  \node[right,inner sep=3pt] at (.88,-1.94) {$\widehat R_{K,\mu}$};
  \draw[fill=white,line width=.65pt] (1,1) circle[radius=1.5pt];
  \draw[fill=white,line width=.65pt] (1,-1) circle[radius=1.5pt];
  \draw[fill=white,line width=.65pt] (1,-2.694126853) circle[radius=1.5pt];
\end{tikzpicture}
\caption{The model Riccati function $R_{K,\mu}$ and the tangent continuation $\widehat R_{K,\mu}$ (dashed line) for $K>0$.}

\label{fig:intro-Riccati-convexity}
\end{figure}
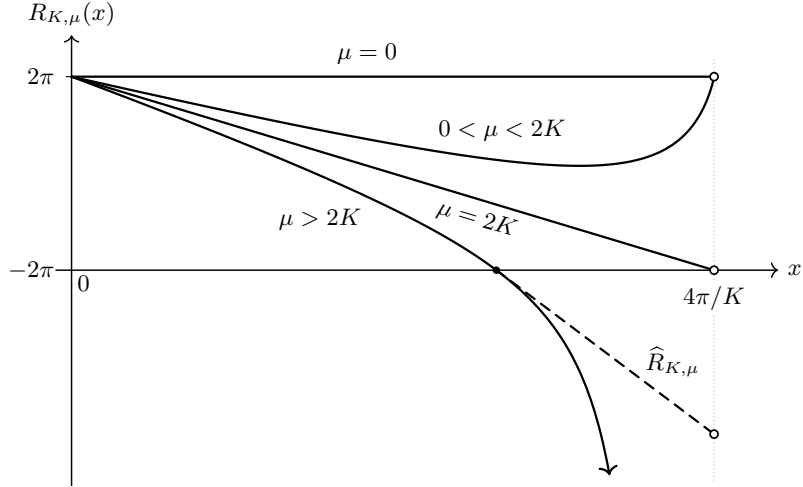

See Figure~\ref{fig:intro-Riccati-convexity}. For $0<\mu<2K$, $R_{K,\mu}$ is strictly convex on $(0,4\pi/K)$; for $\mu = 2K$, $R_{K,2K}(x)=2\pi-Kx$ is affine;
for $\mu>2K$, it is strictly concave throughout the interval on
which the radial part of the first-angular Helmholtz--Steklov
eigenfunction remains positive. The dashed line is the tangent continuation from the first point
where $R_{K,\mu}=-2\pi$.
The resulting function $\widehat R_{K,\mu}$, defined in
\eqref{eq:model-Rhat-definition} and used in the proof, is
well-defined and $C^1$ on $(0,4\pi/K)$ and retains the
(not necessarily strict) concavity of $R_{K,\mu}$.

For $K>0$ and $-4\pi<\alpha\leq4\pi$, the parameter
$\mu$ arising in the proof is the unique positive solution of
\[
\widehat R_{K,\mu}(m_1)+\widehat R_{K,\mu}(m_2)+\alpha=0,
\]
where $m_1,m_2>0$ are the areas of the two Green regions arising
in the proof, with $m_1+m_2=A$.
Since
\[
\widehat R_{K,2K}(m_1)+\widehat R_{K,2K}(m_2)
=4\pi-K(m_1+m_2)=4\pi-KA,
\]
strict monotonicity in $\mu$ gives
\[
\begin{aligned}
\mu=2K
&\quad\Longleftrightarrow\quad \alpha=-4\pi+KA,\\
\mu>2K
&\quad\Longleftrightarrow\quad \alpha>-4\pi+KA,\\
0<\mu<2K
&\quad\Longleftrightarrow\quad -4\pi<\alpha<-4\pi+KA,
\end{aligned}
\]
independently of the area split.

For each fixed $\alpha$, the final comparison value is obtained by
maximizing $\mu$ over $m_1,m_2\geq0$ with $m_1+m_2=A$.
The change in convexity at $\mu=2K$ therefore determines the
transition in the maximizing splits at $\alpha=-4\pi+KA$.
At the critical value, every split gives $\mu=2K$.
For $-4\pi+KA<\alpha\leq4\pi$, the unique maximizing split is
$m_1=m_2=A/2$. For $-4\pi<\alpha<-4\pi+KA$, the maximum is
attained only at the degenerate splits
$\{m_1,m_2\}=\{0,A\}$; see
Proposition~\ref{prop:model-split-optimization}.

In the range $-4\pi<\alpha<-4\pi+KA$, the condition \[\lambda_{1,\mathrm{ang}}\left(
D_K(A);\frac{\alpha+2\pi}{\ell_K(A)}\right)
=
\lambda_{2,\mathrm{rad}}\left(
D_K(A);\frac{\alpha+2\pi}{\ell_K(A)}\right)\] is characterized by
\begin{equation}\label{eq:intro-angular-radial-crossing}
	R_{K,\mu}(A)+\alpha+2\pi=0,
\end{equation}
on the branch $-KA<\alpha<-2\pi$ and $0<\mu<2K$, where 
\[
\mu=-\frac{(\alpha+2\pi)(\alpha+KA)}{A(4\pi-KA)}.\]
The derivation is given in Appendix~\ref{sec:Appendix-mode-crossing}.

Our second result concerns the Steklov--Helmholtz eigenvalue problem
\begin{equation}\label{eq:intro-SH-problem}
	\begin{cases}
		-\Delta_g u=\lambda u &\text{in }\region,\\
		\partial_\nu u=\sigma u &\text{on }\partial\region.
	\end{cases}
\end{equation}
Away from the Dirichlet spectrum, this is the spectral problem for
the Dirichlet-to-Neumann operator of $-\Delta_g-\lambda$; at
$\lambda=0$ it is the classical Steklov problem. The dependence on
$\lambda$ and the passage through Dirichlet eigenvalues are important
features of this family; see the survey of Grebenkov, Levitin and
Polterovich~\cite{grebenkov2026spectral}. Nigam, Patil and Sun proved
a sharp disk comparison for the second Steklov--Helmholtz eigenvalue
below the first Dirichlet threshold in the planar setting, and their
numerical work on the third eigenvalue suggested degeneration toward
two equal disks~\cite{nigam2025spectrum}.

If $\lambda$ is not in the Dirichlet
spectrum, write $\sigma_j(\region,\lambda)$ for the ordered
Steklov--Helmholtz eigenvalues of \eqref{eq:intro-SH-problem}. Letting
\(
N^D_\Omega(\lambda)
:=\#\{j:\lambda_j^D(\region)<\lambda\}
\)
, we consider the Steklov--Helmholtz eigenvalues with Dirichlet index shift
\[
\Sigma_k(\region,\lambda)=
\begin{cases}
	\sigma_{k-N^D_\Omega(\lambda)}(\region,\lambda),
	&N^D_\Omega(\lambda)<k,\\
	-\infty,&N^D_\Omega(\lambda)\geq k.
\end{cases}
\]
We can extend the definition of $\Sigma_k(\region, \lambda)$ to any $\lambda \in \mathbb R$ by Robin--Steklov duality; see Section~\ref{sec:Robin-Steklov-duality}.
\begin{equation}\label{eq:intro-shifted-SH-definition}
	\Sigma_k(\region,\lambda)
	:=-\sup\left\{\beta\in\mathbb R:
	\lambda_k(\region;\beta)\leq\lambda\right\}.
\end{equation}

Define
\[
\lambda_+(K,A)
:=\lambda_2\left(D_K(A/2);\frac{2\pi}{\ell_K(A/2)}\right),
\qquad
\Lambda_+(K,A):=A\lambda_+(K,A).
\]
For $A>0$, with $KA<4\pi$ when $K>0$, let
\begin{align*}
	\mathcal Q_1
	&:=\{(K,A,\Lambda):K\leq0,
	\ 0\leq\Lambda\leq\Lambda_+(K,A)\}\\
	&\quad\cup\{(K,A,\Lambda):K>0,
	\ \Lambda=0\ \text{or}\ 2KA\leq\Lambda\leq\Lambda_+(K,A)\},\\
	\mathcal Q_2
	&:=\{(K,A,\Lambda):K>0,
	\ 0\leq\Lambda\leq2KA\}.
\end{align*}
These are the counterparts of $\mathcal P_1$ and $\mathcal P_2$
under Robin--Steklov duality.

\begin{theorem}\label{thm:intro-Steklov-Helmholtz}
	Under the geometric assumptions of
	Theorem~\ref{thm:intro-Robin-main}, let $\Lambda\geq0$.
	
	If $(K,A,\Lambda)\in\mathcal Q_1$, then
	\begin{equation}\label{eq:intro-Steklov-Helmholtz-two}
		L\,\Sigma_3\left(\region,\frac{\Lambda}{A}\right)
		<
		2\ell_K(A/2)\,
		\Sigma_2\left(D_K(A/2),\frac{\Lambda}{A}\right),
	\end{equation}
	and the bound is sharp.
	
	If $(K,A,\Lambda)\in\mathcal Q_2$, then
	\begin{equation}\label{eq:intro-Steklov-Helmholtz-one}
		L\,\Sigma_3\left(\region,\frac{\Lambda}{A}\right)
		<
		2\pi+
		\ell_K(A)\,\sigma_{1,\mathrm{ang}}
		\left(D_K(A),\frac{\Lambda}{A}\right),
	\end{equation}
	This inequality is sharp among bounded simply connected smooth domains of area $A$ in $\mathbb{M}_K^2$ whenever
    \begin{equation}\label{eq:shifted-SH-one-disk-sharp}
    \beta:=-\sigma_{1,\rm ang}\left(D_K(A),\frac{\Lambda}{A}\right)<0,
    \qquad
    \lambda_{1,\mathrm{ang}}(D_K(A);\beta)\leq
    \lambda_{2,\mathrm{rad}}\left(D_K(A);\beta\right).
    \end{equation}
\end{theorem}

For $K>0$, the two comparison values in the theorem agree at
$\Lambda=0$ and $\Lambda=2KA$. Section~\ref{sec:Robin-Steklov-duality} proves the
theorem, including the Dirichlet index shift and the passage from
Robin sharpness sequences to the shifted Steklov--Helmholtz
inequalities.

\subsection{Ideas of the proof}
The argument brings together a two-pole Green-level construction,
simultaneous topological orthogonality, and a pointwise Riccati
comparison. The Green-level approach was developed for one-pole
Neumann inequalities by Provenzano and Savo
~\cite{provenzano2026isoperimetric}; the two-pole construction on
closed surfaces introduced by Michetti, Provenzano and Savo
~\cite{michettiprovenzanosavo2026} is particularly close to the
geometry needed here. We adapt two-pole Green functions to domains
with boundary and combine them with a Robin variational problem whose
minimizer depends on the chosen poles.

Let $\psi_a$ be the Dirichlet Green function with pole $a\in\region$.
For distinct $a,b\in\region$, write
$\psi_{ab}:=\psi_a-\psi_b$ and let $e^{i\Theta_{ab}}$ be the
associated single-valued phase. Our trial functions have the form
\begin{equation}\label{eq:intro-trial}
	U_{a,b}[g]:=(g\circ\psi_{ab})e^{i\Theta_{ab}}.
\end{equation}
The sign domains $\region_a=\{\psi_{ab}>0\}$ and
$\region_b=\{\psi_{ab}<0\}$ have areas $m_a,m_b$ with
$m_a+m_b=A$. By the coarea formula the Rayleigh quotient becomes
\begin{equation}\label{eq:intro-g}
	\frac{\displaystyle
		\int_{\mathbb R}\bigl(|g'|^2+4\pi^2|g|^2\bigr)\,dt
		+\alpha|g(0)|^2}
	{\displaystyle\int_{\mathbb R}|g|^2\,d\nu_{a,b}},
\end{equation}
where $\nu_{a,b}$ is the push-forward of area under $\psi_{ab}$.
For each pole pair, $g=g_{a,b}$ is chosen as the positive normalized
minimizer of this quotient. The remaining steps are as follows.

\emph{First, simultaneous orthogonality.}
The trial function must be orthogonal to the first two real Robin
eigenfunctions. For a real function $w$, consider
\[
F_w(a,b):=\int_\region U_{a,b}[g_{a,b}]\,w\,dA_g.
\]
A key analytic issue is the extension of this map to the compactified
pole space. When one pole reaches the boundary, its Green region
collapses and the minimization on the disappearing half-line contributes
a shift $\alpha\mapsto\alpha+2\pi$. After conformal uniformization,
the extension has a pole-exchange symmetry and boundary factorization
\begin{equation}\label{eq:W-boundary-data}
	F(a,b)=J(a)\overline{G(b)}
	\qquad (a\in\mathbb S^1,\ b\in\overline\disk).
\end{equation}
The following degree theorem supplies the two orthogonality
conditions.

Let $M_0:=\overline\disk\times\overline\disk$, oriented as a
product; in the topological proof below, this parameter space is also
denoted by $M$. For a nonvanishing continuous map
$J:\mathbb S^1\to\mathbb C\setminus\{0\}$, put
$\degmap J:=\wind(J,0)$.
For an integer $d$, define $\mathcal W_d$ to be the class of
continuous maps $F:M_0\to\mathbb C$ satisfying
\begin{enumerate}[(a)]
	\item $F(b,a)=\overline{F(a,b)}$ for all $(a,b)\in\partial M_0$;
	\item the factorization \eqref{eq:W-boundary-data} holds for
	continuous $G:\overline\disk\to\mathbb C$ and
	$J:\mathbb S^1\to\mathbb C\setminus\{0\}$;
	\item $\degmap J=d$.
\end{enumerate}
\begin{theorem}\label{thm:topo}
	Suppose $F_k\in\mathcal W_{d_k}$ for $k=1,2$, and let
	$\Phi:=(F_1,F_2):M_0\to\mathbb C^2$. If $\Phi$ has no zero on
	$\partial M_0$, then the normalized boundary map
	$\widehat\Phi:=\Phi/\|\Phi\|:\partial M_0\to\mathbb S^3$
	has degree
	\begin{equation}\label{eq:degree-formula}
		\deg(\widehat\Phi)=-2d_1d_2.
	\end{equation}
\end{theorem}
Consequently, whenever $d_1d_2\neq0$, the two components have a
common zero, either on the boundary or, by nonzero degree, in the
interior. Unlike the reflection-symmetric degree construction for
planar Robin eigenvalues~\cite{kimlaugesen2025}, the theorem is
adapted to the pole-exchange symmetry and boundary factorization of
the two-pole Green functions.

\emph{Second, comparison to the constant curvature model.}
In the area variable, the minimum in \eqref{eq:intro-g} takes the
form of a coupled one-dimensional Robin problem:
\begin{equation}\label{eq:intro-Robin-kappa-hat}
	\widehat\kappa_\alpha(G_a,G_b)
	:=\inf_{\substack{f_j\in\mathcal F(G_j)\\
			f_a(m_a)=f_b(m_b)=:c}}
	\frac{\displaystyle
		\sum_{j=a,b}\int_0^{m_j}
		\left(G_j|f_j'|^2+\frac{4\pi^2}{G_j}|f_j|^2\right)dx
		+\alpha|c|^2}
	{\displaystyle\sum_{j=a,b}\int_0^{m_j}|f_j|^2\,dx}.
\end{equation}
Here $G_j$ denotes the Green-level weight of the corresponding
region, and $\mathcal F(G_j)$ is the admissible class of regular
finite-energy profiles, specified in Section~\ref{sec:main}.
For the analogous two-pole variational problem on a closed simply
connected surface, with no Robin term, Michetti, Provenzano and Savo~\cite[Lemma~4.3]{michettiprovenzanosavo2026}
proved monotonicity with respect to the Green weight. Their
Feynman--Hellmann argument relies on $R^2\leq4\pi^2$, where
$R=Gf'/f$. However, in our setting, positive Robin parameters may lead to
$R(m)<-2\pi$.

Instead, we compare Riccati functions pointwise. The direct
comparison with $R_{K,\mu}$ applies until the model first reaches
$-2\pi$. Continuing the model tangentially there by a line of slope
$-\mu$ gives $\widehat R_{K,\mu}$ and extends the comparison to the
whole area interval; see Lemma~\ref{lem:Robin-tangent-comparison}.
This continuation gives an upper bound, but sharpness requires
identifying the final comparison value with a genuine model
eigenvalue. In the two-disk range $\mathcal P_1$, the final optimized value
$\mu^*:=\mu_K^{\mathrm{cmp}}(A/2,A/2;\alpha)$ satisfies
\[
\widehat R_{K,\mu^*}(A/2)=-\frac\alpha2\geq-2\pi
\qquad(\alpha\leq4\pi).
\]
Thus $\widehat R_{K,\mu^*}(A/2)=R_{K,\mu^*}(A/2)$, and $\mu^*$
is precisely the second Robin eigenvalue on the right-hand side of
\eqref{eq:intro-Robin-main}. Thus, the condition $\alpha\leq4\pi$ ensures that the tangential
continuation causes no loss in the sharpness of the bound.

The comparison can now be summarized as follows. Set
$\kappa:=\widehat\kappa_\alpha(G_a,G_b)$ and
\[
R_j:=G_j\frac{f_j'}{f_j},
\qquad R_j'=\frac{4\pi^2-R_j^2}{G_j}-\kappa,
\qquad R_j(0+)=2\pi.
\]
The curvature comparison gives
\[
G_j(x)\geq G_K(x):=4\pi x-Kx^2
\quad\Longrightarrow\quad
R_j(x)\leq\widehat R_{K,\kappa}(x).
\]
Define $\mu_K^{\mathrm{cmp}}(m_a,m_b;\alpha)$ as the unique
nonnegative root $\mu$ of
\[
\widehat R_{K,\mu}(m_a)+\widehat R_{K,\mu}(m_b)+\alpha=0.
\]
Then $\lambda_3(\region;\alpha/L)\leq\kappa
\leq\mu_K^{\mathrm{cmp}}(m_a,m_b;\alpha)$; see Steps~1--2 of the
main proof in Subsection~\ref{sec:main-proof}.

\emph{Third, optimization of the area split.}
The optimal area split is determined by the convexity or concavity
of $\widehat R_{K,\mu}$; see
Propositions~\ref{prop:model-area-concavity}
and~\ref{prop:model-split-optimization}.
To prove these properties, we compare $R_{K,\mu}$ with
$R_{K,2K}(x)=2\pi-Kx$. For $K>0$, the latter corresponds to a
degree-one spherical harmonic, whose radial part is proportional
to $\ell_K(x)$. Green's identity with this explicit comparison
function yields
\[
R_{K,\mu}(x)-(2\pi-Kx)
=\frac{2K-\mu}{B(x)}\int_0^x B(y)\,dy,
\]
where 
\(
B(x):=\ell_K(x)f_{K,\mu}(x);
\) see Lemma~\ref{lem:model-R-critical-identity}.
The quotient in this identity is strictly convex for every
$\mu>0$ wherever $B>0$, regardless of the sign of $2K-\mu$.

This follows from a Sturm--Liouville convexity lemma:
if $p>0$ and $(x^2p(x))'>0$, then every positive regular solution
of $u''+pu=0$ with $u(0)=0$ and $u'(0)>0$ has a strictly convex
quotient $\int_0^x u(y)\,dy/u(x)$
(Lemma~\ref{lem:Sturm-integral-ratio-convexity}). In our case,
\[
B''+\frac{\mu}{G_K}B=0,
\qquad
\left(\frac{\mu x^2}{G_K(x)}\right)'
=\frac{4\pi\mu}{(4\pi-Kx)^2}>0.
\]
The difference formula therefore makes $R_{K,\mu}''$ have the sign
of $2K-\mu$, while tangential continuation preserves the resulting
non-strict convexity or concavity.

\emph{Finally, boundary oscillations and sharpness.}
The Robin convergence results in
Appendix~\ref{sec:Appendix_Robin_convergence} show that spectral convergence
depends on both the convergence of the domains and that of their
boundary measures. Boundary oscillations allow us to change the
latter without changing the limiting domains. In particular,
Lemma~\ref{lem:oscillating-Robin-spectrum} constructs domains
$D_n\to D$ with
\[
ds_{\partial D_n}\rightharpoonup c\,ds_{\partial D},
\qquad
\lambda_k(D_n;\gamma)\longrightarrow\lambda_k(D;c\gamma),
\qquad c\geq1.
\]
Thus a single Robin coefficient on an approximating domain can
produce different effective coefficients on two limiting disks.
Indeed, let $D_i=D_K(m_i)$ and $\ell_i=|\partial D_i|$.
Oscillating the boundary of $D_2$ with multiplier $c$ and joining
it to $D_1$ by a short neck gives a limiting Robin spectrum on
$D_1\sqcup D_2$ with respective coefficients
\[
\beta_1=\frac{\alpha}{\ell_1+c\ell_2},
\qquad
\beta_2=\frac{c\alpha}{\ell_1+c\ell_2}.
\]
In the sharpness constructions, we choose $c$ to redistribute the
fixed total Robin boundary mass
$\beta_1\ell_1+\beta_2\ell_2=\alpha$ between the disks
so that their first angular eigenvalues agree; see
Proposition~\ref{prop:sharpness-coefficient-approximation}.

For the one-disk sharpness sequence, the large and small disks
carry limiting Robin boundary masses $\alpha+2\pi$ and $-2\pi$,
respectively. This explains the restriction $\alpha\leq-2\pi$:
above it the required masses have opposite signs, which a single
constant Robin coefficient cannot produce. The angular--radial
condition in \eqref{eq:intro-Robin-one-disk-sharp} ensures that
the second eigenvalue of the large disk is its first angular
eigenvalue, so the limiting third eigenvalue reaches the stated
bound. The spectral argument first takes the oscillation and neck
limits at fixed disk areas, and only then lets the small disk area
tend to zero. This iterated limit yields both the spectral limit
and the boundary concentration described above.

\subsection{Organization.}
Section~\ref{sec:Green} recalls the Green functions and their phase
functions. Section~\ref{sec:main} develops the two-pole test functions,
the coarea reduction, and the reduction of the Robin inequalities to
the fixed-split comparison and orthogonality results.
Sections~\ref{sec:kappa} and~\ref{sec:Robin-model-concavity} prove the
Riccati comparison and optimize the model area split, respectively.
Section~\ref{sec:Robin-orthogonality-test} constructs the orthogonal
trial functions, and Section~\ref{sec:topo} proves the degree theorem.
Section~\ref{sec:sharpness} treats sharpness and counterexamples.
Section~\ref{sec:Robin-Steklov-duality} establishes the
Steklov--Helmholtz consequences. The appendices provide the
supporting Robin spectrum convergence, variational results, basic facts of the endpoint behavior of the Riccati equation, and the angular-radial crossing on spherical disks.

\section{Preliminaries on Green's Functions}\label{sec:Green}

Let $\Omega\subset M$ be a \emph{Jordan domain}, i.e. an open domain  whose boundary is a Jordan curve and is homeomorphic to $\disk$. 
By uniformization, 
it is conformally equivalent to the planar disk, and the homeomorphism extends to a homeomorphism between $\overline{\Omega}$ and $\overline{\disk}$. 
Let $a\in\Omega$. We use the sign
convention in which the \emph{Dirichlet Green function $\psi_a$ with
pole at $a$} solves
\[
\begin{cases}
\Delta\psi_a=-\delta_a & \text{in }\Omega,\\
\psi_a=0 & \text{on }\partial\Omega.
\end{cases}
\]

The one-form $-2\pi *d\psi_a$ is closed on $\Omega\setminus\{a\}$,
and, for a positively oriented simple loop $\gamma$ enclosing $a$
with interior $D_\gamma$, its period is
\[
-2\pi\int_\gamma *d\psi_a
=-2\pi\int_{D_\gamma}\Delta\psi_a
=2\pi.
\]
Thus there is a multivalued smooth angle function $\Theta_a$ on
$\Omega\setminus\{a\}$ satisfying $d\Theta_a=-2\pi *d\psi_a$,
whose branches differ by integer multiples of $2\pi$.
The associated \emph{phase}
\[
e^{i\Theta_a}:\Omega\setminus\{a\}\longrightarrow\mathbb S^1
\]
is single-valued and unique up to multiplication by a constant
in $\mathbb S^1$.

\subsection{Conformal invariance and one-pole phases}

As usual,  ``conformal'' means orientation-preserving conformal.	 
The following conformal invariance is a key property of Green's functions in dimension two.

\begin{proposition}
	\label{prop:conformal}
	Let $\region_1\subset M_1$ and $\region_2\subset M_2$ be Jordan domains with metrics $g_1,g_2$ on $M_1,M_2$, respectively. Let $\Phi:(\region_1,g_1)\to(\region_2,g_2)$ be a conformal bijection that extends to a homeomorphism from $\overline\region_1$ to $\overline\region_2$. Suppose $a\in\region_1$ and $b=\Phi(a)$. If  $\psi_b$ is the Green function of  $\region_2$ with pole at $b$, and $e^{i\Theta_b}$ is a phase associated with $\psi_b$, then
	\begin{enumerate}[(1)]
		\item $\psi_{a} := \Phi^*(\psi_b)$ 
	is the Green function of $\region_1$ with pole at $a$,
	    \item $e^{i\Theta_a} := \Phi^*(e^{i\Theta_b})$ is a phase associated with $\psi_a$.
	\end{enumerate}
\end{proposition}
\begin{proof} 
By invariance of two-dimensional Dirichlet energy, for any $\phi \in C_c^\infty(\region_1)$, 
$$
\int_{\region_1} \langle d(\psi_b\circ\Phi), d\phi\rangle_{g_1}\,dA_{g_1} 
=
\int_{\region_2}\langle d\psi_b, d(\phi\circ\Phi^{-1})\rangle_{g_2} \,dA_{g_2} = \phi\circ\Phi^{-1}(b)= \phi(a).
$$
This implies
$$
\Delta_{g_1}(\psi_b\circ\Phi) = -\delta_a\quad\text{on }\region_1
$$
in the weak sense. 
Moreover, the assumption that $\Phi$ extends to a homeomorphism 	$\overline{\Phi}:\overline{\Omega}_1
	\longrightarrow \overline{\Omega}_2$ 
together with the fact $\psi_b|_{\partial\Omega_2}=0$ implies  $\left.\psi_a\right|_{\partial\Omega_1}=0$. So $\psi_a$ is the Dirichlet Green function of $\Omega_1$ with pole at $a$.
	
The second assertion follows from the fact that  an orientation-preserving conformal map  commutes with the Hodge star on one-forms. Let $\Theta_a=\Phi^*\Theta_b$. Then  
$$
d\Theta_a=\Phi^*(d\Theta_b)=-2\pi\Phi^*(*d\psi_b)=-2\pi *d\psi_a.
$$
Hence $e^{i\Theta_a}  =\Phi^*(e^{i\Theta_b})$ is a  phase associated with $\psi_a$.  
\end{proof}

For the unit disk $\disk\subset\mathbb R^2=\mathbb C$, the Green function with pole at $0$ is
\[
\psi_0(z)=-\frac{1}{2\pi}\log|z|,
\]
and a corresponding phase is $e^{i\theta}=z/|z|$, where $\theta$ denotes the polar angle. For $a\in\disk$, the M\"obius transformation $M_a: \disk \to \disk$ defined by  
\[
M_a(z):=\frac{z-a}{1-\conj a z}
\]
 is an automorphism of $\disk$ sending $a$ to $0$. Hence, by conformal invariance,
\begin{equation}\label{eq:psi}
\psi_a(z)=-\frac{1}{2\pi}\log\left|\frac{z-a}{1-\conj a z}\right|, 
\end{equation}
with a phase function 
\begin{equation*}
	e^{i\theta_a}(z)=\frac{M_a(z)}{|M_a(z)|}.
\end{equation*}

Using conformal invariance, we now transfer these formulas to a general Jordan domain and record the properties needed below.

\begin{proposition}\label{prop:one-pole-properties}
Let $\region\subset M$ be a Jordan domain and let $\Phi:(\region,g)\to(\disk,g_0)$ be a conformal bijection, where $g$ is the metric on $M$ and $g_0$ is the standard Euclidean metric. We use the same symbol $\Phi$ for its homeomorphic extension to the closures. For $a\in\region$, let $\psi_a$ denote the Green function on $\region$ with   pole at $a$, and set
\[
H_a:=M_{\Phi(a)}\circ\Phi,
\qquad
e^{i\Theta_a}:=\frac{H_a}{|H_a|}.
\]
Then:
\begin{enumerate}[(1)]
\item $e^{i\Theta_a}$ is a phase associated with $\psi_a$.  
\item $d\psi_a\neq0$ on $\region\setminus\{a\}$.
\item $\psi_a>0$ in $\region\setminus\{a\}$, 
and, for every $t>0$, the level set $\psi_a^{-1}(t)$ is a smooth Jordan curve enclosing $a$. 
\item   For any $a_0 \in \partial\Omega$,   set $e^{i\Theta_{a_0}}(z)=-\Phi(a_0)$. Then the map 
\[
(a,z)\longmapsto e^{i\Theta_a}(z)
\]
is continuous on $
\{(a,z)\in\overline\region\times\overline\region:a\neq z\}.$
\end{enumerate}
\end{proposition}

\begin{proof}
By Proposition~\ref{prop:conformal} and the disk formulas above,
\[
\psi_a=-\frac{1}{2\pi}\log|H_a|,
\qquad
e^{i\Theta_a}=\frac{H_a}{|H_a|},
\]
so the second identity defines a phase associated with $\psi_a$, proving part~(1).

The map $H_a:\region\to\disk$ is a conformal diffeomorphism sending $a$ to $0$, hence has nowhere-vanishing differential. Since $d\log|z|\neq0$ on $\disk\setminus\{0\}$, the chain rule gives $d\psi_a\neq0$ on $\region\setminus\{a\}$, proving part~(2).

Because $0<|H_a|<1$ in $\region$ and $|H_a|=1$ on $\bdr\region$, we have $\psi_a>0$ in $\region$ and $\psi_a^{-1}(0)=\bdr\region$. Moreover,
\[
\psi_a^{-1}(t)
=H_a^{-1}\bigl(\{|z|=e^{-2\pi t}\}\bigr),
\qquad t>0,
\]
which proves part~(3).

For part~(4), consider the disk formula
\[
(c,w)\longmapsto
\frac{M_c(w)}{|M_c(w)|}
=\frac{(w-c)/(1-\conj c w)}{|(w-c)/(1-\conj c w)|}.
\]
It extends continuously to
\[
\{(c,w)\in\overline\disk\times\overline\disk:c\neq w\}.
\]
Indeed, if $c\in\bdr\disk$ and $w\neq c$, then
\[
1-\conj c w=-\conj c\,(w-c),
\]
so the extended value is the constant $-c$. Composing with $(a,z)\mapsto(\Phi(a),\Phi(z))$ gives precisely the stated extension on $\region$.
\end{proof}

From now on, $e^{i\Theta_a}$ always denotes the extension in Proposition~\ref{prop:one-pole-properties}, with the pole $a$ allowed to lie in $\overline\region$.

\subsection{Two-pole Green's function}

 For $a,b\in\overline\region$, we also set
\begin{equation}\label{eq:extended-two-pole-phase}
	e^{i\Theta_{ab}}(z)
	:=e^{i\Theta_a}(z)\,\overline{e^{i\Theta_b}(z)},
	\qquad
	z\in\overline\region\setminus\{a,b\}.
\end{equation}
Thus
\[
(a,b,z)\longmapsto e^{i\Theta_{ab}}(z)
\]
is continuous whenever $z\neq a,b$. 
Notice   that $e^{i\Theta_{aa}}=1$ away from $a$.

\begin{definition}\label{defn:double_green}
Let $\region\subset M$ be a Jordan domain and let $a,b\in\region$ be distinct. The \emph{Green function with a positive pole at $a$ and a negative pole at $b$} is the unique solution $\psi_{ab}$ of
\begin{equation*}
\begin{cases}
\Delta_g\psi_{ab}=-\delta_a+\delta_b & \text{in }\region,\\
\restrict{\psi_{ab}}{\bdr\region}=0. &
\end{cases}
\end{equation*}
The function $e^{i\Theta_{ab}}$ defined by \eqref{eq:extended-two-pole-phase} is called the \emph{two-pole phase}.
\end{definition}

By linearity,
\[
\psi_{ab}=\psi_a-\psi_b.
\]
Moreover, locally on $\region\setminus\{a,b\}$,
\[
d\Theta_{ab}
=d\Theta_a-d\Theta_b
=-2\pi*(d\psi_a-d\psi_b)
=-2\pi*d\psi_{ab},
\]
so $e^{i\Theta_{ab}}$ is a phase associated with $\psi_{ab}$. By Proposition~\ref{prop:one-pole-properties}, $\psi_{ab}$ extends continuously to $\overline\region\setminus\{a,b\}$ and vanishes on $\bdr\region$.

The following proposition describes the nodal decomposition of $\psi_{ab}$. This decomposition will be used in the construction of the test functions in Section~\ref{sec:main}.

\begin{proposition}\label{prop:double-green}
Let $\region$ be a Jordan domain with metric $g$ and let $a,b\in\region$ be distinct. Set
\[
\region_a:=\{\psi_{ab}>0\},
\qquad
\region_b:=\{\psi_{ab}<0\}.
\]
Then:
\begin{enumerate}[(1)]
\item $\region_a$ and $\region_b$ are Jordan domains, with $a\in\region_a$ and $b\in\region_b$.
\item The common interior boundary
\[
\gamma:=\bdr\region_a\cap\bdr\region_b
\]
is a simple arc with endpoints on $\bdr\region$, and is analytic in the interior of $\region$.
\item The restriction $\restrict{\psi_{ab}}{\region_a}$ is the Green function on $\region_a$ with pole at $a$, with a phase  $\restrict{e^{i\Theta_{ab}}}{\region_a}$.
\item  The restriction  $-\restrict{\psi_{ab}}{\region_b}$ is the Green function on $\region_b$ with pole at $b$, with a phase 
  $\restrict{e^{-i\Theta_{ab}}}{\region_b}$.
\end{enumerate}
\end{proposition}

\begin{proof}
Choose a conformal bijection from $(\disk,g_0)$ to $(\region,g)$ together with its homeomorphic extension to the closures. Composing with a disk automorphism, we 
obtain a conformal bijection $\Phi:(\overline\disk,g_0)\to(\overline\region,g)$ and a number $x\in(0,1)$ such that
\[
\Phi(x)=a,
\qquad
\Phi(-x)=b.
\]
We claim that
\begin{equation}\label{eq:double-green-halfdisks}
\Phi^{-1}(\region_a)=\disk^+:=\disk\cap\{\operatorname{Re}z>0\},
\qquad
\Phi^{-1}(\region_b)=\disk^-:=\disk\cap\{\operatorname{Re}z<0\}.
\end{equation}

By conformal invariance,
\[
\Phi^*(\psi_{ab})=\psi_+-\psi_-,
\]
where $\psi_+$ and $\psi_-$ are the Green functions on $\disk$ with  poles at $x$ and $-x$, respectively. Let
\[
R(z):=-\conj z
\]
denote reflection across the imaginary axis. By symmetry,
\[
\psi_-(R(z))=\psi_+(z),
\]
and hence
\[
\Phi^*(\psi_{ab})(z)=\psi_+(z)-\psi_+(R(z)).
\]
Thus it remains to show that
\[
\psi_+(z)>\psi_+(R(z))
\qquad(z\in\disk^+).
\]

Fix $z\in\disk^+$ and set $t:=\psi_+(z)>0$. The superlevel set
\[
D_t:=\{w\in\disk:\psi_+(w)>t\}
\]
is the disk
\[
D_t
=\left\{w\in\disk:
\left|\frac{w-x}{1-xw}\right|<v\right\},
\qquad
v:=e^{-2\pi t}\in(0,1).
\]
As a Euclidean disk, $D_t$ has center
\[
c_t=\frac{x(1-v^2)}{1-x^2v^2}>0.
\]
Since $z\in\partial D_t$ and $\operatorname{Re}z>0$, reflection across the imaginary axis increases the distance from $z$ to the center $c_t$. Thus $R(z)\notin\overline{D_t}$, and therefore
\[
\psi_+(R(z))<t=\psi_+(z).
\]
This proves \eqref{eq:double-green-halfdisks}.

Because $\Phi:\overline\disk\to\overline\region$ is a homeomorphism, \eqref{eq:double-green-halfdisks} shows that $\region_a$ and $\region_b$ are Jordan domains containing $a$ and $b$, respectively. Their common interior boundary is
\[
\gamma=\Phi\bigl(i[-1,1]\bigr),
\]
which is a simple arc with endpoints on $\bdr\region$. It is analytic away from the endpoints because $\Phi$ is holomorphic with nonvanishing derivative in $\disk$. This proves parts~(1) and~(2).

On $\region_a$ and $\region_b$, respectively,
\[
\Delta\psi_{ab}=-\delta_a,
\qquad
\Delta\psi_{ab}=\delta_b,
\]
and $\psi_{ab}$ vanishes on the corresponding boundaries. This proves part~(3) and ~(4).
\end{proof}

\section{Test Functions and Proof of the Main Theorem}\label{sec:main}
\subsection{Test functions and basic calculation.}\label{sec:calc}
For the Robin problem with Robin parameter $\alpha$, it is natural to start with the quadratic form
\[\mathfrak{q}_{\alpha}[u] := \int_\region|\nabla u|_g^2\,dA_g +\frac\alpha L\int_{\partial\region}|u|^2\,ds_g.\] 
For $a,b\in \overline\region$, consider the test function 
\begin{equation}\label{eq:Robin-dipole-trial}
U_{a,b}[g]
:=(g\circ\psi_{ab})e^{i\Theta_{ab}},
\end{equation}
where $g\in H^1(\mathbb{R})$ is to be determined. For a single pole $a\in \overline\region$, we also define
$$
U_{a}[g]:=(g\circ\psi_a)e^{i\Theta_a}.
$$
The test function is well-defined for any $a,b\in\overline\region$ but degenerates in some cases. We call $a\neq b\in\region$ the \textit{general case}, $a\in\bdr\region, b\in\region$ or $a\in\region, b\in\bdr\region$ the \textit{boundary case}, and $a=b$ or $a,b\in\bdr\region$ the \textit{collapsed case}. In the boundary case, if $a\in\bdr\region$ and $b\in\region$, then the boundary phase convention in Proposition~\ref{prop:one-pole-properties} gives
\[
e^{i\Theta_{ab}}=-\Phi(a)\overline{e^{i\Theta_b}},
\qquad
\psi_{ab}=-\psi_b<0,
\]
and hence
\[
U_{a,b}[g]
=-\Phi(a)\,g(-\psi_b)\,\overline{e^{i\Theta_b}}.
\]
If instead $a\in\region$ and $b\in\bdr\region$, then
\[
U_{a,b}[g]= -\overline{\Phi(b)}\,U_a[g].
\]
Moreover, $U_{a,b}[g]\equiv g(0)$ when $a=b$, while
$U_{a,b}[g]\equiv g(0)\Phi(a)\overline{\Phi(b)}$ when
$a,b\in\bdr\region$.

Now we calculate the Rayleigh quotient. First consider the general case. Recall from Proposition~\ref{prop:double-green} that
\[
\region_a=\{\psi_{ab}>0\},
\qquad
\region_b=\{\psi_{ab}<0\}
\]
are the two Green regions determined by the dipole $\psi_{ab}$, and that $\psi_{ab}|_{\region_a}$ and $-\psi_{ab}|_{\region_b}$ are the corresponding one-pole Green functions. Since $d\Theta_{ab}=-2\pi *d\psi_{ab}$, the one-forms $d\psi_{ab}$ and $d\Theta_{ab}$ are pointwise orthogonal and $|d\Theta_{ab}|_g=2\pi|d\psi_{ab}|_g$. Hence by \eqref{eq:Robin-dipole-trial}, 
\[
|\nabla U_{a,b}[g]|_g^2
=\left(|g'(\psi_{ab})|^2+4\pi^2|g(\psi_{ab})|^2\right)|d\psi_{ab}|_g^2.
\]
For $t>0$, integration of the Green formula over $\{\psi_{ab}>t\}$ gives
\[
\int_{\psi_{ab}=t}|d\psi_{ab}|_g\,ds_g=1,
\]
and the same identity holds for $t<0$ after integrating over the corresponding sublevel region. Therefore, the coarea formula gives
\begin{equation}\label{eq:Robin-universal-energy}
\int_\region|\nabla U_{a,b}[g]|_g^2\,dA_g
=\int_{\mathbb R}\left(|g'|^2+4\pi^2|g|^2\right)dt.
\end{equation}
For $g\in H^1(\mathbb{R})$, since $H^1(\mathbb{R})\subset C_0(\mathbb{R})$, the function $g$ tends to zero at both infinite ends. Hence $U_{a,b}[g]$ is bounded, tends to zero at the interior poles, and has finite Dirichlet energy by \eqref{eq:Robin-universal-energy}; thus $U_{a,b}[g]\in H^1(\region)$. Similarly,
\begin{equation}\label{eq:Robin-universal-mass}
\int_\region|U_{a,b}[g]|^2\,dA_g
=\int_{\mathbb R}|g|^2\,\left(
\int_{\psi_{ab}=t}\frac{1}{|d\psi_{ab}|_g}\,ds_g
\right)dt.
\end{equation}
On $\partial\region$, the continuous boundary value of $U_{a,b}[g]$ has modulus $|g(0)|$. Since $\region$ is Lipschitz, the Sobolev trace of an $H^1(\region)$ function that extends continuously to the boundary agrees $ds_g$-almost everywhere with its pointwise boundary value. Therefore
\begin{equation}\label{eq:Robin-universal-boundary}
\frac\alpha L\int_{\partial\region}|U_{a,b}[g]|^2\,ds_g
=\alpha|g(0)|^2.
\end{equation}
The same trace observation applies to the one-pole test functions. For the boundary case $a\in \partial\region, b\in \region$, the Green function $\psi_{ab}=-\psi_b$ is negative in $\region$. Equation~\eqref{eq:Robin-universal-boundary} still holds, and equations~\eqref{eq:Robin-universal-energy}--\eqref{eq:Robin-universal-mass} hold if the integral domain on the right-hand side is changed from $\mathbb R$ to $\mathbb R^-=(-\infty,0)$. For the other boundary case $a\in\region, b\in\partial\region$, one only needs to replace $\mathbb R$ by $\mathbb R^+=(0,\infty)$.

We now rewrite these Rayleigh quotients in terms of the one-dimensional function $g$. Let
\[
\nu_{a,b}:=(\psi_{ab})_\#(dA_g)
\]
be the push-forward of the area measure $dA_g$ to $\mathbb R$. Thus the total mass $\abs{\nu_{a,b}}= A$. In either the general case or the boundary case, $\nu_{a,b}$ is absolutely continuous with $d\nu_{a,b}=\rho(t)dt$, where
$$
\rho(t) = \int_{\psi_{ab}=t}\frac{1}{|d\psi_{ab}|_g}\,ds_g.
$$
In particular, when $a\in \bdr\region$, the measure $\nu_{a,b} =: \nu_b^-$ depends only on $b\in\region$ and is supported on $(-\infty,0]$. When $b\in\bdr\region$, the measure $\nu_{a,b} =: \nu_a^+$  depends only on $a\in\region$ and is supported on $[0,\infty)$. In the collapsed case, $\nu_{a,b} = A\delta_0$. 

Define 
\begin{equation*}
E_{\alpha}[g] := 
\int_{\mathbb R}\left(|g'|^2+4\pi^2|g|^2\right)dt
+\alpha|g(0)|^2,\qquad g\in H^1(\mathbb{R})
\end{equation*}
and
\begin{equation*}
E_{\beta}^+[g] := 
\int_{\mathbb R^+}\left(|g'|^2+4\pi^2|g|^2\right)dt
+\beta|g(0)|^2,\qquad g\in H^1(\mathbb{R}^+).
\end{equation*}
Then  \eqref{eq:Robin-universal-energy}--\eqref{eq:Robin-universal-boundary} give 
\begin{lemma}
	 For double poles $a,b\in\region$ with $a\neq b$, 
	 \begin{equation}\label{eq:Rayleigh-double}
	 	\frac{\mathfrak{q}_{\alpha}[U_{a,b}[g]]}
	 	{ \int_{\region}\abs{U_{a,b}[g]}^2dA_g}
	 	=
	 	\frac{E_\alpha[g]}
	 	{ \int_{\mathbb R}|g|^2\,d\nu_{a,b}},\qquad \text{for } g\in H^1(\mathbb{R}),
	 \end{equation}
	 while for a single pole $a\in\region$,  
	 \begin{equation}\label{eq:Rayleigh-single}
	 	\frac{\mathfrak{q}_{\beta}[U_{a}[g]]}
	 	{ \int_{\region}\abs{U_{a}[g]}^2dA_g}
	 	=
	 	\frac{E_\beta^+[g]}
	 	{ \int_{\mathbb R^+}|g|^2\,d\nu_{a}^+}, \qquad\text{for }g\in H^1(\mathbb{R}^+).
	 \end{equation}
\end{lemma}  
Notice that the numerator is independent of the measure $\nu$.

For eigenvalue comparison after, it is more convenient to use the area enclosed by Green levels as variables. For $t>0$ set 
\[
x_a(t):=|\{z\in\region_a:\psi_{ab}(z)>t\}|_g,
\qquad
x_b(t):=|\{z\in\region_b:-\psi_{ab}(z)>t\}|_g.
\]
Both functions are decreasing. Denote $m_a=|\Omega_a|$ and $m_b=|\Omega_b|$. Let
\[
\tau_a:(0,m_a]\longrightarrow[0,\infty),
\qquad
\tau_b:(0,m_b]\longrightarrow[0,\infty)
\]
denote the inverses of $x_a$ and $x_b$ respectively, and define the Green weights 
\begin{align}
G_a(x)
&:=\int_{\psi_{ab}=\tau_a(x)}\frac1{|d\psi_{ab}|_g}\,ds_g,
\qquad 0<x<m_a,
\label{eq:Robin-Ga}\\
G_b(x)
&:=\int_{\psi_{ab}=-\tau_b(x)}\frac1{|d\psi_{ab}|_g}\,ds_g,
\qquad 0<x<m_b.
\label{eq:Robin-Gb}
\end{align}
The coarea formula gives
\begin{equation*}
\tau_a'(x)=-\frac1{G_a(x)},
\qquad
\tau_b'(x)=-\frac1{G_b(x)}.
\end{equation*}
In particular,
\begin{equation}\label{eq:Robin-beta-integral}
\tau_j(x)=\int_x^{m_j}\frac{dy}{G_j(y)},
\qquad j=a,b.
\end{equation}
Now for a continuous representative of $g\in H^1(\mathbb{R})$, define
\begin{equation*}
f_a(x):=g(\tau_a(x)),
\qquad
f_b(x):=g(-\tau_b(x)).
\end{equation*}
Then the two functions have the same endpoint value
\begin{equation*}
f_a(m_a)=f_b(m_b)=g(0)=:c.
\end{equation*}
The change of variables gives
\begin{equation*}
    E_{\alpha}[g] = \displaystyle
\sum_{j=a,b}\int_0^{m_j}
\left[G_j|f_j'|^2+\frac{4\pi^2}{G_j}|f_j|^2\right]dx
+\alpha|c|^2.
\end{equation*}
The coarea formula gives
\begin{equation*}
    \int_{\mathbb{R}}\abs{g}^2 \,d\nu_{a,b}= \sum_{j=a,b}\int_{0}^{m_j}\abs{f_j}^2 \,dx.
\end{equation*}

\subsection{The variational problem.}\label{sec:vairation-def} 
We say  $G\in C^\infty(0,m)$  is an admissible weight on $(0,m)$ if 
$$
G > 0, \quad G(x) \sim 4\pi x \, \text{ as }x\to0, \quad \text{and}\quad \int_{x}^m \frac{1}{G(y)}\,dy < +\infty  \text{ for }x>0.
$$
For example, by \eqref{eq:Robin-beta-integral}, the Green weights $G_j$ defined in \eqref{eq:Robin-Ga}--\eqref{eq:Robin-Gb} are admissible. 
We shall denote by $ \mathcal{G}_m $ the space of all admissible weights on $(0,m)$. 

Let $G$ be an admissible weight  on $(0,m)$. Consider the space
\[
\mathcal F(G)
:=\left\{f\in H^1_{\rm loc}(0,m):
 \sqrt G\,f'\in L^2,
 f/\sqrt G\in L^2\right\}
\]
equipped with norm $\|{f}\|_{G}^2:= \|{\sqrt G f'}\|_{L^2}^2 + \|{f/\sqrt{G}}\|_{L^2}^2$. Let $\tau(x) =\int_{x}^m\frac{1}{G(y)}\,dy$, the change of variables $f(x)=g(\tau(x))$ gives an isomorphism between $\mathcal F(G)$ and $H^1(\mathbb R^+)$.  
 
{For nonzero finite positive Borel measures $\nu$ on $\mathbb{R}$ and $\nu^+$ on $\mathbb{R}^+$, define}
\begin{equation*} 
{\kappa_{\alpha}(\nu) := \inf_{\stackrel{g\in H^1(\mathbb{R}),}{  \int_{\mathbb{R}}\abs{g}^2d\nu \,>\,0}} \frac{E_{\alpha}[g]}{\int_{\mathbb{R}}\abs{g}^2d\nu},}
\end{equation*} 
and the corresponding variational problem on $\mathbb R^+$
\begin{equation*} 
{\kappa_{\beta}^+(\nu^+) := \inf_{\stackrel{g\in H^1(\mathbb{R}^+)}{\int_{\mathbb{R^+}}\abs{g}^2d\nu^+ > 0}} \frac{E_{\beta}^+[g]}{\int_{\mathbb{R^+}}\abs{g}^2d\nu^+}.}
\end{equation*}   

Extending $\nu^+$ by zero to $\mathbb R$, we have
\begin{equation}\label{eq:shift}
\kappa_{\alpha}(\nu^+) = \kappa_{\alpha + 2\pi}^+(\nu^+).
\end{equation}
Indeed, for every $c\in\mathbb R$, Lemma~\ref{lem:appendix-Robin-trace} gives
\[
 \inf_{\substack{h\in H^1(\mathbb R^-)\\ h(0)=c}}
 \int_{\mathbb R^-}(|h'|^2+4\pi^2|h|^2)\,dt
 =2\pi|c|^2,
\]
with the infimum attained uniquely by $h(t)=ce^{2\pi t}$. Minimizing over the negative half-line therefore yields the parameter shift in \eqref{eq:shift}. Moreover, if $g_*\in H^1(\mathbb R)$ is a minimizer of $\kappa_\alpha(\nu^+)$, then $g_*(t)=g_*(0)e^{2\pi t}$ for $t\leq0$, and $\restrict{g_*}{[0,\infty)}$ is a minimizer of $\kappa_{\alpha+2\pi}^+(\nu^+)$. 

For two weights $G_1\in\mathcal{G}_{m_1}, G_2\in \mathcal{G}_{m_2}$, define
\begin{equation*}
\widehat\kappa_\alpha(G_1,G_2)
:=\inf_{\substack{f_j\in\mathcal{F}(G_j) \\ f_1(m_1) = f_2(m_2)=:c}}
\frac{\displaystyle
\sum_{j=1,2}\int_0^{m_j}
\left[G_j|f_j'|^2+\frac{4\pi^2}{G_j}|f_j|^2\right]dx
+\alpha|c|^2}
{\displaystyle
\sum_{j=1,2}\int_0^{m_j}|f_j|^2\,dx}.\end{equation*}
Also, for a single weight $G\in \mathcal{G}_m$, define
\begin{equation*}
\widehat\kappa_\beta(G)
:=\inf_{0\neq f\in\mathcal F(G)}
\frac{\displaystyle
\int_0^m\left[G|f'|^2+\frac{4\pi^2}{G}|f|^2\right]dx
+\beta|f(m)|^2}
{\displaystyle\int_0^m|f|^2dx}.
\end{equation*}
The change of variables and isomorphism between $\mathcal F(G)$ and $H^1(\mathbb R^+)$ shows that, in the general case, with $\nu_{a,b}, G_a, G_b$ defined in  \eqref{eq:Robin-Ga}--\eqref{eq:Robin-Gb},
\begin{equation}\label{eq:Robin-kappa-identification-double}
\kappa_{\alpha}(\nu_{a,b}) = \widehat{\kappa}_\alpha(G_a, G_b).
\end{equation}
Also, in the boundary case, with $\nu_{a}^+, G_a$ defined by a single pole $a$,
\begin{equation}\label{eq:Robin-kappa-identification-single} 
{\kappa_\beta^+(\nu_a^+)} = \widehat{\kappa}_\beta(G_a).
\end{equation}

After changing to the $x$-variable, Appendix~\ref{sec:Appendix_Robin_minimizer} shows that, for $\alpha\geq -4\pi$, the variational infimum defining $\widehat\kappa_{\alpha}(G_a,G_b)$ is nonnegative and is attained by a unique normalized pair $(f_a,f_b)$ satisfying
\[
f_a>0,\qquad f_b>0,\qquad f_a(m_a)=f_b(m_b)=:c>0.
\]
Set $\kappa:=\widehat\kappa_{\alpha}(G_a,G_b)$. 
The corresponding Euler--Lagrange equations are
\begin{equation*}
-(G_jf_j')'+\frac{4\pi^2}{G_j}f_j
=\kappa f_j
\qquad\text{in }(0,m_j),\qquad j=a,b,
\end{equation*}
while the endpoint condition is
\begin{equation}\label{eq:Robin-Kirchhoff}
\lim_{x\uparrow m_a}G_a(x)f_a'(x)
+
\lim_{x\uparrow m_b}G_b(x)f_b'(x)
+\alpha c=0.
\end{equation}
Define the \emph{Riccati function}
\begin{equation*}
R_j(x):=G_j(x)\frac{f_j'(x)}{f_j(x)},
\qquad j=a,b.
\end{equation*}
By Lemma~\ref{lem:Robin-regular-end-flux}, each $R_j$ extends continuously to $m_j$.
Dividing \eqref{eq:Robin-Kirchhoff} by $c$ yields
\begin{equation}\label{eq:Robin-flux-balance}
R_a(m_a)+R_b(m_b)+\alpha=0.
\end{equation}
On the other hand, by Lemma~\ref{lem:Robin-singular-flux} in Appendix~\ref{sec:Appendix_Riccati}, we have 
\begin{equation}\label{eq:Rj'Rj0+}
R_j'
=
\frac{4\pi^2-R_j^2}{G_j}-\kappa,
\qquad
R_j(0+)=2\pi.
\end{equation}
Equation \eqref{eq:Robin-flux-balance}--\eqref{eq:Rj'Rj0+} is the starting point to estimate $\kappa$: we compare $R_j$ pointwise with the Riccati function in the constant-curvature case, then the monotonicity of the Riccati function in the spectral parameter $\kappa$ implies the required comparison. 

Similarly, consider the boundary one-pole case. When $\beta=\alpha+2\pi\geq -2\pi,$ the variational infimum $\widehat\kappa_{\beta}(G)$ is attained by a unique normalized positive minimizer $f$ satisfying
\[
-(Gf')'+\frac{4\pi^2}{G}f
=
\widehat\kappa_{\beta}(G) f
\qquad\text{in }(0,m),
\]
together with the endpoint condition
\[
\lim_{x\uparrow m}G(x)f'(x)+\beta f(m)=0.
\]
Define $R:=Gf'/f$, which extends continuously to $m$ by Lemma~\ref{lem:Robin-regular-end-flux}. We may rewrite the endpoint condition as
\[
R(m)+2\pi+\alpha=0.
\]
Note that in view of \eqref{eq:Rj'Rj0+}, one may regard this one pole endpoint condition as a special case of the two pole endpoint condition  \eqref{eq:Robin-flux-balance} with $m_b \to 0+$. 

For the standard constant-curvature model, set
$G_K(x):=4\pi x-Kx^2$. Then $G_K$ is positive on
\[
I_K=\begin{cases}
    (0, \frac{4\pi}{K}),&K>0,\\[1mm]
    (0,+\infty),&K\leq 0.
\end{cases}
\] 
{For $\mu\in\mathbb{R}$, let $f_{K,\mu}$ be a regular solution, chosen positive near $x=0$, of}
\begin{equation}\label{eq:model-area-EL}
-(G_Kf_{K,\mu}')'+\frac{4\pi^2}{G_K}f_{K,\mu}
=\mu f_{K,\mu},
\qquad
\lim_{x\downarrow0}G_K(x)\frac{f_{K,\mu}'(x)}{f_{K,\mu}(x)}=2\pi.
\end{equation}
Let $x_{\rm max} := \sup\,\{x:f_{K,\mu} > 0 \text{ in }(0,x)\}$. In $I_{K}\cap (0,x_{\rm max})$, define the constant-curvature Riccati function
\begin{equation}\label{eq:model-R-definition}
R_{K,\mu}(x):=G_K(x)\frac{f_{K,\mu}'(x)}{f_{K,\mu}(x)}.
\end{equation}
Throughout the paper, $R_{K,\mu}$ is considered only on this interval,
where $f_{K,\mu}>0$. Then
\begin{equation}\label{eq:model-riccati}
R_{K,\mu}'
=\frac{4\pi^2-R_{K,\mu}^2}{G_K}-\mu,
\qquad
R_{K,\mu}(0+)=2\pi.
\end{equation}
If $x_{\rm max}<\sup I_K$, then $f_{K,\mu}(x_{\rm max})=0$ and $f_{K,\mu}'(x_{\rm max})<0$. Hence $R_{K,\mu}(x_{\rm max}-)=-\infty$. 

We consider a tangential continuation of $R_{K,\mu}$. Let $x_{K}(\mu)\in I_K$ be the smallest positive root of $R_{K,\mu}(x) = -2\pi$. If no such root exists, just set $x_{K}(\mu) := \sup I_K$. Then define
\begin{equation}\label{eq:model-Rhat-definition}
\widehat R_{K,\mu}(x):=
\begin{cases}
R_{K,\mu}(x),&0 < x< x_K(\mu),\\[1mm]
-2\pi-\mu\bigl(x-x_K(\mu)\bigr),&x\geq x_K(\mu).
\end{cases}
\end{equation}
The function $\widehat R_{K,\mu}$ is well-defined in $I_K$. A direct computation shows that $\widehat R_{K,0} =R_{K,0}\equiv2\pi$. When $\mu > 0$, at $x=x_K(\mu)<\sup I_K$, \eqref{eq:model-riccati} gives
$R_{K,\mu}'(x_K(\mu))=-\mu<0$, so $\widehat R_{K,\mu}\in C^1(I_K)$. Moreover, $R_{K,\mu}'(x) \leq -\mu$ if $x\geq x_K(\mu)$, hence $R_{K,\mu}\leq \widehat R_{K,\mu}$. We set $\widehat R_{K,\mu}(0) := \widehat {R}_{K,\mu}(0+)=2\pi$.
 
For $m_1,m_2\geq0$, $m_1+m_2>0$, {with $K(m_1+m_2)<4\pi$ when $K>0$,} and $\alpha>-4\pi$, let $\mu_K^{\mathrm{cmp}}(m_1,m_2;\alpha)$ be the unique positive solution of
\begin{equation}\label{eq:model-comparison-root-early}
\widehat R_{K,\mu}(m_1)+\widehat R_{K,\mu}(m_2)+\alpha=0.
\end{equation}
The existence and uniqueness are proved in Subsection~\ref{sec:fixed-split-comparison}. At $\alpha=-4\pi$, since $\widehat R_{K,0} \equiv 2\pi$, we just set
\[
\mu_K^{\mathrm{cmp}}(m_1,m_2;-4\pi):=0.
\]

The Riccati function $R_{K,\mu}$ also relates the first angular Robin and
Helmholtz--Steklov eigenvalues.
Indeed, the regular first angular eigenfunction in geodesic polar
coordinates is $f_{K,\mu}(x(r))e^{i\theta}$. Its outward normal
derivative divided by its boundary value on $\partial D_K(x)$ is
\[
\ell_K(x)\frac{f_{K,\mu}'(x)}{f_{K,\mu}(x)}
=\frac{R_{K,\mu}(x)}{\ell_K(x)}.
\]
The Robin boundary condition therefore gives
\begin{equation}\label{eq:model-angular-boundary-early}
\lambda_{1,\mathrm{ang}}\bigl(D_K(m);\gamma\bigr)=\mu
\quad\Longleftrightarrow\quad
R_{K,\mu}(m)+\gamma\ell_K(m)=0,
\end{equation}
where the positivity of $f_{K,\mu}$ ensures that $\mu$ is the first angular
Robin eigenvalue. Similarly, the Helmholtz--Steklov boundary
condition in \eqref{eq:intro-SH-problem} with $\lambda=\mu$ gives
\begin{equation}\label{eq:model-angular-SH}
R_{K,\mu}(x)
=\ell_K(x)\,\sigma_{1,\mathrm{ang}}\bigl(D_K(x),\mu\bigr),
\end{equation}

\subsection{Proof of the main theorem}\label{sec:main-proof}
 
The proof uses the following three results, stated here in the order in which they enter. {Their proofs are given in Sections~\ref{sec:kappa}, \ref{sec:Robin-model-concavity}, and \ref{sec:Robin-orthogonality-test}, respectively.}

\begin{theorem}\label{thm:Robin-fixed-split-comparison}
Assume $K_g\leq K$ on $\region$, with $KA<4\pi$ when $K>0$, and let $-4\pi\leq\alpha\leq4\pi$. For two-pole Green weights $G_a,G_b$ with areas $m_a,m_b$,
\[
\widehat\kappa_\alpha(G_a,G_b)
\leq
\mu_K^{\rm cmp}(m_a,m_b;\alpha).
\]
For a one-pole Green weight $G$ with area $m$,
\begin{equation}\label{eq:Robin-fixed-split-one-pole}
\widehat\kappa_{\alpha+2\pi}(G)
\leq
\mu_K^{\rm cmp}(0,m;\alpha).
\end{equation}
\end{theorem}

\begin{proposition}\label{prop:model-split-optimization}
Let $m_1,m_2\geq0$ and $M:=m_1+m_2>0$, with $KM<4\pi$ when $K>0$. If $(K,M,\alpha)\in\mathcal P_1$, then
\begin{equation}\label{eq:model-equal-split}
\mu_K^{\rm cmp}(m_1,m_2;\alpha)
\leq
\mu_K^{\rm cmp}(M/2,M/2;\alpha).
\end{equation}
If $(K,M,\alpha)\in\mathcal P_2$, then
\begin{equation}\label{eq:model-endpoint-split}
\mu_K^{\rm cmp}(m_1,m_2;\alpha)
\leq
\mu_K^{\rm cmp}(0,M;\alpha).
\end{equation}
The inequality in \eqref{eq:model-endpoint-split} is strict when $-4\pi<\alpha<-4\pi+KM$ and $m_1,m_2>0$.
\end{proposition}

\begin{theorem}\label{thm:Robin-orthogonality-test} 
	{For $\alpha\geq-4\pi$, let} $\phi_1,\phi_2$ be real $L^2(\region)$-orthonormal Robin eigenfunctions associated with the first two eigenvalues, with $\phi_1>0$. Then at least one of the following holds:
	\begin{enumerate}[(1)]
		\item there exist distinct $a,b\in\region$ such that, with $g_{a,b}$ the positive normalized minimizer of ${E_\alpha[g]}/\|g\|_{L^2(\nu_{a,b})}^2$,
		\[
		\int_\region U_{a,b}[g_{a,b}]\phi_j\,dA_g=0,
		\qquad j=1,2;
		\]
		\item there exists $a\in\region$ such that, with $\beta=\alpha+2\pi$ and $g_a$ the positive normalized minimizer of ${E_\beta^+[g]}/\|g\|_{L^2(\nu_a^+)}^2$,
		\[
		\int_\region U_a[g_a]\phi_j\,dA_g=0,
		\qquad j=1,2.
		\]
	\end{enumerate}
\end{theorem}

\begin{proof}[Proof of the inequalities in Theorem~\ref{thm:intro-Robin-main}]
Let $\phi_1,\phi_2$ be the first two real $L^2(\region)$-orthonormal Robin eigenfunctions, with $\phi_1>0$.

\medskip
\noindent\textbf{Step 1: orthogonality.}
By Theorem~\ref{thm:Robin-orthogonality-test}, one of the following two cases occurs.

First, there exist distinct $a,b\in\region$ such that
$U_{a,b}[g_{a,b}]$ is orthogonal to $\phi_1$ and $\phi_2$. By the
min--max principle, \eqref{eq:Rayleigh-double}, and
\eqref{eq:Robin-kappa-identification-double},
\begin{align*}
\lambda_3\left(\region;\frac{\alpha}{L}\right)
&\leq
\frac{\mathfrak q_\alpha[U_{a,b}[g_{a,b}]]}
{\displaystyle\int_\region |U_{a,b}[g_{a,b}]|^2\,dA_g}\\
&=
\frac{E_\alpha[g_{a,b}]}
{\displaystyle\int_{\mathbb R}|g_{a,b}|^2\,d\nu_{a,b}}
=\kappa_\alpha(\nu_{a,b})
=\widehat\kappa_\alpha(G_a,G_b).
\end{align*}
In this case, set
\[
\kappa:=\widehat\kappa_\alpha(G_a,G_b),
\qquad
m_1:=|\region_a|_g,
\qquad
m_2:=|\region_b|_g.
\]

Second, there exists $a\in\region$ such that $U_a[g_a]$ is orthogonal
to $\phi_1$ and $\phi_2$. Since $\mathfrak q_\alpha[u]\leq\mathfrak q_{\alpha+2\pi}[u]$ for every $u\in H^1(\region)$, the min--max
principle, \eqref{eq:Rayleigh-single}, and
\eqref{eq:Robin-kappa-identification-single} give
\begin{align*}
\lambda_3\left(\region;\frac{\alpha}{L}\right)
&\leq
\frac{\mathfrak q_\alpha[U_a[g_a]]}
{\displaystyle\int_\region |U_a[g_a]|^2\,dA_g}\\
&\leq
\frac{\mathfrak q_{\alpha+2\pi}[U_a[g_a]]}
{\displaystyle\int_\region |U_a[g_a]|^2\,dA_g}\\
&=
\frac{E_{\alpha+2\pi}^+[g_a]}
{\displaystyle\int_{\mathbb R^+}|g_a|^2\,d\nu_a^+} 
={\kappa_{\alpha+2\pi}^+(\nu_a^+)}
=\widehat\kappa_{\alpha+2\pi}(G_a).
\end{align*}
In this case, set
\[
\kappa:=\widehat\kappa_{\alpha+2\pi}(G_a),
\qquad
m_1:=0,
\qquad
m_2:=A.
\]
Thus, in both cases,
\begin{equation}\label{eq:Robin-main-bridge}
\lambda_3\left(\region;\frac{\alpha}{L}\right)\leq\kappa,
\qquad
m_1+m_2=A.
\end{equation}

\medskip
\noindent\textbf{Step 2: comparison with fixed area split.}
Theorem~\ref{thm:Robin-fixed-split-comparison} gives
\begin{equation}\label{eq:Robin-main-fixed-split}
\kappa\leq\mu_K^{\rm cmp}(m_1,m_2;\alpha).
\end{equation}

\medskip
\noindent\textbf{Step 3: optimization of the split.}
Suppose first that $(K,A,\alpha)\in\mathcal P_1$. Proposition~\ref{prop:model-split-optimization} gives
\[
\mu_K^{\rm cmp}(m_1,m_2;\alpha)
\leq
\mu_K^{\rm cmp}(A/2,A/2;\alpha).
\]
Set $\nu:=\mu_K^{\rm cmp}(A/2,A/2;\alpha)$. By \eqref{eq:model-comparison-root-early},
\[
\widehat R_{K,\nu}(A/2)=-\frac{\alpha}{2}\geq-2\pi,
\]
so $\widehat R_{K,\nu}(A/2)=R_{K,\nu}(A/2)$.
By \eqref{eq:model-angular-boundary-early}, $\nu$ is the first angular
Robin eigenvalue of $D_K(A/2)$ with coefficient $\alpha/(2\ell_K(A/2))$.
By~\cite[Proposition~6.1]{langfordlaugesen2025robin}, this is also the
second Robin eigenvalue. Hence
\[
\mu_K^{\rm cmp}(A/2,A/2;\alpha)
=
\lambda_2\left(D_K(A/2);\frac{\alpha}{2\ell_K(A/2)}\right).
\]
Together with \eqref{eq:Robin-main-bridge} and \eqref{eq:Robin-main-fixed-split}, this proves the non-strict form of \eqref{eq:intro-Robin-main}.

Now suppose that $(K,A,\alpha)\in\mathcal P_2$. Proposition~\ref{prop:model-split-optimization} gives
\[
\mu_K^{\rm cmp}(m_1,m_2;\alpha)
\leq
\mu_K^{\rm cmp}(0,A;\alpha).
\]
Set $\nu:=\mu_K^{\rm cmp}(0,A;\alpha)$. Since $\alpha\leq-4\pi+KA<0$, \eqref{eq:model-comparison-root-early} gives
\[
\widehat R_{K,\nu}(A)=-(\alpha+2\pi)\geq-2\pi,
\]
so $\widehat R_{K,\nu}(A)=R_{K,\nu}(A)$. By \eqref{eq:model-angular-boundary-early},
\[
\mu_K^{\rm cmp}(0,A;\alpha)
=
\lambda_{1,\mathrm{ang}}\left(D_K(A);\frac{\alpha+2\pi}{\ell_K(A)}\right).
\]
Combining this with \eqref{eq:Robin-main-bridge} and \eqref{eq:Robin-main-fixed-split} proves the non-strict form of \eqref{eq:intro-Robin-one-disk}.

\medskip
\noindent\textbf{Strictness.}
In the one-pole case, Step~1 is strict since
\[
\mathfrak q_{\alpha+2\pi}[U_a[g_a]]
-\mathfrak q_\alpha[U_a[g_a]]=2\pi g_a(0)^2>0.
\]
In the two-pole case, equality would make $U_{a,b}[g_{a,b}]$ a
Robin eigenfunction, hence smooth in the interior. If $\alpha\ne0$,
the jump condition \eqref{eq:appendix-Robin-Kirchhoff-potential} gives
\[
g_{a,b}'(0+)-g_{a,b}'(0-)=\alpha g_{a,b}(0)\ne0,
\]
so $U_{a,b}[g_{a,b}]$ is not $C^1$ across a regular interior arc of
$\{\psi_{ab}=0\}$, a contradiction. This also covers $\alpha=-4\pi$,
where $g_{a,b}(t)=Ce^{-2\pi|t|}$.

If $\alpha=0$, then $\kappa>0$. Equality in Step~2 forces equality
in both endpoint Riccati comparisons, so
Lemma~\ref{lem:Robin-tangent-comparison} gives $G_a=G_K$ and $G_b=G_K$
on their respective intervals.
By Proposition~\ref{prop:Robin-Green-weight}, both Green regions
are model disks. This can happen only if $K>0$ and these model disks are complementary caps with total area $4\pi /K$, contradicts the condition $KA < 4\pi$. Thus both bounds are strict.
\end{proof}

This section proves the upper bounds in Theorem~\ref{thm:intro-Robin-main}. The sharpness statements, as well as the failure of the equal two-disk bound for $-4\pi<\alpha<-4\pi+KA$, are proved in Section~\ref{sec:sharpness}.

\section{Green-weight and Riccati Comparison}\label{sec:kappa}

This section contains the first analytic comparison. Its structure is parallel to classical comparison geometry:
\[
K_g\leq K
\Longrightarrow
G\geq G_K
\Longrightarrow
R\leq\widehat R_{K,\mu}
\Longrightarrow
\text{spectral comparison}.
\]
The first implication is geometric, the second is a one-dimensional Riccati comparison, and the last uses the endpoint condition of the minimizing one-dimensional problem.

\subsection{Comparison of Green weights}

Recall the constant-curvature model weight
\begin{equation*}
G_K(x):=4\pi x-Kx^2.
\end{equation*}
Set
\[
I_K:=
\begin{cases}
(0,\infty),&K\leq0,\\
(0,4\pi/K),&K>0.
\end{cases}
\]
Since $KA<4\pi$ when $K>0$, every area interval $(0,m)$ occurring below, with $0<m\leq A$, is contained in $I_K$.

The following Green function comparison was first proved by  Provenzano and Savo for spherical domains in \cite[Section~3]{provenzano2026isoperimetric}. In view of 
the isoperimetric estimate in~\cite[Isoperimetric inequality~(I)]{chavelfeldman1980},  for any simply-connected domain $\Omega$ with $K_g \le K$,  where $K$ is any real number, one also has 
\begin{equation}\label{eq:Green-level-isoperimetric}
	|\partial \Omega|_g^2 	\geq 4\pi|\Omega|_g-K|\Omega|_g^2,  	
\end{equation} 
with equality if and only if $\Omega$ is a geodesic disk with constant curvature $K$. 
Applying this to the Green superlevel set $E_x$ with $|E_x|_g=x$ for each $0<x<m=|\Omega|_g$, one can easily adapt Provenzano and Savo's proof to give the following more general result:    
\begin{proposition}[\cite{provenzano2026isoperimetric}]\label{prop:Robin-Green-weight}
Assume $K_g\leq K$ on $\region$. For every pole $a\in\region$, its Green weight $G_a$ satisfies
\begin{equation}\label{eq:Robin-Green-weight-comparison}
G_a(x)\geq G_K(x),\qquad0<x<m.
\end{equation}
Moreover, the equality holds for all $x \in (0,m)$ if and only if $K_g = K$ and $\Omega$ is a geodesic disk, with $a$ its center. 
\end{proposition}
 We remark that closely related one- and two-pole versions also appear in \cite[Sections~3.1 and~4.1]{michettiprovenzanosavo2026}.

\subsection{Riccati comparison}

Before the comparison, we prove the monotonicity of $\widehat R_{K,\mu}$ in $\mu$.

\begin{lemma}\label{lem:model-R-mu-monotonicity}
For every fixed $x\in I_K$, the map 
\[
\mu\longmapsto\widehat R_{K,\mu}(x)
\]
is continuous and strictly decreasing on $[0, {\infty})$. Moreover,
\begin{equation}\label{eq:model-Rhat-mu-minus-infty}
\widehat R_{K,\mu}(x)\longrightarrow-\infty
\qquad(\mu\to\infty).
\end{equation}
\end{lemma}

\begin{proof}
If $\widehat R_{K,\mu}(x)=R_{K,\mu}(x)$, normalize the regular solution by $f_{K,\mu}(x)=1$. Smooth dependence on $\mu$ and the Lagrange identity give
\begin{equation*}
\frac{\partial}{\partial\mu}R_{K,\mu}(x)
=-\int_0^x f_{K,\mu}(t)^2\,dt<0.
\end{equation*}
Suppose $x_K(\mu)$ lies in the interior of $I_K$. Since $\partial_xR_{K,\mu}(x_K(\mu))=-\mu$, the implicit function theorem shows that $x_K$ is differentiable near $\mu$ and
\[
x_K'(\mu)=\frac{\partial_\mu R_{K,\mu}(x_K(\mu))}{\mu}<0.
\]
For $x>x_K(\mu)$, \eqref{eq:model-Rhat-definition} gives
\[
\frac{\partial}{\partial\mu}\widehat R_{K,\mu}(x)
=-(x-x_K(\mu))+\mu x_K'(\mu)<0.
\]
The two derivative formulas agree at $x=x_K(\mu)$, and the continuous dependence gives continuity when $x=x_K(\mu)$.

It remains to prove \eqref{eq:model-Rhat-mu-minus-infty}. Fix $x\in I_K$ and choose $x_0\in (0,x)$. Let $\lambda_D(x_0)$ be the first Dirichlet eigenvalue of the one-dimensional problem
$$
(-G_Kf')'+\frac{4\pi^2}{G_K}f=\lambda_D f
$$
with $f\in\mathcal{F}(G_K)$ and the Dirichlet condition $f(x_0)=0$. Let $u_D$ be the first Dirichlet eigenfunction. We claim that whenever $\mu>\lambda_D(x_0)$, the regular solution $f_{K,\mu}$ has a zero in $(0,x_0)$. Indeed, if $f_{K,\mu}>0$ in $(0,x_0)$, integration by parts gives
\[
\left[
G_K\bigl(u_D f_{K,\mu}'-f_{K,\mu}u_D'\bigr)
\right]_{0}^{x_0}
=
\bigl(\lambda_D(x_0)-\mu\bigr)
\int_0^{x_0} u_D f_{K,\mu}\,dt<0.
\]
Note that by Lemma \ref{lem:Robin-singular-flux}, $G_K(0+)f'_{K,\mu}(0+)=G_K(0+)u'_D(0+)=0$. 
Together with $u_D(0+)=u_D(x_0)=f_{K,\mu}(0)=0$, the left-hand side equals
\[
-G_K(x_0)f_{K,\mu}(x_0)u_D'(x_0)\geq0,
\]
a contradiction. Thus $f_{K,\mu}$ has a first zero at some point $x_1\leq x_0$. Since $R_{K,\mu}(t)\to-\infty$ as $t\uparrow x_1$, the value $-2\pi$ is attained before $x_1$. Hence
\[
x_K(\mu)<x_0 < x.
\]
Therefore, for every $\mu>\lambda_D(x_0)$, the second line of \eqref{eq:model-Rhat-definition} applies at $x$ and yields
\[
\widehat R_{K,\mu}(x)
=-2\pi-\mu\bigl(x-x_K(\mu)\bigr)\leq -\mu(x-x_0).
\]
Letting $\mu\to\infty$ proves \eqref{eq:model-Rhat-mu-minus-infty}.
\end{proof}

\begin{lemma}\label{lem:Robin-tangent-comparison}
Let $m\in I_K$, $G\in C^\infty(0,m)$ satisfying $G\geq G_K$, and let $\mu>0$. Suppose that $R\in C^1(0,m)\cap C(0,m]$ solves
\begin{equation}\label{eq:general-Riccati-comparison}
R'
=
\frac{4\pi^2-R^2}{G}-\mu,
\qquad
\lim_{x\downarrow0}R(x)=2\pi.
\end{equation}
Then
\begin{equation}\label{eq:R-vs-model-R}
R(x)\leq\widehat R_{K,\mu}(x),
\qquad 0<x\leq m.
\end{equation}
Moreover, if $R(m) = \widehat R_{K,\mu}(m)$, then $G \equiv G_K$ on $(0, m)$.
\end{lemma}
\begin{proof}
Fix $x\in(0,m)$. First suppose that $\widehat R_{K,\mu}(x)\geq-2\pi$. By definition, $\widehat R_{K,\mu}(t)=R_{K,\mu}(t)\geq-2\pi$ for $t\in(0,x]$. Since $R_{K,\mu}(t)<R_{K,0}(t)=2\pi$ by Lemma \ref{lem:model-R-mu-monotonicity} , it follows that $|R_{K,\mu}(t)|\leq2\pi$ for $t\in(0,x]$. Let $w:=R-R_{K,\mu}$. Subtracting the two Riccati equations gives
\begin{equation*}
w'
=-\frac{R+R_{K,\mu}}{G}\,w
+\bigl(4\pi^2-R_{K,\mu}^2\bigr)
\left(\frac1G-\frac1{G_K}\right).
\end{equation*}
The second term is nonpositive in $(0,x]$ because $|R_{K,\mu}|\leq2\pi$ and $G\geq G_K$. Write
\[
a(t):=\frac{R(t)+R_{K,\mu}(t)}{G(t)},
\qquad
b(t):=\bigl(4\pi^2-R_{K,\mu}(t)^2\bigr)
\left(\frac1{G(t)}-\frac1{G_K(t)}\right)\leq0.
\]
Then $w'+aw=b$. For $0<\varepsilon<x$, the Duhamel formula on $[\varepsilon,x]$ gives
\begin{equation*}
w(x)
=
\exp\left(-\int_\varepsilon^x a(s)\,ds\right)w(\varepsilon)
+
\int_\varepsilon^x
\exp\left(-\int_t^x a(s)\,ds\right)b(t)\,dt.
\end{equation*}
Since $b\leq0$, we obtain
\begin{equation}\label{eq:Duhamel-ineq}
w(x)
\leq
w(\varepsilon)
\exp\left(-\int_\varepsilon^x a(s)\,ds\right).
\end{equation}
Since $a(t)>0$ for all sufficiently small $t$ and is continuous on every compact subinterval of $(0,m)$, it is bounded below on $(0,x]$ by some $C\in\mathbb R$. Hence
\[
0\leq \exp\left(-\int_\varepsilon^x a(s)\,ds\right)\leq e^{|C|x}
\qquad\text{for }0<\varepsilon<x.
\]
Moreover, $w(\varepsilon)\to0$. Letting $\varepsilon\downarrow0$ gives $w(x)\leq0$. Therefore
\[
R(x)\leq R_{K,\mu}(x)=\widehat R_{K,\mu}(x)
\]
as long as $\widehat R_{K,\mu}(x)\geq-2\pi$.

If $\widehat R_{K,\mu}(x)<-2\pi$, then $x_K(\mu)<x$, and the preceding inequality at $x_K(\mu)$ gives
\[
R(x_K(\mu))\leq\widehat R_{K,\mu}(x_K(\mu))=-2\pi.
\]
Whenever $R\leq-2\pi$, \eqref{eq:general-Riccati-comparison} gives
\[
R'\leq-\mu<0.
\]
Hence $R\leq-2\pi$ and $R'\leq-\mu$ for $x\geq x_K(\mu)$. Therefore
\[
R(x)
\leq-2\pi-\mu\bigl(x-x_K(\mu)\bigr)
=\widehat R_{K,\mu}(x).
\]
Together with the first part, this proves \eqref{eq:R-vs-model-R} for $x<m$.
Letting $x\uparrow m$ gives the endpoint inequality.

Suppose equality holds at $m$. If $m>x_K(\mu)$, then $R<-2\pi$ and
$(R-\widehat R_{K,\mu})'=(4\pi^2-R^2)/G<0$ on $(x_K(\mu),m)$,
contradicting endpoint equality. Thus
$\widehat R_{K,\mu}=R_{K,\mu}$ on $(0,m]$.
For each $\varepsilon\in(0,m)$, letting $x\uparrow m$ in \eqref{eq:Duhamel-ineq} gives
\[
0=w(m)\leq w(\varepsilon)\exp\left(-\int_{\varepsilon}^m a(t)\,dt\right)\leq0.
\]
The exponential factor is positive, so $w(\varepsilon)=0$. Hence $w\equiv0$
and $b\equiv0$ on $(0,m)$; since $|R_{K,\mu}|<2\pi$ there, we obtain $G\equiv G_K$.
\end{proof}

\subsection{The comparison for a fixed split}\label{sec:fixed-split-comparison}
 
{For $m_1,m_2\geq0$ with $m_1+m_2>0$, suppose $K(m_1+m_2)<4\pi$ when $K>0$. We first justify the definition of $\mu_K^{\rm cmp}$.} For $\alpha>-4\pi$, the left-hand side of \eqref{eq:model-comparison-root-early} equals $4\pi+\alpha>0$ at $\mu=0$. Since $m_1+m_2>0$, Lemma~\ref{lem:model-R-mu-monotonicity} shows that it is strictly decreasing and tends to $-\infty$. Hence \eqref{eq:model-comparison-root-early} has a unique positive solution.

\begin{proof}[Proof of Theorem~\ref{thm:Robin-fixed-split-comparison}]
At $\alpha=-4\pi$ both sides are zero, so assume $\alpha>-4\pi$.

In the two-pole case set $\kappa:=\widehat\kappa_\alpha(G_a,G_b)$. By \eqref{eq:Robin-flux-balance},
\[
R_a(m_a)+R_b(m_b)+\alpha=0.
\]
By Lemma~\ref{lem:Robin-regular-end-flux}, $R_j$ extends continuously to $m_j$.
Proposition~\ref{prop:Robin-Green-weight} and Lemma~\ref{lem:Robin-tangent-comparison} give
\[
0
\leq
\widehat R_{K,\kappa}(m_a)+\widehat R_{K,\kappa}(m_b)+\alpha.
\]
Since the left-hand side of \eqref{eq:model-comparison-root-early} decreases strictly in $\mu$ and vanishes when $\mu=\mu_K^{\rm cmp}(m_a,m_b;\alpha)$,
\[
\kappa\leq\mu_K^{\rm cmp}(m_a,m_b;\alpha).
\]

In the one-pole case set $\kappa:=\widehat\kappa_{\alpha+2\pi}(G)$. The endpoint condition $R(m)+\alpha+2\pi=0$, together with $\widehat R_{K,\kappa}(0)=2\pi$, gives
\[
0
\leq
\widehat R_{K,\kappa}(0)+\widehat R_{K,\kappa}(m)+\alpha,
\]
and the same argument yields \eqref{eq:Robin-fixed-split-one-pole}.
\end{proof}

\section{Optimization of the Model Area Split}\label{sec:Robin-model-concavity}

Section~\ref{sec:kappa} reduces the comparison to the number $\mu_K^{\rm cmp}(m_1,m_2;\alpha)$ with $m_1+m_2=A$. The optimizing split depends on whether the corresponding spectral parameter lies below or above $2K$.

\subsection{Convexity and concavity of the model Riccati function}
Recall
\[
I_K:=
\begin{cases}
(0,\infty),&K\leq0,\\
(0,4\pi/K),&K>0.
\end{cases}
\]
  
Let
\[
s_K(r):=
\begin{cases}
\sin(\sqrt K r)/\sqrt K,&K>0,\\
r,&K=0,\\
\sinh(\sqrt{-K}\,r)/\sqrt{-K},&K<0.
\end{cases}
\] 
\begin{lemma}\label{lem:model-R-critical-identity}
For $x\in I_K\cap(0,x_{\rm max})$, if we denote $B(x):=\ell_K(x)f_{K,\mu}(x)$, then 
\begin{equation}\label{eq:model-R-critical-identity}
R_{K,\mu}(x)=2\pi-Kx
+\frac{2K-\mu}{B(x)}\displaystyle\int_0^x B(y)\,dy.
\end{equation}
\end{lemma}

\begin{proof}
Consider geodesic polar coordinates $(r,\theta)$ on
$D:=D_K(x)$, and let
\[
u_\mu(r,\theta)=f_{K,\mu}(x(r))e^{i\theta},
\qquad u_*(r,\theta)=s_K(r)e^{i\theta}.
\]
Then
\[
-\Delta u_\mu=\mu u_\mu,
\qquad
\partial_\nu u_\mu=\frac{R_{K,\mu}(x)}{\ell_K(x)}u_\mu,
\]
and
\[
-\Delta u_*=2Ku_*,
\qquad
\partial_\nu u_*=\frac{2\pi-Kx}{\ell_K(x)}u_*.
\]
In particular, for $\mu=2K$ we have
$R_{K,2K}(x)=2\pi-Kx$.
For $\mu\ne2K$, Green's identity gives
\[
\frac{R_{K,\mu}(x)-(2\pi-Kx)}{2K-\mu}
=\ell_K(x)
\frac{\displaystyle\int_D u_\mu\overline{u_*}\,dA}
{\displaystyle\int_{\partial D}u_\mu\overline{u_*}\,ds}.
\]
Since $\ell_K(x(r))=2\pi s_K(r)$ and
$dx/dr=\ell_K(x)$, the right-hand side equals the quotient $\frac{1}{B(x)}\int_0^x B(y)\,dy$. 
\end{proof}

\begin{lemma}\label{lem:Sturm-integral-ratio-convexity}
Let $p\in C^1(0,m)$ be positive and satisfy
$(x^2p(x))'>0$. Suppose that
$u\in C^1([0,m))\cap C^2((0,m))$ satisfies
\[
u''+p(x)u=0,
\qquad u(0)=0,\quad u'(0)>0,\quad u>0\ \text{on }(0,m).
\]
Then $\displaystyle\frac{1}{u(x)}\int_0^x u(y)\,dy$ is strictly convex on $(0,m)$.
\end{lemma}

\begin{proof}
Since $u''=-pu<0$ and $u(0)=0$,
\begin{equation}\label{eq:Sturm-dilation-tangent}
u(x)-xu'(x)>0.
\end{equation}
For $0<t<1$, let $u_t(x):=u(tx)$. The equations at $x$
and $tx$ give
\begin{equation}\label{eq:Sturm-dilation-Wronskian}
(uu_t'-u'u_t)'=
\bigl(p(x)-t^2p(tx)\bigr)uu_t.
\end{equation}
Since $(x^2p(x))'>0$, we have $p(x)>t^2p(tx)$.
Thus $uu_t'-u'u_t$ is strictly increasing and vanishes
at $x=0$, so
\begin{equation}\label{eq:Sturm-dilation-ratio}
\left(\frac{u_t}{u}\right)'>0.
\end{equation}
Using the same equations, a direct calculation gives
\[
\left(x\frac{u_t}{u}\right)''
=2\left(1-\frac{xu'}u\right)
\left(\frac{u_t}{u}\right)'
+x\bigl(p(x)-t^2p(tx)\bigr)\frac{u_t}{u}>0,
\]
by \eqref{eq:Sturm-dilation-tangent} and
\eqref{eq:Sturm-dilation-ratio}. Hence
\[
\frac{\int_0^x u(y)\,dy}{u(x)}
=\int_0^1 x\frac{u(tx)}{u(x)}\,dt
\]
is strictly convex.
\end{proof}

\begin{proposition}\label{prop:model-area-concavity}
Let $\mu>0$.
\begin{enumerate}[(1)]
\item If $K>0$ and $0<\mu<2K$, then $x_{K}(\mu) \geq 4\pi/K$ and $\widehat R_{K,\mu} = R_{K,\mu}$ is strictly convex on $I_K$.
\item If $K>0$ and $\mu=2K$, then
\begin{equation}\label{eq:lambda-2K-linear}
R_{K,2K}(x)=2\pi-Kx\quad(x\in I_K).
\end{equation}
\item If $K>0$ and $\mu>2K$, then $\widehat R_{K,\mu}$ is concave on $I_K$.
\item If $K\leq0$, then $\widehat R_{K,\mu}$ is concave on $(0,\infty)$.
\end{enumerate}
\end{proposition}

\begin{proof}
Since $G_K(x)=\ell_K(x)^2=4\pi x-Kx^2$,
\[
\left(\frac{\mu x^2}{G_K(x)}\right)'
=\frac{4\pi\mu}{(4\pi-Kx)^2}>0.
\]
A direct calculation from \eqref{eq:model-area-EL} shows that the
function $B$ in Lemma~\ref{lem:model-R-critical-identity} satisfies
\[
B''+\frac{\mu}{G_K}B=0,\qquad B(0)=0,\quad B'(0)>0
\]
after positive normalization. Lemma~\ref{lem:Sturm-integral-ratio-convexity}
therefore gives
\[
\left(\frac{\int_0^x B(y)\,dy}{B(x)}\right)''>0
\]
wherever $f_{K,\mu}>0$. By Lemma~\ref{lem:model-R-critical-identity},
$R_{K,\mu}''$ consequently has the strict sign of $2K-\mu$
on this interval, unless $\mu=2K$, when
$R_{K,2K}(x)=2\pi-Kx$.

Since $\widehat R_{K,\mu}$ is a tangential continuation of $R_{K,\mu}$, it inherits the non-strict convexity or concavity of $R_{K,\mu}$. This proves (3) and (4). 

For (1), if $K>0$ and $0<\mu<2K$, Lemma~\ref{lem:model-R-mu-monotonicity}
gives
\[
\widehat R_{K,\mu}(x)>
\widehat R_{K,2K}(x)=2\pi-Kx>-2\pi,
\qquad x\in I_K.
\]
Thus the tangential continuation is never entered, that is, $x_K(\mu)\geq 4\pi/K$ and $\widehat R_{K,\mu} = R_{K,\mu}$. This proves (1).
\end{proof}
For later use, the regular-end expansion of \eqref{eq:model-riccati} is
\begin{equation}\label{eq:model-R-expansion}
R_{K,\mu}(x)=2\pi-\frac{\mu}{2}x
+\frac{\mu(2K-\mu)}{48\pi}x^2+O(x^3).
\end{equation}

\subsection{Optimization of the split}

\begin{proof}[Proof of Proposition~\ref{prop:model-split-optimization}]
At $\alpha=-4\pi$, $\mu_K^{\rm cmp}=0$ for every split. If $K>0$ and $\alpha=-4\pi+KM$, then \eqref{eq:lambda-2K-linear} gives
\[
\widehat R_{K,2K}(m_1)+\widehat R_{K,2K}(m_2)+\alpha=0,
\]
so uniqueness yields $\mu_K^{\rm cmp}(m_1,m_2;\alpha)=2K$ for every split. Thus assume $\alpha>-4\pi$, and in the $\mathcal P_2$ case also $\alpha<-4\pi+KM$.

Set
\[
\mu:=\mu_K^{\rm cmp}(m_1,m_2;\alpha).
\]
Then
\begin{equation}\label{eq:model-comparison-root-at-mu}
\widehat R_{K,\mu}(m_1)+\widehat R_{K,\mu}(m_2)+\alpha=0.
\end{equation}

Suppose $(K,M,\alpha)\in\mathcal P_1$. If $K>0$, then $\alpha\geq-4\pi+KM$, and \eqref{eq:lambda-2K-linear} gives
\[
\widehat R_{K,2K}(m_1)+\widehat R_{K,2K}(m_2)+\alpha
=4\pi-KM+\alpha\geq0.
\]
Hence $\mu\geq2K$. Proposition~\ref{prop:model-area-concavity} gives the required concavity both for $K>0$ and for $K\leq0$, and therefore
\[
2\widehat R_{K,\mu}(M/2)+\alpha\geq0.
\]
Since $\mu\mapsto\widehat R_{K,\mu}(M/2)$ is strictly decreasing, this proves \eqref{eq:model-equal-split}.

Now suppose $(K,M,\alpha)\in\mathcal P_2$. Then $K>0$ and $-4\pi<\alpha<-4\pi+KM$. At $\mu=0$ the left-hand side of \eqref{eq:model-comparison-root-early} is positive, while at $\mu=2K$ it equals
\[
4\pi-KM+\alpha<0.
\]
Thus $0<\mu<2K$. Proposition~\ref{prop:model-area-concavity} gives
\[
R_{K,\mu}(m_1)+R_{K,\mu}(m_2)
\leq
R_{K,\mu}(0)+R_{K,\mu}(M)
=2\pi+R_{K,\mu}(M),
\]
with strict inequality when $m_1,m_2>0$. Since $R_{K,\mu}>R_{K,2K}=2\pi-Kx>-2\pi$ on $(0,M]$, the hats in \eqref{eq:model-comparison-root-at-mu} are unnecessary. Hence
\[
2\pi+R_{K,\mu}(M)+\alpha\geq0,
\]
and strict decrease in $\mu$ proves \eqref{eq:model-endpoint-split}, with the stated strictness.
\end{proof}

\section{A topological orthogonality  theorem}\label{sec:Robin-orthogonality-test}

We first give two elementary convergence facts. They will be used to show that the minimizing test function depends continuously on the pole pair, including the cases that a pole reaches the boundary and that the two poles collapse. For the meaning of $\nu_n\rightharpoonup\nu$ in the following lemma, see Proposition~\ref{prop:measure-eigenpair-stability}.

\begin{lemma}\label{lem:orthogonality-pushforward}
Let $\psi_n,\psi:\region\to\mathbb R$ be measurable and finite almost everywhere, and assume that $\psi_n\to\psi$ almost everywhere. If
\[
\nu_n:=(\psi_n)_\#(dA_g),
\qquad
\nu:=\psi_\#(dA_g),
\]
then $\nu_n\rightharpoonup\nu$ on $\mathbb R$.
\end{lemma}

\begin{proof}
For every $h\in C_b(\mathbb R)$,
\[
\int_{\mathbb R}h\,d\nu_n
=
\int_\region h(\psi_n)\,dA_g
\longrightarrow
\int_\region h(\psi)\,dA_g
=
\int_{\mathbb R}h\,d\nu
\]
by dominated convergence.
\end{proof}

\begin{lemma}\label{lem:orthogonality-composition}
Assume $\psi_n\to\psi$ almost everywhere in $\region$, $\chi_n\to\chi$ almost everywhere in $\region$, and $|\chi_n|\leq1$. If $g_n\to g$ strongly in $H^1(\mathbb R)$, then
\[
 g_n(\psi_n)\chi_n
 \longrightarrow
 g(\psi)\chi
 \qquad\text{strongly in }L^2(\region).
\]
\end{lemma}

\begin{proof}
The embedding $H^1(\mathbb R)\subset C_b(\mathbb R)$ gives $g_n\to g$ uniformly on $\mathbb R$. Hence
\begin{align*}
\|g_n(\psi_n)\chi_n-g(\psi)\chi\|_{L^2}
&\leq
\|(g_n-g)(\psi_n)\chi_n\|_{L^2}\\
&\quad+
\|(g(\psi_n)-g(\psi))\chi_n\|_{L^2}\\
&\quad+
\|g(\psi)(\chi_n-\chi)\|_{L^2}.
\end{align*}
The first term is at most $A^{1/2}\|g_n-g\|_{L^\infty}$ and tends to zero. The other two terms tend to zero by dominated convergence.
\end{proof}
 
{For the rest of this section fix $\alpha\geq-4\pi$.} We use the natural almost-everywhere extension of the Green potentials to boundary poles: $\psi_a\equiv0$ when $a\in\partial\region$, and
\[
 \psi_{ab}:=\psi_a-\psi_b,
 \qquad a,b\in\overline\region.
\]
Set
\[
 \nu_{a,b}:=(\psi_{ab})_\#(dA_g).
\]
For every $(a,b)\in\overline\region\times\overline\region$, let $g_{a,b}$ be the positive normalized minimizer of $\kappa_\alpha(\nu_{a,b})$.
For $\alpha>-4\pi$, this minimizer is given by Proposition~\ref{prop:appendix-Robin-minimizer}. For $\alpha=-4\pi$, it is the corresponding normalization of $e^{-2\pi|t|}$.

For $b\in\region$, also set
\[
 \nu_b^+:=(\psi_b)_\#(dA_g),
\]  
and let $\widehat g_b$ be the positive normalized minimizer of $\kappa_\alpha(\nu_b^+)$. By \eqref{eq:shift} and the argument following it, with $\beta:=\alpha+2\pi$, the restriction $g_b:=\widehat g_b|_{[0,\infty)}$ is the positive normalized minimizer of $\kappa_\beta^+(\nu_b^+)$. 

For a real function $w\in L^2(\region)$ define
\begin{equation*}
 F_w(a,b)
 :=\int_\region U_{a,b}[g_{a,b}]\,w\,dA_g,
 \qquad
 (a,b)\in\overline\region\times\overline\region.
\end{equation*}
The phase function $e^{i\Theta_{ab}}$ in $U_{a,b}$ (see \eqref{eq:Robin-dipole-trial} for its definition) is determined by a fixed uniformization $\Phi: \region \to \disk$ as in Proposition~\ref{prop:one-pole-properties}.

\begin{proposition}\label{prop:orthogonality-map}
For every real $w\in L^2(\region)$, the function $F_w$ has the following properties:
\begin{enumerate}[(1)]
\item $F_w$ is continuous on $\overline\region\times\overline\region$;
\item $F_w(b,a)=\overline{F_w(a,b)}$ for any $a,b\in\overline\region$;
\item there is a continuous map 
 $G_w:\overline\region\to\mathbb C$
such that
\begin{equation}\label{eq:orthogonality-boundary-factorization}
 F_w(a,b)=-\Phi(a)\overline{G_w(b)},
 \qquad a\in\partial\region,\quad b\in\overline\region;
\end{equation}
\end{enumerate}
\end{proposition}

\begin{proof}
We first prove (1). Let $(a_n,b_n)\to(a,b)$ in $\overline\region\times\overline\region$, and write
\[
 \psi_n:=\psi_{a_nb_n},
 \qquad
 \psi:=\psi_{ab},
 \qquad
 \chi_n:=e^{i\Theta_{a_nb_n}},
 \qquad
 \chi:=e^{i\Theta_{ab}}.
\]
By \eqref{eq:psi}, $\psi_n\to\psi$ almost everywhere in $\region$.
Proposition~\ref{prop:one-pole-properties} and \eqref{eq:extended-two-pole-phase} similarly give
$\chi_n\to\chi$ almost everywhere in $\region$ and $|\chi_n|=1$.
Lemma~\ref{lem:orthogonality-pushforward} therefore gives
\begin{equation}\label{eq:orthogonality-measure-convergence}
 \nu_{a_n,b_n}\rightharpoonup\nu_{a,b}.
\end{equation}
If $\alpha>-4\pi$, Proposition~\ref{prop:measure-eigenpair-stability} yields
\[
 g_{a_n,b_n}\longrightarrow g_{a,b}
 \qquad\text{strongly in }H^1(\mathbb R).
\]
When $\alpha=-4\pi$, write
\[
 g_{a_n,b_n}(t)=c_n e^{-2\pi|t|},
 \qquad
 g_{a,b}(t)=c e^{-2\pi|t|}.
\]
The normalization gives
\[
 c_n^{-2}=\int_{\mathbb R}e^{-4\pi|t|}\,d\nu_{a_n,b_n}(t),
 \qquad
 c^{-2}=\int_{\mathbb R}e^{-4\pi|t|}\,d\nu_{a,b}(t),
\]
so \eqref{eq:orthogonality-measure-convergence} gives $c_n\to c$, and again $g_{a_n,b_n}\to g_{a,b}$ in $H^1(\mathbb R)$. Lemma~\ref{lem:orthogonality-composition} now gives
\[
 U_{a_n,b_n}[g_{a_n,b_n}]
 \longrightarrow
 U_{a,b}[g_{a,b}]
 \qquad\text{strongly in }L^2(\region).
\]
Hence $F_w(a_n,b_n)\to F_w(a,b)$ by Cauchy--Schwarz.

For (2), let $R(t):=-t$. Then
\[
 \psi_{ba}=-\psi_{ab},
 \qquad
 \nu_{b,a}=R_\#\nu_{a,b},
 \qquad
 e^{i\Theta_{ba}}=\overline{e^{i\Theta_{ab}}}.
\]
Since $E_\alpha[g\circ R]=E_\alpha[g]$, uniqueness of the positive normalized minimizer gives
\[
 g_{b,a}(t)=g_{a,b}(-t).
\]
Thus
\[
 U_{b,a}[g_{b,a}]
 =\overline{U_{a,b}[g_{a,b}]},
\]
and (2) follows because $w$ is real.

We next prove (3). 
First suppose that $a\in\partial\region$ and $b\in\region$. Then $\psi_{ab}=-\psi_b$ and
\[
 \nu_{a,b}=R_\#\nu_b^+,
\]
which depends only on $b$. Hence
\[
 g_{a,b}(t)=\widehat g_b(-t).
\]
Moreover, by \eqref{eq:extended-two-pole-phase},
\[
 e^{i\Theta_{ab}}
 =-\Phi(a)\,\overline{e^{i\Theta_b}}.
\]
It follows that
\[
 U_{a,b}[g_{a,b}]
 =-\Phi(a)\,\overline{U_b[g_b]}.
\]
Define, for $b\in\region$,
\begin{equation}\label{eq:orthogonality-G-interior}
 G_w(b):=\int_\region U_b[g_b]w\,dA_g.
\end{equation}
Then \eqref{eq:orthogonality-boundary-factorization} holds in the boundary case.

It remains to consider $a,b\in\partial\region$, which is the collapsed case. Here $\nu_{a,b}=A\delta_0$, and the normalization gives $g_{a,b}(0)=A^{-1/2}$. Also
\[
 e^{i\Theta_{ab}}=\Phi(a)\overline{\Phi(b)}.
\]
Therefore
\[
 F_w(a,b)
 =A^{-1/2}\Phi(a)\overline{\Phi(b)}
 \int_\region w\,dA_g.
\]
Define
\begin{equation*}
 G_w(b):=-A^{-1/2}\Phi(b)\int_\region w\,dA_g,
 \qquad b\in\partial\region.
\end{equation*}
Then \eqref{eq:orthogonality-boundary-factorization} also holds in the collapsed case. The continuity of $G_w$ on $\overline\region$ follows from (1) and \eqref{eq:orthogonality-boundary-factorization}: for any fixed $a_0\in\partial\region$,
\[
 G_w(b)=\overline{\frac{F_w(a_0,b)}{-\Phi(a_0)}}.
\]
This completes the proof. 
\end{proof}

In particular, after transporting the pole variables to $\overline\disk\times\overline\disk$ by $\Phi\times\Phi$, Proposition~\ref{prop:orthogonality-map} shows that $F_w$ belongs to the class $\mathcal W_1$ defined before Theorem~\ref{thm:topo}.

Since $\region$ is a connected bounded Lipschitz domain, the first Robin eigenvalue is simple and a corresponding eigenfunction may be chosen strictly positive in $\region$ for every real Robin parameter; see \cite[Proposition~2.6 and the discussion in Section~1]{daners2013principal}.

\begin{proof}[Proof of Theorem~\ref{thm:Robin-orthogonality-test}]
Apply Proposition~\ref{prop:orthogonality-map} with $w=\phi_1$ and $w=\phi_2$. After uniformization, the two functions belong to $\mathcal W_1$. With $J= \restrict{-\Phi \circ \Phi^{-1}}{\mathbb{S}^1}$ we have $\degmap J=1$. Thus Theorem~\ref{thm:topo} gives a common zero $(a,b)\in\overline\region\times\overline\region$.

The zero cannot be in the collapsed case. Indeed, if $a=b$, then $\nu_{a,a}=A\delta_0$, $g_{a,a}(0)=A^{-1/2}$, and $U_{a,a}[g_{a,a}]\equiv A^{-1/2}$. Hence
\[
 F_{\phi_1}(a,a)
 =A^{-1/2}\int_\region\phi_1\,dA_g>0.
\]
If $a,b\in\partial\region$, then
\[
 F_{\phi_1}(a,b)
 =A^{-1/2}\Phi(a)\overline{\Phi(b)}
 \int_\region\phi_1\,dA_g\neq0.
\]
Thus only the general and boundary cases remain.

If $a,b\in\region$, then $a\neq b$ and the common-zero condition is exactly alternative~(1). If one pole lies on the boundary, symmetry allows us to assume $a\in\partial\region$ and $b\in\region$. By \eqref{eq:orthogonality-boundary-factorization},
\[
 0=F_{\phi_j}(a,b)
 =J(a)\overline{G_{\phi_j}(b)},
 \qquad j=1,2.
\]
Since $J(a)\neq0$, we have $G_{\phi_j}(b)=0$. By \eqref{eq:orthogonality-G-interior}, this is
\[
 \int_\region U_b[g_b]\phi_j\,dA_g=0,
 \qquad j=1,2,
\]
which is alternative~(2).
\end{proof}

\section{Proof of Theorem~\ref{thm:topo}}\label{sec:topo}

We first prove the following lemma.
\begin{lemma}\label{lem:u}
	For each $k$ there is a real number $u_k$ such that
	\[
	G_k(x)=u_kJ_k(x)
	\qquad(x\in \mathbb{S}^1).
	\]
	If $\Phi$ has no zero on $\bd M$, then $(u_1,u_2)\neq(0,0)$.

\end{lemma}

\begin{proof}
	Fix $k$ and take $x,y\in \mathbb{S}^1$.  Then 
	\[
	J_k(y)\,\overline{G_k(x)}
	=F_k(y,x)
	=\overline{F_k(x,y)}
	=\overline{J_k(x)}\,G_k(y).
	\]
	So the function $r_k:=\frac{G_k}{J_k}:\mathbb{S}^1\to\C$ satisfies  
	\[
	\overline{r_k(x)}=r_k(y)
	\qquad(x,y\in \mathbb{S}^1).
	\]
	Thus $r_k$ is a real constant, say $u_k\in\R$.
	
	Write
	$g:=(G_1,G_2):\overline\disk\longrightarrow\C^2$. Now suppose $\Phi$ is nonzero on $\bd M$. On the face $\mathbb S^1\times\overline\disk$ we have
	\[
	\Phi(x,y)
	=
	\bigl(J_1(x)\overline{G_1(y)},\,
	J_2(x)\overline{G_2(y)}\bigr).
	\]
	Because $J_1$ and $J_2$ are nonzero, a zero of $g$ at some $y\in\overline\disk$ would give a zero of $\Phi$ at every $(x,y)\in \mathbb{S}^1 \times\{y\}$.  Hence $g$ is nonvanishing.  In particular,
	\[
	g|_{\mathbb S^1}=(u_1J_1,u_2J_2)
	\]
	is nonvanishing, so $(u_1,u_2)\neq(0,0)$.
\end{proof}

For $d_1,d_2\in\mathbb Z$, let $\mathcal V_{d_1d_2}$ be the space of continuous maps
\[
\Phi=(F_1,F_2):M\longrightarrow\C^2
\]
such that $F_1\in\mathcal{W}_{d_1}$, $F_2\in\mathcal W_{d_2}$, and
\begin{equation}\label{eq:admissible-pair-boundary}
	\Phi(a,b)\neq0
	\qquad\text{for every }(a,b)\in\bdr M.
\end{equation}

\begin{lemma}\label{lem:V-path-connected}
	The space $\mathcal V_{d_1d_2}$ is path connected. 
\end{lemma}

\begin{proof}
	Given any map $\Phi^k=(F_1^k,F_2^k)\in\mathcal V_{d_1d_2}$, $k=0,1$,  let $(G_j^k,J_j^k)$ be boundary data for $F_j^k$ as in \eqref{eq:W-boundary-data}, where $j=1,2$. Then, by Lemma~\ref{lem:u}, there exists $u_j^k\in\mathbb{R}$ such that $G_j^k(a)=u_j^kJ_j^k(a)$ for every $a\in\bdr\disk$. Set
	\[
	G^k:=(G_1^k,G_2^k):\overline\disk\longrightarrow\C^2,
	\qquad
	u^k:=(u_1^k,u_2^k)\in\R^2.
	\]
	For $a\in\mathbb S^1$ and $b\in\overline\disk$ we have
	\begin{equation}\label{eq:endpoint-boundary-pair}
		\Phi^k(a,b)
		=\bigl(J_1^k(a)\overline{G_1^k(b)},
		J_2^k(a)\overline{G_2^k(b)}\bigr).
	\end{equation}
	Since the factors $J_j^k$ never vanish, condition
	\eqref{eq:admissible-pair-boundary} implies that $G^k(b)\ne(0,0)$ for every $b\in\overline\disk$: otherwise \eqref{eq:endpoint-boundary-pair} would vanish on the boundary face $\mathbb S^1\times\{b\}$. In particular, $u^k\neq(0,0)$. 
	
	Choose a path
	\[
	u^s=(u_1^s,u_2^s)\in\R^2\setminus\{0\},
	\qquad s\in[0,1],
	\]
	joining $u^0$ to $u^1$.  Since $\deg J_j^0=\deg J_j^1=d_j$, 
	there exist homotopies
	\[
	J_j^s:\mathbb S^1\longrightarrow\C\setminus\{0\},
	\qquad j=1,2,
	\]
	joining $J_j^0$ to $J_j^1$ and satisfying $\deg J_j^s=d_j$ for all $s$.
	
	To construct a nowhere-zero homotopy between $G^0$ and $G^1$ with the required boundary values, define on
	\[
	\partial\bigl([0,1]\times\overline\disk\bigr)
	=\bigl(\{0,1\}\times\overline\disk\bigr)
	\cup\bigl([0,1]\times\mathbb S^1\bigr)
	\]
	the map
	\begin{equation*}
		G(s,a):=\begin{cases}
		\bigl(u_1^sJ_1^s(a),u_2^sJ_2^s(a)\bigr), & \text{ on } [0,1] \times \mathbb S^1,\\
		G^i, & \text{ on } \{i\}\times\overline{\disk}\quad(i=0,1).
		\end{cases}
			\end{equation*}
    These formulas agree on the corner circles $\{i\} \times \mathbb S^1$ and define a continuous map into
	$\C^2\setminus\{0\}$.  Since $\partial\bigl([0,1]\times\overline\disk\bigr)$ is a two-sphere and
	$\pi_2(\C^2\setminus\{0\})=\pi_2(\mathbb S^3)=0$,
	the map $G$ extends to a continuous map
	\begin{equation*}
		G=(G_1,G_2):[0,1]\times\overline\disk
		\longrightarrow\C^2\setminus\{0\}.
	\end{equation*}
	
	With $G$ and $J_j^s$ at hand, we define $\Phi: \big([0,1]\times\bdr M\big) \cup \big(\{0,1\} \times M\big) \to \mathbb C^2$ by
	\begin{equation*}
		\Phi(s,a,b)
		:=\begin{cases}
			 \bigl(J_1^s(a)\overline{G_1(s,b)},
			 J_2^s(a)\overline{G_2(s,b)}\bigr), &
			  \text{ on } [0,1] \times \mathbb S^1 \times \overline\disk,  \\
			  \bigl(\overline{J_1^s(b)}G_1(s,a),
			  \overline{J_2^s(b)}G_2(s,a)\bigr), & \text{ on } [0,1] \times \overline\disk \times \mathbb S^1, \\
			  \Phi^i(a,b), & \text{ on } \{i\} \times M\ (i=0,1).
		\end{cases},
	\end{equation*}
Then $\Phi$ is continuous, maps $[0,1] \times \bdr M$ into $\mathbb C^2 \setminus \{0\}$, and satisfies 
	\[
	\Phi(s,b,a)=\overline{\Phi(s,a,b)}
	\quad\text{on } [0,1] \times \bdr M.
	\]
By the Tietze extension theorem, there exists a continuous extension to all of
	$[0,1]\times M$. 
	For every $s$, the resulting slice $\Phi^s=\Phi(s,\cdot,\cdot)$ belongs to $\mathcal V_{d_1d_2}$.
\end{proof}

For any $\Phi\in\mathcal V_{d_1d_2}$, the Brouwer degree at the origin is given by
\[
\deg(\Phi,M,0):=\deg(\widehat\Phi)\in\mathbb Z,
\] 
where  
\[
\widehat\Phi:=\frac{\Phi}{|\Phi|}:\bdr M\longrightarrow\mathbb S^3
\]
is the normalized boundary map. For basic properties of the Brouwer degree, we refer to \cite[Chapter~4]{milnor1997topology}.

For any $t\in(0,1)$, consider the disk automorphisms
\[M_t(z):=\frac{z-t}{1-tz}.\]
For $d\in \mathbb{Z}$, define $Q_d:\overline\disk \to \mathbb{C}$ by
\[
Q_d(z)=
\begin{cases}
z^d,&d>0,\\
1,&d=0,\\
\overline z^{-d},&d<0,
\end{cases}
\]
and define the composition $M_t^d = Q_d\circ M_t$.

Taking $G=M_t^d$, $J=M_t^d|_{\mathbb S^1}$, and $u=1$, we see that the map $(a,b)\mapsto M_t^d(a)\overline{M_t^d(b)}$ belongs to $\mathcal W_d$. Consequently, the map
\begin{equation}\label{eq:model-pair}
	\widehat\Phi_t(a,b)
	:=\bigl(M_t^{d_1}(a)\overline{M_t^{d_1}(b)},
	M_{-t}^{d_2}(a)\overline{M_{-t}^{d_2}(b)}\bigr) \in \mathcal V_{d_1d_2}. 
\end{equation}
Indeed, $\widehat{\Phi}_t(a,b)\neq 0$ on $\bdr(\disk\times\disk)$, because $M_t^{d}$ is nowhere zero on $\bdr\disk$ and $M_{t}^{d_1}$ and $M_{-t}^{d_2}$ have no common zeros on $\overline\disk$.

\begin{lemma}\label{lem:degree-negative2}
	 $\deg(\widehat\Phi_t,M,0)=-2d_1d_2.$
\end{lemma}

\begin{proof}
If $d_1d_2=0$, then one component of $\widehat\Phi_t$ is identically
equal to $1$ and hence $\widehat\Phi_t$ has no zeros in $M$. Therefore
$\deg(\widehat\Phi_t,M,0)=0$, and the conclusion follows. We may thus
assume that $d_1d_2\neq0$. In this case the two components of \eqref{eq:model-pair} vanish
simultaneously precisely at
\begin{equation*}
    \widehat\Phi_t^{-1}(0)
    =
    \{(t,-t),(-t,t)\}.
\end{equation*}

Recall first that $\deg_0(Q_d)=d$.
Since $M_t$ and $M_{-t}$ are conformal diffeomorphisms near $t$ and $-t$, respectively,
it follows that
\[
\deg_t(M_t^d)=
\deg_{-t}(M_{-t}^d)=d.
\]
Moreover, complex conjugation reverses the orientation of
$\mathbb C\simeq\mathbb R^2$, and hence
\[
\deg_t\!\left(\overline{M_t^d}\right)=
\deg_{-t}\!\left(\overline{M_{-t}^d}\right)=-d.
\]

Consider first the zero $(t,-t)$. In a sufficiently small product
neighborhood of this point, the factors
$\overline{M_t^{d_1}(b)}$ and $M_{-t}^{d_2}(a)$ are nonzero. A homotopy changes these nonzero factors to nonzero constants without changing local degrees. Since multiplication by a nonzero complex
constant is orientation preserving, the local degree is the same as
that of
\[
(a,b)\longmapsto
\left(M_t^{d_1}(a),\,\overline{M_{-t}^{d_2}(b)}\right).
\]
The standard product property of the local Brouwer degree,
\[
\deg_{(x_1,x_2)}(f_1\times f_2)
=
\deg_{x_1}(f_1)\deg_{x_2}(f_2),
\]
gives
\[
\deg_{(t,-t)}(\widehat\Phi_t)
=
d_1(-d_2)
=
-d_1d_2.
\]

At the zero $(-t,t)$, the local map, up to multiplication by nonzero complex constants, is
\[
(a,b)\longmapsto
\left(\overline{M_t^{d_1}(b)},\,M_{-t}^{d_2}(a)\right).
\]
The exchanging map $(a,b)\mapsto (b,a)$ is orientation preserving, since $(-1)^{2\cdot2}=1$. Hence the product property gives
\[
\deg_{(-t,t)}(\widehat\Phi_t)
=
(-d_1)d_2
=
-d_1d_2.
\]

The Brouwer degree equals the sum of local degrees over the zero set, so
\[
\deg(\widehat\Phi_t,M,0)
=
\deg_{(t,-t)}(\widehat\Phi_t)
+
\deg_{(-t,t)}(\widehat\Phi_t)
=
-2d_1d_2.
\]
\end{proof}

We can now prove the theorem.

\begin{proof}[Proof of Theorem~\ref{thm:topo}]
Set $\Phi:=(F_1,F_2)$. If $\Phi$ vanishes at some point of $\bdr M$, then
$F_1$ and $F_2$ already have a common zero. Assume henceforth that $\Phi$ is
nonvanishing on $\bdr M$, so that $\Phi\in\mathcal V_{d_1d_2}$.

Fix $t\in(0,1)$. By Lemma~\ref{lem:V-path-connected}, $\Phi$ is homotopic,
through maps in $\mathcal V_{d_1d_2}$, to the model map $\widehat\Phi_t$
from \eqref{eq:model-pair}. Homotopy invariance of the Brouwer degree and
Lemma~\ref{lem:degree-negative2} therefore give
\[
\deg(\Phi,M,0)
=\deg(\widehat\Phi_t,M,0)
=-2d_1d_2.
\]
This proves \eqref{eq:degree-formula}.
\end{proof}

\section{Spectral Convergence and Sharpness}\label{sec:sharpness}

\subsection{Three approximation procedures}\label{sec:sharpness-approximation}

We shall use three approximation procedures. In all statements below, spectral
convergence is understood for every fixed eigenvalue index. 

\medskip
\noindent\textbf{Boundary oscillation.}
Let $D_K(m)\subset\mathbb M_K^2$ be a geodesic disk, let $c\geq1$, and write $\ell_K(m)=|\partial D_K(m)|$. There exist smooth domains $D_\varepsilon$ converging to $D_K(m)$ in the Hausdorff sense such that
\begin{equation*}
|\partial D_\varepsilon|\longrightarrow c\ell_K(m),
\end{equation*}
and, for every fixed $k$ and compact interval $I\subset \mathbb{R}$,
\begin{equation*}
\lambda_k(D_\varepsilon;\gamma)
\longrightarrow
\lambda_k(D_K(m);c\gamma),\qquad\text{uniformly for }\gamma \in I.
\end{equation*}
Indeed, parametrize $\partial D_K(m)$ by arclength $s\in\mathbb R/\ell_K(m)\mathbb Z$ and, in Fermi coordinates $(s,r)$, replace the boundary by the normal graph
\[
r_n(s)=\frac{q}{k_n}\sin(k_ns),
\qquad
k_n=\frac{2\pi n}{\ell_K(m)}.
\]
Its amplitude tends to zero, whereas its boundary Jacobian converges weakly to the constant
\[
c(q)=\frac1{2\pi}\int_0^{2\pi}\sqrt{1+q^2\cos^2\theta}\,d\theta.
\]
Since $c(q)$ ranges over $[1,\infty)$, $q$ may be chosen so that $c(q)=c$. The spectral convergence, including negative $\gamma$, follows from the signed uniform-domain lemma in Appendix~\ref{sec:Appendix_Robin_convergence}, applied in Lemma~\ref{lem:oscillating-Robin-spectrum}. This is the finite-boundary-multiplier mechanism considered in Arrieta--Bruschi~\cite{arrietabruschi2006}.

\medskip
\noindent\textbf{Short-neck connection.}
Let $\Omega_1,\Omega_2$ be two smooth domains with disjoint interiors whose
boundaries meet only at a nondegenerate tangency $p$ (their local
boundary graphs have strictly positive quadratic separation).
There exist smooth connected domains
$\widehat\Omega_\varepsilon$ such that
\begin{equation}\label{eq:short-neck-smooth-geometry}
 |\widehat\Omega_\varepsilon\triangle(\Omega_1\cup\Omega_2)|\longrightarrow0,
 \qquad
 |\partial\widehat\Omega_\varepsilon|
 \longrightarrow
 |\partial\Omega_1|+|\partial\Omega_2|,
\end{equation}
and, for every fixed $k\geq1$ and compact interval $I\subset\mathbb R$,
\begin{equation}\label{eq:sharpness-neck-spectrum}
 \lambda_k(\widehat\Omega_\varepsilon;\gamma)
 -\lambda_k(\Omega_1\sqcup\Omega_2;\gamma)
 \longrightarrow0,\qquad \text{uniformly for }\gamma \in I.
\end{equation}
If $\Omega_1$ and $\Omega_2$ are simply connected, the approximating domains may be
chosen simply connected.

To construct them, make the two domains overlap slightly near $p$, and
cut their union along the short segment joining the two intersection
points in local coordinates.  The cut pieces have uniformly Lipschitz
boundaries, whereas the segment shrinks to $p$.  Finally, round the two corners on a smaller scale;
the Robin spectra and perimeters converge under this approximation.
Details are given in Lemma~\ref{lem:short-neck-spectral-convergence}.

\medskip
\noindent\textbf{Area correction.}
Let $\Omega_\varepsilon$ be smooth domains with
$|\Omega_\varepsilon|\to M>0$ and
$\lambda_1(\Omega_\varepsilon;-b)\ge-C_b$ for every $b>0$
and all sufficiently small $\varepsilon$. Suppose also that their boundaries contain
arcs of a fixed positive length and uniformly bounded curvature,
with embedded two-sided normal collars of a fixed width.
Lemma~\ref{lem:area-correction} gives smooth domains
$\widetilde\Omega_\varepsilon$ diffeomorphic to $\Omega_\varepsilon$,
with $|\widetilde\Omega_\varepsilon|=M$ and
$|\partial\widetilde\Omega_\varepsilon|-|\partial\Omega_\varepsilon|\to0$,
such that for every fixed $k$,
\[
\lambda_k(\widetilde\Omega_\varepsilon;\gamma)
-
\lambda_k(\Omega_\varepsilon;\gamma)
\to0,\qquad\text{locally uniformly for } \gamma\in\mathbb{R}.
\]

\subsection{Approximation of the two-disks spectra}
\label{sec:sharpness-root-approximation}

For two smooth domains $D_1,D_2$ with two Robin coefficients
$\beta_1,\beta_2$, let
\[
\lambda_k(D_1,D_2;\beta_1,\beta_2)
\]
denote the $k$-th eigenvalue, counted with multiplicity, of the
disconnected problem on $D_1\sqcup D_2$, with Robin parameter 
$\beta_i$ on $\partial D_i$. Thus its spectrum is the ordered union of
the Robin spectra of $(D_1,\beta_1)$ and $(D_2,\beta_2)$.

Let $m_1,m_2>0$, with $K(m_1+m_2)<4\pi$, and set
\[
D_i:=D_K(m_i),
\qquad
\ell_i:=\ell_K(m_i),
\qquad i=1,2.
\]

\begin{proposition}\label{prop:sharpness-coefficient-approximation}
For any $r_1,r_2>0$, there exist smooth simply connected domains
$\Omega_\varepsilon\subset\mathbb M_K^2$ with
\(
|\Omega_\varepsilon|=m_1+m_2
\)
such that, for every fixed $k$ and every compact interval
$I\subset\mathbb R$,
\begin{equation}\label{eq:sharpness-spectrum-approximation}
\lambda_k\left(
\Omega_\varepsilon;
\frac{a}{|\partial\Omega_\varepsilon|}
\right)
\longrightarrow
\lambda_k\left(
D_1,D_2;
\frac{a}{r_1+r_2}\frac{r_1}{\ell_1},
\frac{a}{r_1+r_2}\frac{r_2}{\ell_2}
\right)
\qquad\text{uniformly for }a\in I.
\end{equation}
\end{proposition}

\begin{proof}
Interchanging the indices if necessary, we may assume that
\[
\frac{r_1}{\ell_1}\geq\frac{r_2}{\ell_2}.
\]
Set
\[
c:=\frac{r_1\ell_2}{r_2\ell_1}\geq1.
\]
By Lemma \ref{lem:oscillating-Robin-spectrum}, through oscillating the boundary $\partial D_1$, we obtain simply
connected domains $D_{1,\varepsilon}$ converging to $D_1$
in the Hausdorff sense, with
$|\partial D_{1,\varepsilon}|\to c|\partial D_1|$,
such that, for every fixed $k$,
\[
\lambda_k(D_{1,\varepsilon};\gamma)
\longrightarrow\lambda_k(D_1;c\gamma)
\]
locally uniformly for $\gamma\in\mathbb R$.

For each fixed $\varepsilon$, we may position
$D_{1,\varepsilon}$ and $D_2$ so that their interiors
are disjoint and their boundaries are nondegenerately tangent at a point $p$.
Connecting them by a short neck, Lemma \ref{lem:short-neck-spectral-convergence} gives smooth simply
connected domains $\Omega_{\varepsilon,\delta}$ such that 
\[
\lim_{\delta\downarrow 0}|\Omega_{\varepsilon,\delta}\triangle (D_{1,\varepsilon}\cup D_2)|= 0,
\qquad
\lim_{\delta\downarrow0}|\partial\Omega_{\varepsilon,\delta}|
=|\partial D_{1,\varepsilon}|+|\partial D_2|,
\]
and, for every fixed $k$,
\[
\lambda_k(\Omega_{\varepsilon,\delta};\gamma)
\longrightarrow
\lambda_k(D_{1,\varepsilon}\sqcup D_2;\gamma)
\]
locally uniformly for $\gamma\in\mathbb R$.

Choose $\delta=\delta(\varepsilon)$ sufficiently small that
\[
\sup_{\substack{k\leq\varepsilon^{-1}\\
|\gamma|\leq\varepsilon^{-1}}}
\left|
\lambda_k(\Omega_{\varepsilon,\delta(\varepsilon)};\gamma)
-\lambda_k(D_{1,\varepsilon}\sqcup D_2;\gamma)
\right|<\varepsilon.
\]
Moreover, \(
|\Omega_{\varepsilon,\delta(\varepsilon)}\triangle (D_{1,\varepsilon}\cup D_2)|\to0\) and \(
|\partial\Omega_{\varepsilon,\delta(\varepsilon)}|\to c\peri{D_1} + \peri{D_2}
\) as $\varepsilon\downarrow 0$.

Set $\widehat\Omega_\varepsilon
:=\Omega_{\varepsilon,\delta(\varepsilon)}$.
Then, for every fixed $k$,
\[
\lambda_k(\widehat\Omega_\varepsilon;\gamma)
\longrightarrow
\lambda_k(D_1,D_2;c\gamma,\gamma)
\]
locally uniformly for $\gamma\in\mathbb R$.

For $a\in\mathbb R$, set
\(
\gamma_\varepsilon
:={a}/{|\partial\widehat\Omega_\varepsilon|}.
\)
Since
\[
\gamma_\varepsilon
\longrightarrow
\frac{a}{r_1+r_2}\frac{r_2}{\ell_2},
\qquad
c\gamma_\varepsilon
\longrightarrow
\frac{a}{r_1+r_2}\frac{r_1}{\ell_1},
\]
we obtain
\[
\lambda_k\left(
\widehat\Omega_\varepsilon;
\frac{a}{|\partial\widehat\Omega_\varepsilon|}
\right)
\longrightarrow
\lambda_k\left(
D_1,D_2;
\frac{a}{r_1+r_2}\frac{r_1}{\ell_1},
\frac{a}{r_1+r_2}\frac{r_2}{\ell_2}
\right).
\]
This convergence is locally uniform for $a\in\mathbb R$.
For every $b>0$, we also have
\[
\lambda_1(\widehat\Omega_\varepsilon;-b)
\longrightarrow\lambda_1(D_1,D_2;-cb,-b)>-\infty.
\]
Away from the neck, $\partial\widehat\Omega_\varepsilon$ contains arcs satisfying the uniform collar assumptions of
Lemma~\ref{lem:area-correction}. So there exist smooth simply connected domains $\Omega_\varepsilon$ with
$|\Omega_\varepsilon|=m_1+m_2$ and
$|\partial\Omega_\varepsilon|\to c\ell_1+\ell_2$.
Since  $a/|\partial\widehat\Omega_\varepsilon|$ and $a/|\partial\Omega_\varepsilon|$ have the same limit,  the difference of eigenvalues tends to zero locally uniformly for $a\in\mathbb{R}$, proving \eqref{eq:sharpness-spectrum-approximation}.
\end{proof}

\subsection{Sharpness and counterexamples.}

The approximation geometries used below are shown in
Figure~\ref{fig:sharpness-geometries}. First we show that the bound \eqref{eq:intro-Robin-main} is sharp among bounded simply connected Lipschitz domains of area $A$ in $\mathbb M_K^2$.

\begin{proposition}\label{prop:two-disk-sharpness}
If $(K,A,\alpha)\in\mathcal P_1$, then there exist  bounded simply connected smooth domains  $\region_\varepsilon \subset \mathbb M_K^2$ with $A=|\region_\varepsilon|$,  such that 
\[
	\lambda_3\left(\region_\varepsilon;\frac{\alpha}{|\partial\region_\varepsilon|}\right)
	\to
	\lambda_2\left(D_K(A/2);\frac{\alpha}{2\ell_K(A/2)}\right)
    \]
as $\varepsilon \to 0+$.
\end{proposition}

\begin{proof}
Take two isometric copies of $D_K(A/2)$ with tangent boundaries. By the short-neck construction of Subsection~\ref{sec:sharpness-approximation}, followed by the area correction, there are smooth simply connected domains $\Omega_\varepsilon$ of area $A$ such that
\[
|\partial\Omega_\varepsilon|
\longrightarrow
2\ell_K(A/2)
\]
and
\[
\lambda_3\left(
\Omega_\varepsilon;\frac{\alpha}{|\partial\Omega_\varepsilon|}
\right)
\longrightarrow
\lambda_3\left(
D_K(A/2)\sqcup D_K(A/2);\frac{\alpha}{2\ell_K(A/2)}
\right).
\]
The two components have identical spectra, so the last quantity equals
\[
\lambda_2\left(
D_K(A/2);\frac{\alpha}{2\ell_K(A/2)}
\right).
\]
Together with \eqref{eq:intro-Robin-main}, this proves sharpness.
\end{proof}
In the case when $(K,A,\alpha)\in \mathcal P_1$ and $\alpha=-4\pi+KA$, the sharp upper bound \eqref{eq:intro-Robin-main} admits extremizing sequences with unequal limiting disk areas, and thus is not confined to the equal-area split.   

\begin{proposition}\label{prop:critical-sharpness-instability}
Let $K>0$, $0 <KA < 4\pi$ and $\alpha=-4\pi+KA$. For $\delta > 0$ such that $m_-:=A/2-\delta > 0$ and $m_+:= A/2 + \delta < 2\pi/K$, there exist smooth simply connected domains
$\Omega_\varepsilon\subset\mathbb M_K^2$ of area $A$, obtained by
joining $D_K(m_+)$ to an oscillated disk converging to $D_K(m_-)$,
such that
\[
\lambda_3\left(\Omega_\varepsilon;
\frac{\alpha}{|\partial\Omega_\varepsilon|}\right)
\longrightarrow 2K=\lambda_2\left(
D_K(A/2);\frac{\alpha}{2\ell_K(A/2)}
\right).
\]
Moreover, 
\[
|\partial \Omega_{\varepsilon}|\longrightarrow\frac{-\alpha \ell_K(m_+)}{2\pi-Km_+}=\frac{-\alpha\ell_K(m_+)}{R_{K,2K}(m_+)}.
\]
\end{proposition}

\begin{proof}
Set
\[
\begin{gathered}
r_1:=2\pi-Km_-,\qquad r_2:=2\pi-Km_+,\\
\ell_1:=\ell_K(m_-),\qquad \ell_2:=\ell_K(m_+).
\end{gathered}
\]
Then $r_1\geq r_2>0$, $\ell_1\leq\ell_2$, and hence
\[
\frac{r_1}{\ell_1}\geq\frac{r_2}{\ell_2},
\qquad r_1+r_2=-\alpha.
\]
Let $c:=r_1\ell_2/(r_2\ell_1)\geq1$. Proposition~\ref{prop:sharpness-coefficient-approximation} gives smooth simply connected domains $\Omega_\varepsilon$ obtained by joining $D_K(m_+)$ to an oscillated disk converging to $D_K(m_-)$ with
\[
|\Omega_{\varepsilon}| = A,\qquad
|\partial\Omega_\varepsilon|
\longrightarrow c\ell_1+\ell_2
=\frac{r_1+r_2}{r_2}\ell_2
\]
and
\begin{equation}\label{eq:critical-sharpness-spectrum}
\lambda_3\left(\Omega_\varepsilon;
\frac{\alpha}{|\partial\Omega_\varepsilon|}\right)
\longrightarrow
\lambda_3\left(
D_K(m_-),D_K(m_+);
-\frac{r_1}{\ell_1},-\frac{r_2}{\ell_2}
\right).
\end{equation}
By $R_{K,2K}(x)=2\pi-Kx$ and
\eqref{eq:model-angular-boundary-early},
\[
\lambda_{1,\mathrm{ang}}\left(D_K(m_-);-\frac{r_1}{\ell_1}\right)
=
\lambda_{1,\mathrm{ang}}\left(D_K(m_+);-\frac{r_2}{\ell_2}\right)
=2K.
\]
Since both disk areas are smaller than $2\pi/K$, these are
the second Robin eigenvalues of the respective disks
by~\cite[Proposition~6.1]{langfordlaugesen2025robin}.
Each disk has precisely one eigenvalue below $2K$, so the
right-hand side of \eqref{eq:critical-sharpness-spectrum} equals
$2K$.
\end{proof}

As $\alpha$ drops below $-4\pi + KA$, a similar construction gives counterexamples to the two-disk inequality \eqref{eq:intro-Robin-main}.

\begin{proposition}\label{prop:equal-counterexample}
Let $K>0$, $0<KA<4\pi$, and $-4\pi<\alpha<-4\pi+KA$. For $\delta > 0$ such that $m_-:=A/2-\delta > 0$ and $m_+:= A/2 + \delta < 2\pi/K$, there exists a smooth simply connected domain $\region\subset\mathbb M_K^2$ of area $A$, obtained by joining the standard disk $D_K(m_+)$ to an oscillated disk approximating $D_K(m_-)$, such that
\begin{equation}\label{eq:equal-counterexample}
\lambda_3\left(\region;\frac{\alpha}{|\partial\region|}\right)
> \mu_{\rm eq}:=
\lambda_2\left(
D_K(A/2);\frac{\alpha}{2\ell_K(A/2)}
\right).
\end{equation}
\end{proposition}

\begin{proof}
Since $A/2<2\pi/K$, $\mu_{\rm eq}$ is the first angular eigenvalue of $D_K(A/2)$ by~\cite[Proposition~6.1]{langfordlaugesen2025robin}; hence
\begin{equation*}
2R_{K,\mu_{\rm eq}}(A/2)+\alpha=0.
\end{equation*}
At $\mu=0$ the left-hand side is $4\pi+\alpha>0$, whereas at $\mu=2K$ it is $4\pi-KA+\alpha<0$. Thus
\[
0<\mu_{\rm eq}<2K.
\]
Strict convexity in Proposition~\ref{prop:model-area-concavity} gives
\[
R_{K,\mu_{\rm eq}}(m_-)+R_{K,\mu_{\rm eq}}(m_+)
>
2R_{K,\mu_{\rm eq}}(A/2)
=-\alpha.
\]
Since the left-hand side is strictly decreasing in $\mu$, there is a unique $
\mu_\delta\in(\mu_{\rm eq},2K)$ 
such that
\begin{equation*}
R_{K,\mu_\delta}(m_-)
+R_{K,\mu_\delta}(m_+)
+\alpha=0.
\end{equation*}
Moreover, $R_{K,\mu_\delta}'(x)<0$ on $(0,2\pi/K)$.
Indeed, $R_{K,\mu_\delta}'(0+)=-\mu_\delta/2<0$;
at a first zero of $R'$ in this interval, differentiating
\eqref{eq:model-riccati} would give
$R''=-\mu_\delta G_K'/G_K<0$, a contradiction.
Since $\mu_\delta<2K$, Lemma~\ref{lem:model-R-critical-identity}
gives $R_{K,\mu_\delta}(x)>2\pi-Kx$ on $I_K$. Hence
\[
R_{K,\mu_\delta}(m_-)
>
R_{K,\mu_\delta}(m_+)
>0.
\]
Since also $\ell_K(m_-)<\ell_K(m_+)$, we have
\[
\frac{R_{K,\mu_\delta}(m_-)}{\ell_K(m_-)}
>
\frac{R_{K,\mu_\delta}(m_+)}{\ell_K(m_+)}.
\]
Since both $D_{K}(m_-)$ and $D_K(m_+)$ are smaller than a hemisphere, their second Robin eigenvalues are their first angular eigenvalues and are equal to $\mu_\delta$. As in Proposition~\ref{prop:critical-sharpness-instability}, Proposition~\ref{prop:sharpness-coefficient-approximation} gives $\Omega_{\varepsilon}$ by joining $D_K(m_+)$ to an oscillated $D_K(m_-)$, such that $|\Omega_{\varepsilon}| = A$ and $\lambda_3(\Omega_\varepsilon;\alpha/|\partial \Omega_{\varepsilon}|)\to \mu_{\delta} > \mu_{\rm eq}$. Hence one member of the approximating sequence satisfies \eqref{eq:equal-counterexample}.
\end{proof}

Assume now that $(K,A,\alpha)\in\mathcal P_2$, $\alpha\leq-2\pi$, and
\begin{equation}\label{eq:one-disk-angular-radial-condition}
\lambda_{1,\mathrm{ang}}\left(
D_K(A);\frac{\alpha+2\pi}{\ell_K(A)}
\right)
\leq
\lambda_{2,\mathrm{rad}}\left(
D_K(A);\frac{\alpha+2\pi}{\ell_K(A)}
\right).
\end{equation}
We show that the bound \eqref{eq:intro-Robin-one-disk} is sharp.
\begin{proposition}\label{prop:one-disk-sharpness}
Under these assumptions, there exist smooth simply connected domains
$\Omega_j\subset\mathbb M_K^2$ of area $A$ such that
\[
\lambda_3\left(\Omega_j;\frac{\alpha}{|\partial\Omega_j|}\right)
\longrightarrow \nu:=
\lambda_{1,\mathrm{ang}}\left(
D_K(A);\frac{\alpha+2\pi}{\ell_K(A)}
\right)
=
\lambda_2\left(D_{K}(A);\frac{\alpha + 2\pi}{\ell_K(A)}\right).
\]
If $-4\pi<\alpha<-4\pi+KA$, the domains may be chosen so that
\[
\frac1{|\partial\Omega_j|}
\longrightarrow
\frac{\alpha+2\pi}{\alpha\ell_K(A)}.
\]
\end{proposition}

\begin{proof}
If $\alpha=-4\pi$, both comparison values equal zero. If
$\alpha=-4\pi+KA$, the identity $R_{K,2K}(x)=2\pi-Kx$ shows that
both values equal $2K$. In these cases the claim follows from
Proposition~\ref{prop:two-disk-sharpness}. Thus assume
\[
-4\pi<\alpha<-4\pi+KA.
\]
By definition, $\nu$ is the unique solution of
\[
F_0(\mu):=2\pi+R_{K,\mu}(A)+\alpha=0.
\]
For $m\in(0,A)$, let $\mu(m)$ be the unique solution of
\[
F_m(\mu):=R_{K,\mu}(m)+R_{K,\mu}(A-m)+\alpha=0.
\]
Both $F_0$ and $F_m$ are strictly decreasing, positive at $0$, and negative at $2K$. Hence $0<\nu,\mu(m)<2K$. Since $F_m$ converges pointwise to $F_0$ as $m\downarrow0$, their zeros satisfy $\mu(m)\to\nu$.

Set
\[
r_1(m):=R_{K,\mu(m)}(m),\qquad
r_2(m):=R_{K,\mu(m)}(A-m).
\]
Then $r_1(m)+r_2(m)=-\alpha$. By the singular-end expansion
\eqref{eq:model-R-expansion} and $\mu(m)\to\nu$,
\[
\begin{aligned}
r_1(m)&=2\pi-\frac{\nu}{2}m+o(m),\\
r_2(m)&=-\alpha-2\pi+\frac{\nu}{2}m+o(m).
\end{aligned}
\]
Since $\alpha\leq-2\pi$ and $\nu>0$, both $r_1(m)$ and $r_2(m)$
are positive for all sufficiently small $m$. Moreover,
\[
\frac{r_1(m)}{\ell_K(m)}
>
\frac{r_2(m)}{\ell_K(A-m)},
\]
since the first quotient diverges and the second is bounded.

For each such fixed $m$, Proposition~\ref{prop:sharpness-coefficient-approximation}
provides smooth simply connected domains $\Omega_{m,n}$ of area $A$
such that, as $n\to\infty$,
\begin{equation}\label{eq:lambda_3_limit}
\lambda_3\left(\Omega_{m,n};\frac{\alpha}{|\partial\Omega_{m,n}|}\right)
\longrightarrow
\lambda_3\left(D_K(m),D_K(A-m);
-\frac{r_1(m)}{\ell_K(m)},-\frac{r_2(m)}{\ell_K(A-m)}\right).
\end{equation}
The parameter $c$ in Proposition \ref{prop:sharpness-coefficient-approximation} is
$c_m=r_1(m)\ell_K(A-m)/(r_2(m)\ell_K(m))$, so 
\[
|\partial\Omega_{m,n}|\longrightarrow
L_m:=c_m\ell_K(m)+\ell_K(A-m)
=\frac{-\alpha\,\ell_K(A-m)}{r_2(m)}.
\]

Both disks have first angular eigenvalue $\mu(m)$. On the small disk
this is its second eigenvalue by
\cite[Proposition~6.1]{langfordlaugesen2025robin}. On the large disk,
\eqref{eq:one-disk-angular-radial-condition} and continuity of the
angular and radial eigenvalues yield
\[
\lambda_2\left(D_K(A-m);-\frac{r_2(m)}{\ell_K(A-m)}\right)
\longrightarrow\nu.
\]
Each disk has its first eigenvalue strictly below $\mu(m)$.
Thus the third eigenvalue of the disconnected problem in
\eqref{eq:lambda_3_limit} tends to $\nu$ as $m\downarrow0$.
Choosing $m\downarrow0$ and then $n=n(m)$ sufficiently large gives the
asserted approximating sequence. Moreover, 
\[
\lim_{m\downarrow 0} \frac1{|\partial\Omega_{m,n(m)}|}
=
\lim_{m\downarrow0}\frac{r_2(m)}{-\alpha\,\ell_K(A-m)}
=\frac{\alpha+2\pi}{\alpha\ell_K(A)}.
\]
This proves the proposition.
\end{proof}
\begin{remark} 
{Moreover, for $-4\pi<\alpha<-4\pi+KA$,} the large and small disks in Proposition~\ref{prop:one-disk-sharpness} carry limiting Robin boundary masses $\alpha + 2\pi$ and $-2\pi$, respectively.  The sequence may be chosen to collapse at any
$p\in\partial D_K(A)$. For each fixed $m$, the oscillation and short-neck
constructions give
\[
ds_{\partial\Omega_{m,n}}\rightharpoonup
ds_{\partial D_K(A-m)}+c_m\,ds_{\partial D_K(m)}.
\]
By the definitions of $c_m$ and $L_m$,
\[
\frac{\alpha\ell_K(A-m)}{L_m}
=-r_2(m)\longrightarrow\alpha+2\pi,
\qquad
\frac{\alpha c_m\ell_K(m)}{L_m}
=-r_1(m)\longrightarrow-2\pi.
\]
Thus a diagonal choice yields
\[
\frac{\alpha}{|\partial\Omega_j|}
\,ds_{\partial\Omega_j}
\rightharpoonup
\frac{\alpha+2\pi}{\ell_K(A)}
\,ds_{\partial D_K(A)}
-2\pi\delta_p.
\]
This remains valid for $\alpha=-2\pi$, even though
$L_m\to\infty$.
\end{remark}

The restriction $\alpha\leq-2\pi$ ensures that both effective Robin
coefficients in the construction are negative. Condition
\eqref{eq:one-disk-angular-radial-condition} prevents the second radial
eigenvalue of the large disk from falling below the angular limiting
value. Outside these conditions the upper bound
\eqref{eq:intro-Robin-one-disk} remains valid, but the present
construction does not prove sharpness.

\section{The Robin--Steklov Duality}\label{sec:Robin-Steklov-duality}

For $k\geq1$ and $\lambda\in\mathbb R$, set
\begin{equation}\label{eq:shifted-SH-definition-Robin}
\Sigma_k(\Omega,\lambda)
:=-\sup\left\{\beta\in\mathbb R:\lambda_k(\Omega;\beta)\leq\lambda\right\},
\end{equation}
with the value $-\infty$ if the set is unbounded above. 
Consequently, for every $\lambda\in\mathbb R$ and every $\beta\in\mathbb R$,
\begin{equation}\label{eq:global-Robin-SH-duality}
	\lambda_k(\Omega;\beta)\leq\lambda
	\quad\Longleftrightarrow\quad
	\Sigma_k(\Omega,\lambda)\leq-\beta.
\end{equation}
Moreover, $\lambda_k(\Omega;\beta)<\lambda$ implies
$\Sigma_k(\Omega,\lambda)<-\beta$: by continuity in $\beta$,
$\lambda_k(\Omega;\beta+\varepsilon)<\lambda$ for some
$\varepsilon>0$, and \eqref{eq:global-Robin-SH-duality} gives
$\Sigma_k(\Omega,\lambda)\leq-\beta-\varepsilon$.

Since $\beta\mapsto\lambda_k(\Omega;\beta)$ is continuous and strictly increasing and tends to $\lambda_k^D(\Omega)$ as $\beta\to+\infty$, $\Sigma_k(\Omega,\lambda)$ is finite for $\lambda<\lambda_k^D(\Omega)$ and equals $-\infty$ for $\lambda\geq\lambda_k^D(\Omega)$.

The following proposition identifies $\Sigma_{k}(\region,\lambda)$ as the Helmholtz--Steklov eigenvalue with an index shifted according to the number of Dirichlet eigenvalues crossed.

\begin{proposition}[{\cite[Proposition~5]{arendtmazzeo2012}}]\label{prop:shifted-Robin-SH-duality}
If $\lambda\notin\operatorname{Spec}(-\Delta_\Omega^D)$ and
\[
N^D_\Omega(\lambda):=\#\{j:\lambda_j^D(\Omega)<\lambda\},
\]
then
\begin{equation*}
\Sigma_k(\Omega,\lambda)
=
\begin{cases}
\sigma_{k-N^D_\Omega(\lambda)}(\Omega,\lambda),&N^D_\Omega(\lambda)<k,\\[1mm]
-\infty,&N^D_\Omega(\lambda)\geq k.
\end{cases}
\end{equation*}

\end{proposition}

In particular, for $\Sigma_3$, which is our main concern, we have
\begin{equation*}
\Sigma_3(\Omega,\lambda)
=
\begin{cases}
\sigma_3(\Omega,\lambda), &\lambda < \lambda_1^D(\Omega),
\\[1mm]
\sigma_2(\Omega, \lambda), &\lambda_1^D(\Omega) < \lambda < \lambda_2^D(\Omega),
\\[1mm]
\sigma_1(\Omega, \lambda), &\lambda_2^D(\Omega) < \lambda < \lambda_3^D(\Omega),
\\[1mm]
-\infty, &\lambda > \lambda_3^D(\Omega).
\end{cases}
\end{equation*}
Below the first Dirichlet eigenvalue, we have $\Sigma_k=\sigma_k$. Between the first and second Dirichlet eigenvalues, $\Sigma_3=\sigma_2$; between the second and third, $\Sigma_3=\sigma_1$; and $\Sigma_3=-\infty$ once $\lambda\geq\lambda_3^D$.

We now relate the parameter ranges $\mathcal P_i$ and $\mathcal Q_i$ defined in Subsection~\ref{sec:intro-main-results} through \eqref{eq:global-Robin-SH-duality}.

\begin{lemma}\label{lem:dual}
There is a one-to-one correspondence between $(K,A,\alpha)\in\mathcal P_1$ and $(K,A,\Lambda)\in\mathcal Q_1$ given by
    \begin{equation}\label{eq:two-disk-alpha-lam}
\frac{\Lambda}{A}=\lambda_2\left(D_{K}(A/2);\frac{\alpha}{2\ell_K(A/2)}\right),
\end{equation}
and a one-to-one correspondence between $(K,A,\alpha)\in\mathcal P_2$ and $(K,A,\Lambda)\in\mathcal Q_2$ given by
\begin{equation}\label{eq:one-disk-alpha-lam}
\frac{\Lambda}{A}=\lambda_{1,\mathrm{ang}}\left(D_K(A);\frac{\alpha + 2\pi}{\ell_K(A)}\right).
\end{equation}
\end{lemma}
\begin{proof}
The right-hand side of \eqref{eq:two-disk-alpha-lam} is continuous and strictly increasing in $\alpha$. The value $\Lambda$ is the unique solution of
\[
R_{K,\Lambda/A}(A/2)+\frac{\alpha}{2}=0.
\]
Since $R_{K,0}(A/2) = 2\pi$ and $R_{K,2K}(A/2) = 2\pi - KA/2$, the relation \eqref{eq:two-disk-alpha-lam} gives
\[
\Lambda = \begin{cases}
    0,&\alpha = -4\pi,\\[1mm]
    2KA,&\alpha=-4\pi + KA,\\[1mm]
    A\lambda_+(K,A),&\alpha = 4\pi.
\end{cases}
\]
The last case follows from the definition. Thus, for fixed $K,A$, as $\alpha$ ranges over $\{\alpha:(K,A,\alpha)\in\mathcal P_1\}$, the range of $\Lambda$ is
\[
\begin{cases}
    \{0\}\cup [2KA, A\lambda_+(K,A)], & K> 0,\\
    [0, A\lambda_+(K,A)],& K\leq 0. \\
\end{cases}
\]
Hence $(K,A,\alpha)\mapsto(K,A,\Lambda)$ defined by \eqref{eq:two-disk-alpha-lam} gives a one-to-one correspondence between $\mathcal P_1$ and $\mathcal Q_1$.

Next, the right-hand side of \eqref{eq:one-disk-alpha-lam} is continuous and strictly increasing in $\alpha$. The value $\Lambda$ is the unique solution of 
\[
R_{K,\Lambda/A}(A) + \alpha + 2\pi = 0.
\]
Since $R_{K,0}(A) = 2\pi$ and $R_{K,2K}(A) = 2\pi - KA$, the relation \eqref{eq:one-disk-alpha-lam} gives
\[
\Lambda = \begin{cases}
    0, &\alpha = -4\pi,\\[1mm]
    2KA,&\alpha = -4\pi + KA.
\end{cases}
\]
Thus, as $\alpha$ ranges over $\{\alpha:(K,A,\alpha)\in\mathcal P_2\}$, the range of $\Lambda$ is $[0,2KA]$. Hence $(K,A,\alpha)\mapsto(K,A,\Lambda)$ defined by \eqref{eq:one-disk-alpha-lam} gives a one-to-one correspondence between $\mathcal P_2$ and $\mathcal Q_2$.
\end{proof}

To establish sharpness in \eqref{eq:intro-Steklov-Helmholtz-two} and \eqref{eq:intro-Steklov-Helmholtz-one} from Robin spectral approximations, we use the following monotonicity principle.

For a sequence of smooth domains $\Omega_n$ with $L_n = |\partial \Omega_n|$, consider the functions $F_n(a):=\lambda_k(\Omega_n;a/L_n)$.
\begin{lemma}\label{lem:sharpness-transfer}
Suppose there exists a strictly increasing function $F_*:\mathbb{R}\to\mathbb{R}$ such that $F_n$ converges pointwise to $F_*$. If $\lambda\in\mathbb R$ is such that there exists $a_*\in \mathbb R$ satisfying $\lambda = F_*(a_*)$, then, for all sufficiently large $n$, there exists $a_n$ such that $F_n(a_n) = \lambda$ and
    \begin{equation}\label{eq:shifted-SH-sharpness-transfer}
    L_n\Sigma_k(\Omega_n;\lambda)=-a_n\longrightarrow -a_*.
    \end{equation}
\end{lemma}
\begin{proof}
Fix $\varepsilon>0$. Since $F_*$ is strictly increasing and
$F_*(a_*)=\lambda$, we have
\[
F_*(a_*-\varepsilon)<\lambda<F_*(a_*+\varepsilon).
\]
By pointwise convergence, for all sufficiently large $n$,
\[
F_n(a_*-\varepsilon)<\lambda<F_n(a_*+\varepsilon).
\]
Since $F_n$ is continuous and strictly increasing, there is a unique
$a_n\in(a_*-\varepsilon,a_*+\varepsilon)$ such that $F_n(a_n)=\lambda$. Taking $\varepsilon\downarrow 0$ shows that $a_n\to a_*$. The Robin--Steklov duality then yields
\[
L_n\Sigma_k(\Omega_n,\lambda)=-a_n\longrightarrow-a_*.
\]
\end{proof}

\begin{proof}[Proof of Theorem~\ref{thm:intro-Steklov-Helmholtz} in the $\mathcal Q_1$ range]
Let $\alpha \in \mathbb R$ be given by \eqref{eq:two-disk-alpha-lam}. Lemma~\ref{lem:dual} shows that $(K,A,\alpha)\in \mathcal P_1$. Theorem~\ref{thm:intro-Robin-main} then gives
\[
\lambda_3\left(\Omega;\frac{\alpha}{L}\right)
<\lambda_2\left(D_{K}(A/2);\frac{\alpha}{2\ell_{K}(A/2)}\right)
=\frac{\Lambda}{A}.
\]
By the strict Robin--Steklov implication following \eqref{eq:global-Robin-SH-duality}, this implies
\[L\Sigma_3\left(\Omega,\frac{\Lambda}{A}\right)<-\alpha=2\ell_K(A/2)\Sigma_2\left(D_K(A/2);\frac{\Lambda}{A}\right).\]
Hence the required inequality follows.

For sharpness, Proposition~\ref{prop:two-disk-sharpness} gives
smooth domains $\Omega_n$ such that, for every $a\in\mathbb R$,
\[
\lambda_3\left(\Omega_n;\frac{a}{|\partial\Omega_n|}\right)
\longrightarrow\lambda_2\left(D_{K}(A/2);\frac{a}{2\ell_{K}(A/2)}\right).
\]
The limit is strictly increasing and equals $\lambda$ at $a=\alpha$.
Lemma~\ref{lem:sharpness-transfer} therefore proves sharpness.
\end{proof}

\begin{proof}[Proof of Theorem~\ref{thm:intro-Steklov-Helmholtz} in the $\mathcal Q_2$ range]
Let $\alpha \in \mathbb R$ be given by \eqref{eq:one-disk-alpha-lam}. Lemma~\ref{lem:dual} shows that $(K,A,\alpha)\in \mathcal P_2$. Theorem~\ref{thm:intro-Robin-main} then gives
\[
\lambda_3\left(\Omega;\frac{\alpha}{L}\right)
<\lambda_{1,\mathrm{ang}}\left(D_{K}(A);\frac{\alpha+2\pi}{\ell_{K}(A)}\right)
=\frac{\Lambda}{A}.
\]
The same strict Robin--Steklov implication gives
\[
L\Sigma_3\left(\Omega,\frac{\Lambda}{A}\right)<-\alpha.
\]
Recalling the definitions of $\lambda_{1,\mathrm{ang}},\sigma_{1,\mathrm{ang}}$ and $R_{K,\mu}$, we have $R_{K,\Lambda/A}(A) + \alpha + 2\pi = 0$ and hence
\begin{equation}\label{eq:alpha_lambda_ang}
-\alpha = 2\pi +\ell_{K}(A)\sigma_{1,\mathrm{ang}}\left(D_K(A),\frac{\Lambda}{A}\right).
\end{equation}
The inequality \eqref{eq:intro-Steklov-Helmholtz-one} then follows.

We next prove sharpness. Equation~\eqref{eq:alpha_lambda_ang} gives $\beta = ({\alpha + 2\pi})/{\ell_K(A)}$, and hence $\lambda_{1,\rm ang}(D_K(A);\beta) = \Lambda / A$. The second condition in \eqref{eq:shifted-SH-one-disk-sharp}
ensures that $$\lambda_2(D_K(A);\beta)=\lambda_{1,\mathrm{ang}}(D_K(A);\beta)=\Lambda/ A.$$

For small $m>0$, let $\beta_m$ be the unique Robin coefficient satisfying
\[
\lambda_2(D_K(A-m);\beta_m)=\Lambda/A.
\]
For every fixed $b\in\mathbb R$,
\[
\lim_{m\to0}\lambda_2(D_K(A-m);b)=\lambda_2(D_K(A);b).
\]
Lemma~\ref{lem:sharpness-transfer} therefore gives the existence of
$\beta_m$ for small $m$ and shows that $\beta_m\to\beta<0$.

Set
\[r_1(m):=R_{K,\Lambda/A}(m),\quad r_2(m):=-\beta_m\ell_K(A-m),\quad \alpha_m = -(r_1(m)+r_2(m)).\]
Since $D_K(m)$ lies inside a hemisphere, \cite[Proposition~6.1]{langfordlaugesen2025robin} gives that
\[
\lambda_2\left(D_K(m);\frac{-r_1(m)}{\ell_K(m)}\right)=\lambda_{1,\rm ang}\left(D_K(m);\frac{-r_1(m)}{\ell_K(m)}\right)=\frac{\Lambda}{A}.
\]
By the definition of $\beta_m$, the second eigenvalue of
$D_K(A-m)$ with Robin coefficient $-r_2(m)/\ell_K(A-m)=\beta_m$
also equals $\Lambda/A$. The first eigenvalue of each disk is strictly smaller.

Consequently,
\[
\lambda_3\left(
D_K(m),D_K(A-m);
-\frac{r_1(m)}{\ell_K(m)},
-\frac{r_2(m)}{\ell_K(A-m)}
\right)=\frac{\Lambda}{A}.
\]
Now define
\[
\begin{aligned}
F_m(a):=\lambda_3\bigl(D_K(m),D_K(A-m); -\frac{a\,r_1(m)}{\alpha_m\ell_K(m)},
 -\frac{a\,r_2(m)}{\alpha_m\ell_K(A-m)}\bigr).
\end{aligned}
\]

For each sufficiently small fixed $m > 0$, we have $r_1(m), r_2(m) > 0$. Proposition~\ref{prop:sharpness-coefficient-approximation} then provides smooth simply connected domains $\Omega_{m,n}$ of area $A$.
Writing $L_{m,n}=|\partial\Omega_{m,n}|$, we have, for each $a$,
\[
\lambda_3\left(\Omega_{m,n};\frac{a}{L_{m,n}}\right)
\longrightarrow F_m(a)\qquad(n\to\infty).
\]
The function $F_m$ is strictly increasing and satisfies $F_m(\alpha_m)=\lambda$. Lemma~\ref{lem:sharpness-transfer} therefore gives
\[
L_{m,n}\Sigma_3(\Omega_{m,n},\lambda)
\longrightarrow-\alpha_m\qquad(n\to\infty).
\]
Since $\alpha_m \to -2\pi + \beta\ell_K(A)=\alpha$, a diagonal choice as $m\downarrow0$ then gives
\[
L_{m,n(m)}\Sigma_3(\Omega_{m,n(m)},\lambda)
\longrightarrow-\alpha
=2\pi+\ell_K(A)\sigma_{1,\mathrm{ang}}(D_K(A),\lambda),
\]
proving sharpness.
\end{proof}
Moreover, at the endpoint values $\Lambda=0$ and $\Lambda=2KA$, the right-hand side of \eqref{eq:intro-Steklov-Helmholtz-one} agrees with that of \eqref{eq:intro-Steklov-Helmholtz-two}.  The sharpness condition $\beta < 0$ corresponds to $\alpha < -2\pi$ in the Robin case, since $\beta = (\alpha + 2\pi)/\ell_K(A).$ The sharpness result in this theorem is not obtained directly from Robin sharpness by duality, but the construction of the approximating sequences is similar. We have not established sharpness for $\beta = 0$ in the range $0<\Lambda < 2KA$ .

The additional $2\pi$ in \eqref{eq:intro-Steklov-Helmholtz-one} is $R_{K,\lambda}(0+)$, the contribution of the Green region whose area tends to zero. Thus, the model should be regarded as the degenerate split $(0,A)$.

\appendix

\section{Robin spectral convergence under boundary perturbations}
\label{sec:Appendix_Robin_convergence}
\label{sec:Appendix_boundary_oscillation}

Fix a relatively compact smooth region $U$ containing all domains below.
A family $\{\Omega_n\}$ is \emph{uniformly Lipschitz} if its boundaries
admit Lipschitz graph representations in finitely many common coordinate
neighborhoods of fixed size, with uniformly bounded Lipschitz constants.
We write $\Omega_n\to\Omega$ if
\begin{equation}\label{eq:app-domain-convergence}
 \chi_{\Omega_n}\to\chi_\Omega\text{ in }L^1(U),\qquad
 \overline U\setminus\Omega_n\to\overline U\setminus\Omega
 \text{ in Hausdorff distance}.
\end{equation}

\begin{lemma}\label{lem:uniform-extension-trace}
\begin{enumerate}[(i)]
\item If $\Omega_n\Subset U$ are uniformly Lipschitz, there are
linear extension operators $E_n:H^1(\Omega_n)\to H^1(U)$ such that
\begin{equation}\label{eq:app-HRPT-extension}
 (E_nu)|_{\Omega_n}=u,\qquad
 \|E_nu\|_{H^1(U)}\le C\|u\|_{H^1(\Omega_n)},
\end{equation}
with $C$ independent of $n$.
\item If $\mu_n$ are finite signed Borel measures on
$\overline U$ with $|\mu_n|(B(x,r))\le Cr$ for
$x\in\overline U$ and $0<r\le1$, then, for $1/2<s<1$,
\begin{equation}\label{eq:app-uniform-measure-trace}
 \|Tw\|_{L^2(|\mu_n|)}\le C_s\|w\|_{H^s(U)},
 \qquad w\in H^s(U),
\end{equation}
with $C_s$ independent of $n$.
\end{enumerate}
\end{lemma}

\begin{proof}
By \cite[Remark~4.16(2)]{erikssonbiquegong2021}, the Lipschitz
graph patches have common uniform-domain constants.
The Jones extension theorem recalled in
\cite[Section~3, after Definition~3.1]{hinzrozanovateplyaev2023},
with $p=2$, gives uniform $H^1$ extensions on these patches.
Then the partition of unity yields \eqref{eq:app-HRPT-extension}.

For (ii), apply \cite[Corollary~2.2]{hinzrozanovateplyaev2023}
to $|\mu_n|$, with $n=p=2$, $d=1$, and $\beta=s$.
Its condition $(n-d)/p<\beta$ becomes $s>1/2$.
The growth bound and compactness of $\overline U$ bound
$|\mu_n|(\overline U)$ uniformly.
This gives \eqref{eq:app-uniform-measure-trace}.
The conclusion also extends to surfaces by localization in smooth charts.
\end{proof}

All boundary measures are regarded as finite signed Borel measures on
$\overline U$, supported on the respective boundaries.  Thus
$\mu_n\rightharpoonup\mu$ means
$\int\varphi\,d\mu_n\to\int\varphi\,d\mu$ for every
$\varphi\in C(\overline U)$.  For $\mu$ satisfying
$|\mu|(B(x,r))\le Cr$ for $0<r\le1$, define
\begin{equation}\label{eq:app-Robin-measure-definition}
 \lambda_k(\Omega;\mu)
 =\inf_{\substack{V\subset H^1(\Omega)\\\dim V=k}}
   \ \sup_{0\ne u\in V}
 \frac{\displaystyle\int_\Omega|\nabla u|^2\,dA
       +\int_{\partial\Omega}|Tu|^2\,d\mu}
      {\displaystyle\int_\Omega|u|^2\,dA}.
\end{equation}
Here $T$ denotes the Sobolev trace with respect to $|\mu|$;
$\lambda_k(\Omega;\gamma)=\lambda_k(\Omega;\gamma\,ds_{\partial\Omega})$.

\begin{lemma}\label{lem:Robin-signed-uniform-convergence}
Let $\Omega_n,\Omega\Subset U$ be uniformly Lipschitz, with
$\Omega_n\to\Omega$.  Let $\mu_n\rightharpoonup\mu$ be finite signed
boundary measures satisfying, uniformly in $n$,
\begin{equation}\label{eq:app-HRPT-upper-regular}
 |\mu_n|(B(x,r))\le Cr,\qquad x\in\overline U,\quad 0<r\le1.
\end{equation}
Then, for every fixed $k\ge1$,
\begin{equation}\label{eq:app-HRPT-signed-convergence}
 \lambda_k(\Omega_n;\mu_n)\longrightarrow\lambda_k(\Omega;\mu).
\end{equation}
\end{lemma}

\begin{proof}
Lemma~\ref{lem:uniform-extension-trace} gives
\eqref{eq:app-HRPT-extension} and
\eqref{eq:app-uniform-measure-trace} for $1/2<s<1$.
The weak limit $\mu$ has the same total-variation growth bound.
By approximation with smooth functions and the compact embedding
$H^1(U)\hookrightarrow H^s(U)$, we obtain
\begin{equation}\label{eq:app-HRPT-moving-traces}
 w_n\rightharpoonup w\text{ in }H^1(U)
 \quad\Longrightarrow\quad
 \int|Tw_n|^2\,d\mu_n\to\int|Tw|^2\,d\mu.
\end{equation}
For nonnegative measures this is
\cite[Lemma~2.1]{hinzrozanovateplyaev2023}; the same argument
with total variation covers signed measures.

We claim that, for every $\delta>0$,
\begin{equation}\label{eq:app-HRPT-infinitesimal-trace}
 \int_{\partial\Omega_n}|Tu|^2\,d|\mu_n|
 \le\delta\int_{\Omega_n}|\nabla u|^2\,dA
 +C_\delta\int_{\Omega_n}|u|^2\,dA.
\end{equation}
Otherwise, after normalization there are $u_n\in H^1(\Omega_n)$
(with a subsequence understood) such that
\[
 \int|Tu_n|^2\,d|\mu_n|=1,\qquad
 \int_{\Omega_n}|\nabla u_n|^2\le\delta^{-1},\qquad
 \int_{\Omega_n}|u_n|^2\to0.
\]
By \eqref{eq:app-HRPT-extension} and compactness, after passing to a
subsequence,
\[
 E_nu_n\rightharpoonup w\text{ in }H^1(U),\qquad
 E_nu_n\to w\text{ in }L^2(U),\qquad w|_\Omega=0.
\]
Also $|\mu_n|\rightharpoonup\nu$ for a nonnegative measure supported
on $\partial\Omega$.  Its growth bound implies
$\nu\ll ds_{\partial\Omega}$, so $Tw=0$ $\nu$-a.e.  The nonnegative
version of \eqref{eq:app-HRPT-moving-traces} gives
$1=\int|Tw|^2\,d\nu=0$, a contradiction.  Hence
\begin{equation}\label{eq:app-HRPT-Garding}
 \int_{\Omega_n}|\nabla u|^2\,dA
 +\int_{\partial\Omega_n}|Tu|^2\,d\mu_n
 +C\int_{\Omega_n}|u|^2\,dA
 \ge\tfrac12\|u\|_{H^1(\Omega_n)}^2.
\end{equation}
The analogous estimate holds on $\Omega$; both forms are closed,
lower bounded, and have compact resolvent.

Let $\phi_1,\ldots,\phi_k$ be the first $k$ orthonormal
Robin eigenfunctions on $\Omega$, extended to $H^1(U)$.
By \eqref{eq:app-domain-convergence} and
\eqref{eq:app-HRPT-moving-traces}, the $L^2$, Dirichlet, and
boundary products of their restrictions to $\Omega_n$ converge.
Min--max gives
\begin{equation}\label{eq:app-HRPT-limsup}
 \limsup_n\lambda_k(\Omega_n;\mu_n)\le\lambda_k(\Omega;\mu).
\end{equation}
For the reverse inequality, take orthonormal eigenfunctions
$u_{1,n},\ldots,u_{k,n}$ on $\Omega_n$.  Equations
\eqref{eq:app-HRPT-Garding}--\eqref{eq:app-HRPT-limsup} give
uniform $H^1$ bounds, so along a subsequence
\[
 E_nu_{j,n}\rightharpoonup v_j\text{ in }H^1(U),\qquad
 E_nu_{j,n}\to v_j\text{ in }L^2(U).
\]
Using \eqref{eq:app-domain-convergence},
\begin{equation}\label{eq:app-HRPT-limit-orthonormal}
 \int_\Omega v_iv_j\,dA=\delta_{ij}.
\end{equation}
For $u_n=\sum_{j=1}^ka_ju_{j,n}$ and $v=\sum_{j=1}^ka_jv_j$,
weak lower semicontinuity of
$\chi_{\Omega_n}\nabla E_nu_n$ and
\eqref{eq:app-HRPT-moving-traces} imply
\[
 \begin{aligned}
 \int_\Omega|\nabla v|^2\,dA
 &\le\liminf_n\int_{\Omega_n}|\nabla u_n|^2\,dA,\\
 \int|Tv|^2\,d\mu
 &=\lim_n\int|Tu_n|^2\,d\mu_n.
 \end{aligned}
\]
Since
$\mathfrak q_{\Omega_n,\mu_n}[u_n]\le
\lambda_k(\Omega_n;\mu_n)\|u_n\|_{L^2(\Omega_n)}^2$,
\eqref{eq:app-HRPT-limit-orthonormal} and min--max yield
$\lambda_k(\Omega;\mu)\le\liminf_n\lambda_k(\Omega_n;\mu_n)$.
\end{proof}

\begin{lemma}\label{lem:oscillating-Robin-spectrum}
Let $D$ be a smooth disk and $c\ge1$.  There are smooth simply
connected domains $D_n$ with $|D_n|=|D|$ and
$|\partial D_n|\to c|\partial D|$ such that, for every $k\ge1$,
\begin{equation}\label{eq:app-oscillating-spectrum}
 \lambda_k(D_n;\gamma)\longrightarrow\lambda_k(D;c\gamma)
 \quad\text{locally uniformly for }\gamma\in\mathbb R.
\end{equation}
\end{lemma}

\begin{proof}
Let $p=|\partial D|$.  In a fixed normal collar, use
arclength $s\in\mathbb R/p\mathbb Z$ and signed distance $t$, with
$g=dt^2+a(s,t)^2ds^2$, $a(s,0)=1$.  Replace the boundary $t=0$ by
\[
 t=r_n(s):=\frac q{k_n}\sin(k_ns)+b_n,
 \qquad k_n=\frac{2\pi n}{p}.
\]
There is a choice $b_n=O(k_n^{-2})$ for which $|D_n|=|D|$.
The boundary density satisfies
\begin{equation}\label{eq:app-oscillation-density}
 \begin{aligned}
 ds_{\partial D_n}
 &=\sqrt{a(s,r_n)^2+(r_n')^2}\,ds\\
 &=\bigl(\sqrt{1+q^2\cos^2(k_ns)}+o(1)\bigr)\,ds
 \rightharpoonup c(q)\,ds_{\partial D},
 \end{aligned}
\end{equation}
where
\[
 c(q)=\frac1{2\pi}\int_0^{2\pi}
 \sqrt{1+q^2\cos^2\theta}\,d\theta.
\]
Since $c(0)=1$ and $c(q)\to\infty$, choose $q$ such that
$c(q)=c$.  Moreover,
\[
 \|r_n\|_\infty\to0,\qquad
 \sup_n\|r_n'\|_\infty<\infty,\qquad
 D_n\to D.
\]
The fixed collar gives uniform Lipschitz graphs and
\eqref{eq:app-HRPT-upper-regular} for the arclength measures.
For any $\gamma_n\to\gamma$, apply
Lemma~\ref{lem:Robin-signed-uniform-convergence} to the signed measures
$\gamma_n\,ds_{\partial D_n}\rightharpoonup
c\gamma\,ds_{\partial D}$.  Compactness of $I$ and continuity
in $\gamma$ give the locally uniform convergence.
\end{proof}

We call the tangency of two smooth domains $\Omega_1,\Omega_2$
\emph{nondegenerate} if their interiors are disjoint, their boundaries
meet only at $p$, and there are coordinates $(s,t)$ centered at $p$
with $\Omega_1$ locally below $t=h_1(s)$ and $\Omega_2$ above
$t=h_2(s)$, where
\begin{equation}\label{eq:app-nondegenerate-tangency}
 \begin{aligned}
 h_1(0)&=h_2(0)=h'_1(0)=h'_2(0)=0,\\
 h_2(s)-h_1(s)&=as^2+O(|s|^3),\qquad a>0.
 \end{aligned}
\end{equation}

\begin{lemma}\label{lem:short-neck-spectral-convergence}
Let $\Omega_1,\Omega_2$ be smooth bounded domains with a
nondegenerate tangency.  There exist smooth connected domains
$\widehat\Omega_\varepsilon$ satisfying $|\widehat\Omega_{\varepsilon}\triangle (\Omega_1\cup\Omega_2)|\to 0$, $|\partial\widehat\Omega_{\varepsilon}|\to|\bdr\Omega_1|+|\bdr\Omega_2|$, and, for every $k\ge1$,
\begin{equation}\label{eq:converge}
 \lambda_k(\widehat\Omega_\varepsilon;\gamma)
 \longrightarrow\lambda_k(\Omega_1\sqcup\Omega_2;\gamma)
 \quad\text{locally uniformly for }\gamma\in\mathbb R.
\end{equation}
If $\Omega_1,\Omega_2$ are simply connected, so are the
approximating domains.
\end{lemma}

\begin{proof}
Replace $t$ by $t-\frac12(h_1+h_2)$.  In the new coordinates
$h_1=-h$, $h_2=h$, where
\[
 h(s)=\tfrac a2s^2+O(|s|^3)>0\qquad (0<|s|\ll1).
\]
Choose $\chi\in C_c^\infty((-1,1))$, equal to one near zero, and
$L>0$ fixed and sufficiently large.  In $|s|<L\sqrt\varepsilon$,
replace the two boundary graphs by
\begin{equation}\label{eq:app-overlap-graphs}
 t=-h(s)+\varepsilon\chi\bigl(s/(L\sqrt\varepsilon)\bigr),
 \qquad
 t=h(s)-\varepsilon\chi\bigl(s/(L\sqrt\varepsilon)\bigr).
\end{equation}
The boundaries cross twice, at
\[
 s_\varepsilon^\pm=\pm\sqrt{2\varepsilon/a}+O(\varepsilon).
\]
Their union $\Omega_\varepsilon$ is connected.  The segment
$S_\varepsilon=\{t=0,\ s_\varepsilon^-<s<s_\varepsilon^+\}$
lies inside the overlap and divides $\Omega_\varepsilon$ into
$P_{1,\varepsilon}$ and $P_{2,\varepsilon}$.  Their boundaries near
$p$ are the graphs
\begin{equation}\label{eq:app-overlap-cut-graphs}
 \begin{aligned}
 t&=\min\bigl\{0,-h(s)+\varepsilon\chi(s/(L\sqrt\varepsilon))\bigr\},\\
 t&=\max\bigl\{0, h(s)-\varepsilon\chi(s/(L\sqrt\varepsilon))\bigr\}.
 \end{aligned}
\end{equation}
The graph slopes in \eqref{eq:app-overlap-cut-graphs} are bounded
independently of $\varepsilon$.  Moreover,
\[
 |S_\varepsilon|=O(\sqrt\varepsilon),
 \qquad
 |P_{i,\varepsilon}\triangle\Omega_i|=O(\varepsilon^{3/2}).
\]
Away from $p$ the boundaries are unchanged.  Hence the cut
pieces form uniformly Lipschitz families and
$P_{i,\varepsilon}\to\Omega_i$.  Give their exterior boundaries
Robin measure, but give $S_\varepsilon$ zero measure:
\begin{equation}\label{eq:app-overlap-measures}
 \mu_{i,\varepsilon}
 =\restrict{ds}{\partial P_{i,\varepsilon}
                       \cap\partial\Omega_\varepsilon}.
\end{equation}
Then $\mu_{i,\varepsilon}\rightharpoonup ds_{\partial\Omega_i}$. Moreover, $\mu_{i,\varepsilon}$ satisfy \eqref{eq:app-HRPT-upper-regular} uniformly 
since \eqref{eq:app-overlap-cut-graphs} has uniformly bounded slopes. Then for every $\gamma\in\mathbb R$, Lemma \ref{lem:Robin-signed-uniform-convergence}, applied on $P_{i,\varepsilon}$ with $\gamma\mu_{i,\varepsilon}\rightharpoonup \gamma ds_{\partial\Omega_i}$, gives
\[
\lambda_k(P_{i,\varepsilon};\gamma\mu_{i,\varepsilon})\longrightarrow\lambda_k(\Omega_i; \gamma)
\]
Neumann bracketing gives
\(
\lambda_k\bigl((P_{1,\varepsilon};\gamma\mu_{1,\varepsilon})
  \sqcup(P_{2,\varepsilon};\gamma\mu_{2,\varepsilon})\bigr)
 \le\lambda_k(\Omega_\varepsilon;\gamma).
\)
So
\begin{equation}\label{eq:app-neck-Neumann-bracketing}
 \begin{aligned}
 \liminf_{\varepsilon\downarrow0}\lambda_k(\Omega_\varepsilon;\gamma)
 \ge\lambda_k(\Omega_1\sqcup\Omega_2;\gamma).
 \end{aligned}
\end{equation}

For the opposite inequality, let $u_1,\ldots,u_k$ be the first
$k$ orthonormal eigenfunctions on $\Omega_1\sqcup\Omega_2$.  Let
$\eta_\varepsilon=0$ on $B_{C\sqrt\varepsilon}(p)$, equal to one
outside $B_{\varepsilon^{1/4}}(p)$, and logarithmic in between.
Then
\begin{equation}\label{eq:app-neck-log-capacity}
 \int|\nabla\eta_\varepsilon|^2\,dA
 \le\frac{C}{|\log\varepsilon|}\to0,
 \qquad
 \eta_\varepsilon\to1\quad\text{in }L^2(U).
\end{equation}
Multiply each eigenfunction by $\eta_\varepsilon$ and extend by zero. This gives test functions $u_{j,\varepsilon}\in H^1(\Omega_\varepsilon)$ such that
\[
 \begin{aligned}
 \int_{\Omega_\varepsilon}u_{i,\varepsilon}u_{j,\varepsilon}
  &\to\delta_{ij},\\
 \int_{\Omega_\varepsilon}\nabla u_{i,\varepsilon}\cdot
           \nabla u_{j,\varepsilon}
  &\to\sum_{m=1}^2\int_{\Omega_m}\nabla u_i\cdot\nabla u_j,\\
 \int_{\partial\Omega_\varepsilon}u_{i,\varepsilon}u_{j,\varepsilon}\,ds
  &\to\sum_{m=1}^2\int_{\partial\Omega_m}u_i u_j\,ds.
 \end{aligned}
\]
Consequently, min--max gives
\[
 \limsup_{\varepsilon\downarrow0}\lambda_k(\Omega_\varepsilon;\gamma)
 \le\lambda_k(\Omega_1\sqcup\Omega_2;\gamma),
\]
which, with \eqref{eq:app-neck-Neumann-bracketing}, proves convergence for every $\gamma \in \mathbb R$ before smoothing. The convergence is locally uniform since the Robin eigenvalues are strictly increasing and continuous in $\gamma$.

For fixed $\varepsilon$, round the corners in Lipschitz graph
coordinates on a scale $r\downarrow0$.  The resulting domains
$\Omega_{\varepsilon,r}$ converge to $\Omega_\varepsilon$ as
uniformly Lipschitz domains, with Lipschitz constants allowed to depend
on $\varepsilon$.  Their arclength measures converge weakly and
satisfy \eqref{eq:app-HRPT-upper-regular} uniformly in $r$.
Lemma~\ref{lem:Robin-signed-uniform-convergence}, applied to
$\gamma\,ds_{\partial\Omega_{\varepsilon,r}}$, gives that for every fixed $k$, 
\[
\lambda_k(\Omega_{\varepsilon, r};\gamma)\longrightarrow\lambda_k(\Omega_{\varepsilon};\gamma),\qquad (r\downarrow 0)
\]
for every $\gamma \in \mathbb R$. Since $\lambda_k(\Omega_{\varepsilon};\gamma)$ is strictly increasing in $\gamma$, the convergence is locally uniform in $\gamma\in\mathbb R$.  Choose $r=r(\varepsilon)<\varepsilon$ so that
\[
 \sup_{\substack{1\le k\le\varepsilon^{-1}\\
                  |\gamma|\le\varepsilon^{-1}}}
 \left|\lambda_k(\Omega_{\varepsilon,r(\varepsilon)};\gamma)
       -\lambda_k(\Omega_\varepsilon;\gamma)\right|<\varepsilon,
\]
and $\max\{|\Omega_{\varepsilon,r(\varepsilon)}\triangle\,\Omega_{\varepsilon}|, ||\bdr\Omega_{\varepsilon,r(\varepsilon)}|-|\bdr\Omega_{\varepsilon}||\}<\varepsilon$. Then $\widehat\Omega_\varepsilon
:=\Omega_{\varepsilon,r(\varepsilon)}$ satisfies $|\widehat\Omega_{\varepsilon}\triangle(\Omega_1\cup \Omega_2)|\to 0, |\bdr\widehat \Omega_{\varepsilon}|\to |\bdr\Omega_1|+ |\bdr\Omega_2|$ and \eqref{eq:converge},
with simple connectivity preserved. The convergence is locally uniform for the same reason above.
\end{proof}

We also need an area-correction result.
\begin{lemma}\label{lem:area-correction}
Let $\Omega_\varepsilon\Subset U$ be smooth domains with
$|\Omega_\varepsilon|\to M>0$ and
$\lambda_1(\Omega_\varepsilon;-b)\ge-C_b$ for every $b>0$
and all sufficiently small $\varepsilon$. Suppose also that their boundaries contain
arcs of a fixed positive length and uniformly bounded curvature,
with embedded two-sided normal collars of a fixed width.
Then there are smooth domains $\widetilde\Omega_\varepsilon$
diffeomorphic to $\Omega_\varepsilon$ such that \( |\widetilde\Omega_\varepsilon|=M\),  
\(
 |\partial\widetilde\Omega_\varepsilon|
 -|\partial\Omega_\varepsilon|\longrightarrow0,
\)
and, for every fixed $k$,
\[
 \lambda_k(\widetilde\Omega_\varepsilon;\gamma)
 -\lambda_k(\Omega_\varepsilon;\gamma)\longrightarrow0\qquad\text{locally uniformly for }\gamma.
\]
\end{lemma}

\begin{proof}
Taking $b=\delta^{-1}$, by the lower bound of $\lambda_1$ and the min-max principle, for every
$\delta>0$ and $u\in H^1(\Omega_\varepsilon)$,
\[
 \int_{\partial\Omega_\varepsilon}|u|^2\,ds
 \le\delta\int_{\Omega_\varepsilon}|\nabla u|^2\,dA
   +C_\delta\int_{\Omega_\varepsilon}|u|^2\,dA,
\]
uniformly for sufficiently small $\varepsilon$.

By the uniform collar assumptions, choose smooth vector fields
$X_\varepsilon$ supported in the collars, with uniformly bounded
$C^1$ norms, whose flows have area derivatives uniformly bounded
away from zero for small times. Since $|\Omega_\varepsilon|\to M$,
these flows yield smooth diffeomorphisms $J_\varepsilon:\Omega_\varepsilon\to\widetilde\Omega_\varepsilon$
such that $J_\varepsilon=\mathrm{Id}+o(1)$ in $C^1$ and $|\widetilde\Omega_\varepsilon|=M$.
The perimeter changes by $o(1)$.

For $u\in H^1(\Omega_{\varepsilon})$, let 
\[q_\varepsilon[u]:=\|\nabla u\|_{L^2(\Omega_\varepsilon)}^2 + \gamma \|u\|_{L^2(\bdr\Omega_{\varepsilon})}^2,\quad
\widetilde q_\varepsilon[u] := \|\nabla (u\circ J_{\varepsilon}^{-1})\|_{L^2(\widetilde\Omega_\varepsilon)}^2 + \gamma \|u\circ J_{\varepsilon}^{-1}\|_{L^2(\bdr\widetilde\Omega_\varepsilon)}^2.
\]
The trace estimate shows that, for a fixed sufficiently large
$C$, the shifted form $Q_{\varepsilon}[u]:=q_\varepsilon[u]+C\|u\|_{L^2(\Omega_\varepsilon)}^2$
uniformly controls the squared $H^1$ and boundary $L^2$ norms. We also define $\widetilde Q_{\varepsilon}[u] := \widetilde q_{\varepsilon}[u] + C\|u\circ J_{\varepsilon}^{-1}\|_{L^2(\widetilde\Omega_{\varepsilon})}^2$
Change of variables therefore gives $\rho_{\varepsilon} \to 0$ such that
\begin{align*}
 |\widetilde Q_\varepsilon[u]-Q_\varepsilon[u]|
 &\le \rho_{\varepsilon}Q_{\varepsilon}[u],\\
 |\|u\circ J_\varepsilon^{-1}\|_{L^2(\widetilde\Omega_\varepsilon)}^2-\|u\|^2_{L^2(\Omega_{\varepsilon})}|
 &\leq\rho_{\varepsilon}\|u\|_{L^2(\Omega_\varepsilon)}^2,
\end{align*}
uniformly in $u$.
Test functions supported in the inner collars give uniform
upper bounds for each fixed eigenvalue $\lambda_k(\Omega_{\varepsilon};\gamma)$, while the assumption gives lower bounds.  Min--max applied to the $Q_{\varepsilon}$ and $\widetilde Q_{\varepsilon}$ proves the claim, uniformly on bounded parameter intervals.
\end{proof}

\section{The One-Dimensional Robin Variational Problem}
\label{sec:Appendix_Robin_minimizer}

We prove the existence, positivity, and uniqueness results for the Robin
quotient in Subsection~\ref{sec:vairation-def}. Throughout this appendix,
functions in $H^1(\mathbb R)$ are identified with their continuous
representatives. Fix $\alpha>-4\pi$, and let $\mu$ be a nonzero finite
positive Borel measure. 
For $u\in H^1(\mathbb R)$, set
\begin{equation*}
 E_\alpha[u]
 :=\int_{\mathbb R}\bigl(|u'|^2+4\pi^2|u|^2\bigr)\,dt
 +\alpha|u(0)|^2
\end{equation*}
and
\begin{equation}\label{eq:appendix-Robin-quotient}
 \lambda_\alpha(\mu)
 :=\inf_{0\neq u\in H^1(\mathbb R)}
 \frac{E_\alpha[u]}{\displaystyle\int_{\mathbb R}|u|^2\,d\mu},
\end{equation}
where functions with zero denominator are omitted.

The following trace inequality is essential for $E_{\alpha}$ when $\alpha>-4\pi$.

\begin{lemma}\label{lem:appendix-Robin-trace}
For every $u\in H^1(\mathbb R)$,
\begin{align*}
 2\pi|u(0)|^2
 &\leq\int_0^\infty\bigl(|u'|^2+4\pi^2|u|^2\bigr)\,dt,\\
 2\pi|u(0)|^2 
 &\leq\int_{-\infty}^0\bigl(|u'|^2+4\pi^2|u|^2\bigr)\,dt.
\end{align*}
{Equality in the first inequality holds if and only if $u(t)=Ce^{-2\pi t}$ for $t\geq0$, and equality in the second holds if and only if $u(t)=Ce^{2\pi t}$ for $t\leq0$.}
\end{lemma}

\begin{proof}
Since every function in $H^1(0,\infty)$ tends to zero as $t\to\infty$,
we have
\begin{align*}
 \int_0^\infty|u'+2\pi u|^2\,dt
 &=\int_0^\infty\bigl(|u'|^2+4\pi^2|u|^2\bigr)\,dt
 -2\pi|u(0)|^2,\\
 \int_{-\infty}^0|u'-2\pi u|^2\,dt
 &=\int_{-\infty}^0\bigl(|u'|^2+4\pi^2|u|^2\bigr)\,dt
 -2\pi|u(0)|^2.
\end{align*} 
This proves the two inequalities. {The equality cases follow by solving $u'+2\pi u=0$ and $u'-2\pi u=0$ on the respective half-lines.}
\end{proof}

As a consequence, for every $u\in H^1(\mathbb{R})$,
\begin{equation}\label{eq:appendix-Robin-trace-R}
 4\pi|u(0)|^2
 \leq\int_{\mathbb R}\bigl(|u'|^2+4\pi^2|u|^2\bigr)\,dt.
\end{equation}
When $-4\pi<\alpha<0$, this gives
\[
 E_\alpha[u]\geq
 \left(1+\frac{\alpha}{4\pi}\right)\int_{\mathbb R}\bigl(|u'|^2+4\pi^2|u|^2\bigr)\,dt,
\]
while $E_\alpha[u]$ is bounded above by the same bulk energy. If $\alpha\geq0$,
then $E_\alpha[u]$ is bounded below by the bulk energy, whereas
\eqref{eq:appendix-Robin-trace-R} gives
\[
E_\alpha[u]
\leq
\left(1+\frac{\alpha}{4\pi}\right)
\int_{\mathbb R}\bigl(|u'|^2+4\pi^2|u|^2\bigr)\,dt.
\]
Hence for every $\alpha>-4\pi$ there are constants $c_\alpha,C_\alpha>0$ such that
\begin{equation}\label{eq:appendix-Robin-energy-equivalence}
 c_\alpha\|u\|_{H^1(\mathbb R)}^2
 \leq E_\alpha[u]
 \leq C_\alpha\|u\|_{H^1(\mathbb R)}^2.
\end{equation}
 
The map $H^1(\mathbb R)\to L^2(\mu): f\mapsto f$, defined using continuous representatives of $H^1(\mathbb{R})$, is bounded. Indeed, since there exists $C > 0$ such that the sup-norm $\|u\|_{\infty} \leq C\|u\|_{H^1}$ for every $u\in H^1(\mathbb{R})$, we have
\[
\int_{\mathbb{R}}\abs{u}^2d\mu \leq \mu(\mathbb{R})\|u\|_{\infty}^2\leq C^2\mu(\mathbb{R})\|u\|_{H^1}^2.
\] 
We will also use compactness of this map.

\begin{lemma}\label{lem:appendix-Robin-compact-embedding}
For every finite positive Borel measure $\mu$, the  map
\[
 H^1(\mathbb R)\longrightarrow L^2(\mu),\qquad f\longmapsto f,
\]
is compact.
\end{lemma}

\begin{proof}
Let $(u_n)$ be bounded in $H^1(\mathbb R)$. After passing to a subsequence,
$u_n\rightharpoonup u$ weakly in $H^1(\mathbb R)$. The Sobolev embedding gives
a uniform $L^\infty$ bound. Choose compact intervals $K_N\nearrow\mathbb R$.
By the compact embedding $H^1(K_N)\hookrightarrow C^0(K_N)$ and a diagonal
argument, after passing to a further subsequence we have $u_n\to u$ uniformly
on every $K_N$. The same uniform bound passes to $u$.

Fix $\varepsilon>0$. Since $\mu$ is finite, choose $N$ so that
$\mu(\mathbb R\setminus K_N)<\varepsilon$. Then
\[
\begin{aligned}
\|u_n-u\|_{L^2(\mu)}^2 \leq \mu(K_N)\|u_n-u\|_{L^\infty(K_N)}^2 \quad+\mu(\mathbb R\setminus K_N)
\|u_n-u\|_{L^\infty(\mathbb R)}^2.
\end{aligned}
\]
The first term tends to zero, while the second is bounded by a constant times
$\varepsilon$. Hence $u_n\to u$ in $L^2(\mu)$, proving compactness.
\end{proof}

We can now treat the variational problem. In the following proposition, we consider only measures $\mu$ of one of the following forms:
\begin{equation}\label{eq:appendix-Robin-measure-class}
 \mu=M\delta_0, \qquad\text{or}\qquad d\mu=\rho(t)\,dt, \quad 0\leq\rho\in L^1(\mathbb R)\cap C^\infty(\mathbb{R}\setminus\{0\}).
\end{equation}

\begin{proposition}\label{prop:appendix-Robin-minimizer}
Let $\alpha>-4\pi$ and let $\mu$ satisfy
\eqref{eq:appendix-Robin-measure-class}. Then
$\lambda_\alpha(\mu)>0$ and the infimum
\eqref{eq:appendix-Robin-quotient} is attained. Its minimizer is unique up to multiplication by a nonzero constant and may be chosen strictly positive.
\end{proposition}

\begin{proof} 
By boundedness of the map $H^1(\mathbb R)\to L^2(\mu)$ and
\eqref{eq:appendix-Robin-energy-equivalence}, there is $c>0$ such that
\[
 \frac{E_\alpha[u]}{\int |u|^2\,d\mu}\geq c.
\]
Thus $\lambda_\alpha(\mu)>0$. Choose a minimizing sequence $(u_n)$ normalized by
\[
 \int_{\mathbb R}|u_n|^2\,d\mu=1.
\]
Then \eqref{eq:appendix-Robin-energy-equivalence} gives a uniform $H^1$ bound.
After passing to a subsequence,
\[
 u_n\rightharpoonup u\quad\text{in }H^1(\mathbb R),
 \qquad
 u_n\longrightarrow u\quad\text{in }L^2(\mu)
\]
by Lemma~\ref{lem:appendix-Robin-compact-embedding}. Hence
$\int |u|^2\,d\mu=1$. The quadratic form $E_\alpha$ is weakly lower
semicontinuous, so $u$ attains the infimum.

Suppose first that $\mu=M\delta_0$. Lemma~\ref{lem:appendix-Robin-trace} gives
\[
 E_\alpha[u]\geq(4\pi+\alpha)|u(0)|^2.
\]
By the proof of the same lemma, equality holds precisely when 
\[
u'+2\pi u = 0\quad \text{on}\quad (0,\infty),\qquad u'-2\pi u = 0\quad \text{on} \quad (-\infty,0).
\]
The solution $u(t)=u(0)e^{-2\pi|t|}$ is uniquely determined by $u(0)$ and has a fixed sign, which proves all assertions in the point-mass case.

Now let $d\mu=\rho\,dt$. Suppose $u$ is a minimizer with $u\geq0$, and write
$\lambda=\lambda_\alpha(\mu)$. The weak Euler--Lagrange equation gives
\begin{equation*}
-u''+4\pi^2u=\lambda\rho u
\qquad\text{on }\mathbb R\setminus\{0\},
\end{equation*}
or equivalently $u''=(4\pi^2-\lambda\rho)u$. Since
$\rho\in L^1(\mathbb R)$ and $u\in L^\infty(\mathbb R)$, the right-hand side
is locally integrable up to zero, so $u'$ has finite one-sided limits there.
Integrating by parts in the weak Euler--Lagrange identity on $(-\infty,0)$ and $(0,\infty)$ then yields the jump condition
\begin{equation}\label{eq:appendix-Robin-Kirchhoff-potential}
u'(0-)-u'(0+)+\alpha u(0)=0.
\end{equation}

The coefficient $4\pi^2-\lambda\rho$ belongs to $L^1_{\rm loc}$ on each of $(-\infty,0)$ and $(0,\infty)$. Thus the standard existence--uniqueness theorem for the associated
first-order system applies to Cauchy data at an interior point and, by taking
one-sided limits, at the endpoint $0$; see \cite[Lemma~2.1]{wangzettl2013}. If $u(t_0)=0$ for some $t_0>0$, then $u'(t_0)=0$ because $u\geq0$. Uniqueness gives $u\equiv0$ on $(0,\infty)$, so $u(0)=u'(0+)=0$. Equation~\eqref{eq:appendix-Robin-Kirchhoff-potential} then gives $u'(0-)=0$, and uniqueness on $(-\infty,0)$ gives $u\equiv0$, a contradiction. The case $t_0<0$ is the same. If $u(0)=0$, then
\[
 u'(0-)\leq0\leq u'(0+).
\]
Equation~\eqref{eq:appendix-Robin-Kirchhoff-potential} gives
$u'(0-)=u'(0+)=0$, and the same uniqueness argument again gives
$u\equiv0$. Hence
\begin{equation*}
 u(t)>0\qquad(t\in\mathbb R).
\end{equation*}
Thus every nonnegative minimizer is strictly positive. For any minimizer $u$, the identity $E_\alpha[|u|]=E_\alpha[u]$ shows that $|u|$ is also a minimizer. Since $|u|>0$, every minimizer has a fixed sign on $\mathbb R$.

Let $u$ and $v$ be two nonzero minimizers. Consider 
\[
 w=u-\frac{u(t_0)}{v(t_0)}v.
\]
Linearity of the Euler--Lagrange equation shows that $w$ is also a minimizer. Since $w(t_0)=0$, the strict positivity just proved forces $w\equiv0$. Thus $u$ and $v$ are proportional, which proves uniqueness.
\end{proof}

For $\mu_n,\mu$ satisfying the condition of
Proposition~\ref{prop:appendix-Robin-minimizer} and having the same total mass
$A>0$, let $v_n,v$ denote the positive normalized minimizers corresponding to
$\lambda_\alpha(\mu_n),\lambda_\alpha(\mu)$, respectively, with
\[
\int_{\mathbb R}v_n^2\,d\mu_n=\int_{\mathbb R}v^2\,d\mu=1.
\]
The next proposition gives the continuity needed in
Section~\ref{sec:Robin-orthogonality-test}.

\begin{proposition}\label{prop:measure-eigenpair-stability}
Assume $\mu_n\rightharpoonup\mu$, in the sense that
\begin{equation}\label{eq:appendix-measure-weak}
\int_{\mathbb R}\varphi\,d\mu_n
\longrightarrow
\int_{\mathbb R}\varphi\,d\mu
\qquad\text{for every }\varphi\in C_b(\mathbb R).
\end{equation}
Then
\begin{equation}\label{eq:eigen-convergence}
\lambda_\alpha(\mu_n)\longrightarrow\lambda_\alpha(\mu),
\qquad
v_n\longrightarrow v
\quad\text{strongly in }H^1(\mathbb R).
\end{equation}
\end{proposition}

\begin{proof}
Write $\lambda_n:=\lambda_\alpha(\mu_n)$ and
$\lambda:=\lambda_\alpha(\mu)$. Since $v^2\in C_b(\mathbb R)$,
\eqref{eq:appendix-measure-weak} gives
\[
\int_{\mathbb R}v^2\,d\mu_n\longrightarrow1.
\]
Using $v$ as a test function in the variational characterization therefore
yields
\begin{equation}\label{eq:H1-limsup}
\limsup_{n\to\infty}\lambda_n\leq\lambda.
\end{equation}
By the coercivity \eqref{eq:appendix-Robin-energy-equivalence}, the normalized
minimizers $v_n$ are bounded in $H^1(\mathbb R)$. After passing to a
subsequence, write
\[
v_n\rightharpoonup w\qquad\text{weakly in }H^1(\mathbb R).
\]
Since $H^1(K)$ is compactly embedded into $C^0(K)$ for every compact interval $K$, a diagonal argument allows us to pass to a further subsequence such that $v_n\to w$ uniformly on every compact interval.

We claim that
\begin{equation}\label{eq:H1-claim}
\int_{\mathbb R}v_n^2\,d\mu_n
\longrightarrow
\int_{\mathbb R}w^2\,d\mu.
\end{equation}
Indeed, the $H^1$ bound gives a common $L^\infty$ bound for $v_n$ and $w$.
Fix $\varepsilon>0$. Choose $R>0$ and $\chi_R\in C_b(\mathbb R)$ such that
$0\leq\chi_R\leq1$, $\chi_R=0$ on $[-R,R]$, $\chi_R=1$ on
$\mathbb R\setminus[-R-1,R+1]$, and
\[
\int_{\mathbb R}\chi_R\,d\mu<\varepsilon.
\]
By \eqref{eq:appendix-measure-weak}, for all sufficiently large $n$,
\[
\mu_n\bigl(\mathbb R\setminus[-R-1,R+1]\bigr)
\leq\int\chi_R\,d\mu_n<2\varepsilon.
\]
Consequently,
\[
\int_{\mathbb R}|v_n^2-w^2|\,d\mu_n
\leq
A\|v_n^2-w^2\|_{L^\infty([-R-1,R+1])}+C\varepsilon,
\]
where $C$ is independent of $n$. The first term tends to zero by local uniform
convergence. Since $w^2\in C_b(\mathbb R)$,
\[
\int_{\mathbb R}w^2\,d\mu_n\longrightarrow
\int_{\mathbb R}w^2\,d\mu.
\]
Combining the last two displays and then letting $\varepsilon\downarrow0$
proves \eqref{eq:H1-claim}. In particular, the normalization on the left-hand
side gives $\int w^2\,d\mu=1$.

Weak lower semicontinuity of $E_\alpha$ now gives
\[
\lambda
\leq E_\alpha[w]
\leq\liminf_{n\to\infty}E_\alpha[v_n]
=\liminf_{n\to\infty}\lambda_n.
\]
Together with \eqref{eq:H1-limsup}, this yields $\lambda_n\to\lambda$, and
$w$ is a normalized minimizer for $\mu$. By uniqueness and positivity,
$w=v$. Since the same argument applies to every subsequence, the whole
sequence satisfies $v_n\rightharpoonup v$ in $H^1(\mathbb R)$.

Let $E_\alpha[\cdot,\cdot]$ be the bilinear form associated with $E_\alpha$.
Coercivity gives
\begin{align*}
c_\alpha\|v_n-v\|_{H^1}^2 \leq E_\alpha[v_n-v] =E_\alpha[v_n]+E_\alpha[v]-2E_\alpha[v_n,v].
\end{align*}
Here $E_\alpha[v_n]=\lambda_n\to\lambda=E_\alpha[v]$, while weak
$H^1$ convergence implies $E_\alpha[v_n,v]\to E_\alpha[v,v]$. Hence the
right-hand side tends to zero, proving the strong convergence in
\eqref{eq:eigen-convergence}.
\end{proof}

\section{Endpoint Behavior of the Riccati Variable}\label{sec:Appendix_Riccati}

We record the endpoint facts used in the Riccati comparison of Section~\ref{sec:kappa}.

\begin{lemma}\label{lem:Robin-singular-flux}
Let $m>0$, let $G\in\mathcal G_m$, let $\lambda>0$, and let
$0<f\in\mathcal F(G)$ solve
\[
-(Gf')'+\frac{4\pi^2}{G}f=\lambda f
\qquad\text{in }(0,m).
\]
Then
\[
f(x)\longrightarrow 0,\qquad
G(x)f'(x)\longrightarrow0
\qquad (x\downarrow0).
\]
Moreover, for all sufficiently small $x>0$,
\[
2\pi-\lambda x
\leq
\frac{G(x)f'(x)}{f(x)}
<2\pi.
\]
\end{lemma}

\begin{proof}
Since $f\in\mathcal F(G)$,
\[
\int_0^m |ff'|\,dx
\leq
\left(\int_0^m \frac{f^2}{G}\,dx\right)^{1/2}
\left(\int_0^m G(f')^2\,dx\right)^{1/2}
<\infty.
\]
Thus $f^2$ has a finite limit at $0$. Since
\[
\int_0^m \frac{dx}{G(x)}=\infty
\qquad\text{and}\qquad
\frac{f}{\sqrt G}\in L^2(0,m),
\]
this limit must be zero. Hence
\[
f(x)\longrightarrow0
\qquad (x\downarrow0).
\]

For all sufficiently small $x>0$,
\[
(Gf')'
=
\left(\frac{4\pi^2}{G}-\lambda\right)f
>0.
\]
Since $f>0$ and $f(0+)=0$, it follows that $Gf'>0$ near $0$. Indeed,
if $Gf'(x_0)\leq0$ at some sufficiently small $x_0$, then
$Gf'\leq0$ on $(0,x_0]$, and hence $f'\leq0$ there, contradicting
$f(0+)=0<f(x_0)$. Thus $Gf'$ has a nonnegative limit as $x\downarrow0$.
This limit must be zero, for otherwise $Gf'\geq c>0$ near $0$, and
\[
f(x_0)-f(x)
=
\int_x^{x_0}\frac{G(s)f'(s)}{G(s)}\,ds
\geq
c\int_x^{x_0}\frac{ds}{G(s)}
\longrightarrow\infty,
\]
a contradiction. Therefore
\[
G(x)f'(x)\longrightarrow0
\qquad (x\downarrow0).
\]

Fix $x_0>0$ sufficiently small and define the integrating factor
\[
\rho(x)
:=
\exp\left(
-\int_x^{x_0}\frac{2\pi}{G(s)}\,ds
\right),
\qquad 0<x\leq x_0.
\]
Then $G\rho' = 2\pi \rho$ with $\rho(0+)=0$ and $\rho'>0$. A direct calculation using the
differential equation gives
\[
\bigl[\rho(2\pi f-Gf')\bigr]'
=
\lambda\rho f.
\]
Since
\[
f(0+)=0,
\qquad
Gf'(0+)=0,
\qquad
\rho(0+)=0,
\]
integration from $0$ to $x$ yields
\[
\rho(x)\bigl(2\pi f(x)-G(x)f'(x)\bigr)
=
\lambda\int_0^x \rho(s)f(s)\,ds.
\]
Since $Gf'>0$ near $0$, the function $f$ is increasing there. Together
with the monotonicity of $\rho$, this gives
\[
0
<
\lambda\int_0^x \rho(s)f(s)\,ds
\leq
\lambda x\,\rho(x)f(x).
\]
Therefore
\[
0
<
2\pi f(x)-G(x)f'(x)
\leq
\lambda x f(x).
\]
Dividing by $f(x)>0$ gives
\[
0
<
2\pi-\frac{G(x)f'(x)}{f(x)}
\leq
\lambda x,
\]
and the result follows.
\end{proof}

In particular,
\[
\lim_{x\downarrow0}\frac{G(x)f'(x)}{f(x)}=2\pi.
\]

\begin{lemma}\label{lem:Robin-regular-end-flux}
Let $G\in\mathcal G_m$ and $0<f\in\mathcal F(G)$ solve
\[
-(Gf')'+\frac{4\pi^2}{G}f=\lambda f
\qquad\text{in }(0,m),
\]
with $f(m)>0$. Then $R:=Gf'/f$ extends continuously to $m$.
If $G^{-1}\in C(0,m]$, then $R\in C^1(0,m]$.
\end{lemma}

\begin{proof}
For $0<x<y<m$, Cauchy--Schwarz gives
\[
|f(y)-f(x)|^2
\leq\int_x^y G|f'|^2\,dt
\int_x^y\frac{dt}{G(t)}.
\]
Thus $f(m)$ exists. Since $f(m)>0$, $f$ is bounded near $m$ and
\[
(Gf')'=\left(\frac{4\pi^2}{G}-\lambda\right)f
\in L^1(x_0,m),\qquad x_0>0.
\]
Consequently, $Gf'$ has a finite limit at $m$, and so does $R$.
Finally,
\[
R'=\frac{4\pi^2-R^2}{G}-\lambda.
\]
If $G^{-1}\in C(0,m]$, the right-hand side extends continuously to $m$.
\end{proof}

\section{The Angular--Radial Crossing on Spherical Disks}
\label{sec:Appendix-mode-crossing}

\begin{proposition}\label{prop:Robin-angular-radial-crossing}
Let $K>0$, $0<KA<4\pi$, and $-4\pi<\alpha<-4\pi+KA$. Set
\[
\beta:=\frac{\alpha+2\pi}{\ell_K(A)},\qquad
\mu:=-\frac{(\alpha+2\pi)(\alpha+KA)}{A(4\pi-KA)}.
\]
Then
\[
\lambda_{1,\mathrm{ang}}(D_K(A);\beta)
=\lambda_{2,\mathrm{rad}}(D_K(A);\beta)
\]
if and only if
\begin{equation}\label{eq:appendix-angular-radial-condition}
-KA<\alpha<-2\pi,\qquad 0<\mu<2K,\qquad
R_{K,\mu}(A)+\alpha+2\pi=0.
\end{equation}
In particular, the equality occurs only when $KA>2\pi$.
\end{proposition}

\begin{proof} 
Let $r_0$ be the radius of $D_K(A)$. Then
\[
\ell_K(A)=2\pi s_K(r_0),\qquad
\frac{s_K'(r_0)}{s_K(r_0)}
=\frac{2\pi-KA}{\ell_K(A)},\qquad
\ell_K(A)^2=A(4\pi-KA).
\]
Suppose first that the two eigenvalues coincide, and denote their
common value by $\lambda$. Let $v$ be a second radial eigenfunction,
normalized by $v(0)>0$. It satisfies
\[
v''+\frac{s_K'}{s_K}v'+\lambda v=0,\qquad
v'(0)=0,\qquad v'(r_0)+\beta v(r_0)=0.
\]
Since $v$ has one interior zero, $(s_Kv')'=-\lambda s_Kv$ implies
$\lambda>0$. Differentiating and using $(s_K'/s_K)'=-s_K^{-2}$,
we find that $w:=-v'$ solves
\begin{equation}\label{eq:appendix-angular-from-radial}
w''+\frac{s_K'}{s_K}w'+\left(\lambda-\frac{1}{s_K^2}\right)w=0.
\end{equation}
Moreover, $w(r)=\lambda v(0)r/2+O(r^3)$ at the center.
The regular solutions of this angular equation form a one-dimensional
space. Thus $w$ is a positive multiple of the radial factor of the first
angular eigenfunction. In particular, $w>0$ on $(0,r_0]$ and
$w'(r_0)+\beta w(r_0)=0$.

Consequently, $v'=-w<0$. Since $v$ has exactly one zero,
$v(r_0)<0$ and $\beta=-v'(r_0)/v(r_0)<0$. The two Robin conditions give
\[
0=w'(r_0)+\beta w(r_0)
=\left(\lambda+\beta^2
-\beta\frac{s_K'(r_0)}{s_K(r_0)}\right)v(r_0),
\]
hence
\begin{equation}\label{eq:appendix-angular-radial-mu}
\lambda=\beta\left(\frac{s_K'(r_0)}{s_K(r_0)}-\beta\right)
=-\frac{(\alpha+2\pi)(\alpha+KA)}{A(4\pi-KA)}=\mu.
\end{equation}
Together with $\lambda>0$ and $\beta<0$, this yields
$-KA<\alpha<-2\pi$. By \eqref{eq:model-angular-boundary-early},
\[
R_{K,\lambda}(A)+\alpha+2\pi=0.
\]
Since $\alpha<-4\pi+KA$, the value
$R_{K,\lambda}(A)=-\alpha-2\pi$ exceeds
$R_{K,2K}(A)=2\pi-KA$. Strict monotonicity in the spectral parameter
(Lemma~\ref{lem:model-R-mu-monotonicity}) gives $\lambda<2K$.
Thus \eqref{eq:appendix-angular-radial-condition} holds.

Conversely, assume \eqref{eq:appendix-angular-radial-condition}.
Let $w>0$ be the regular first-angular solution of
\eqref{eq:appendix-angular-from-radial} with $\lambda=\mu$.
The Riccati condition gives $w'(r_0)+\beta w(r_0)=0$.
Define
\[
v:=\frac{w'+(s_K'/s_K)w}{\mu}.
\]
The angular equation gives $v'=-w$, so $v$ is the regular radial
solution with eigenvalue $\mu$ and $v(0)=2w'(0)/\mu>0$.
Using the definition of $\mu$,
\[
v(r_0)=\frac{(s_K'(r_0)/s_K(r_0)-\beta)w(r_0)}{\mu}
=\frac{w(r_0)}\beta<0,\qquad
v'(r_0)+\beta v(r_0)=0.
\]
Thus $v$ is strictly decreasing and has exactly one interior zero,
so $\mu$ is the second radial Robin eigenvalue, as required.
\end{proof}

\section*{Declaration on the Use of Artificial Intelligence}

During the preparation of this manuscript, the authors used ChatGPT for language polishing and grammar checking. 
ChatGPT was also used to identify the reference \cite{arrietabruschi2006}, and to assist with the proofs of Proposition B.4, Lemma C.1 and Proposition D.1. 
The authors independently verified all computations, reviewed and edited all AI-assisted content, and take full responsibility for the content of the manuscript.

\end{document}